\documentclass[11pt,a4paper,leqno]{amsart}
\usepackage{amssymb,amsmath,amsthm}
\usepackage[utf8]{inputenc}
\usepackage[english]{babel}
\usepackage{multicol}
\usepackage{graphicx}
\usepackage[hidelinks]{hyperref} 
\usepackage{xcolor}
\hypersetup{
 colorlinks,
 linkcolor={red!50!black},
 citecolor={blue!50!black},
 urlcolor={blue!80!black}
}
\usepackage{enumerate} % Fancy enumeration lists 
\usepackage{dsfont}
\usepackage[margin=2.0cm]{geometry}

\usepackage{tikz}
\usetikzlibrary{patterns}
\usetikzlibrary[patterns.meta]
\usetikzlibrary{arrows.meta,backgrounds}
\usepackage{subfig}

\usepackage{relsize} 
\usepackage{bigints}
\allowdisplaybreaks
\numberwithin{equation}{section}

\theoremstyle{plain}
\newtheorem{theorem}{Theorem}[section]

\newtheorem{lemma}[theorem]{Lemma}
\newtheorem{corollary}[theorem]{Corollary}
\newtheorem{proposition}[theorem]{Proposition}
\newtheorem{remark}[theorem]{Remark}

\newtheorem{conj}{Conjecture}

\def\ds{\displaystyle}

\def\R{{\mathbb R}}% real numbers
 \def\dd{\mathrm d}
\newcommand{\loc}{\rm loc}

\DeclareMathOperator{\curl}{curl}
\DeclareMathOperator{\diam}{diam}
\newcommand{\dist}{\operatorname{dist}}
\renewcommand{\div}{\operatorname{div}}
\newcommand{\supp}{\operatorname{supp}}

\renewcommand{\leq}{\leqslant}\renewcommand{\le}{\leqslant}
\renewcommand{\geq}{\geqslant}\renewcommand{\ge}{\geqslant}

\def\cal#1{\mathcal{#1}}

\newcommand{\Ind}{\mathds{1}}

\def\eps{\varepsilon}

\renewcommand{\epsilon}{\varepsilon}

\renewcommand{\tilde}{\widetilde}
\newcommand{\der}[2]{\frac{\dd #1}{\dd #2}}

\date\today

\title[]{On the dynamics of leapfrogging vortex rings}

\author[M. Donati]{Martin Donati}

\author[L.E. Hientzsch]{Lars Eric Hientzsch}

\author[C. Lacave]{Christophe Lacave}

\author[E. Miot]{Evelyne Miot}

\email{\newline martin.donati@univ-grenoble-alpes.fr \newline lars.hientzsch@kit.edu \newline Christophe.Lacave@univ-smb.fr \newline evelyne.miot@univ-grenoble-alpes.fr}

\address[M. Donati $\&$ E. Miot]{Univ. Grenoble Alpes, CNRS, Institut Fourier, F-38000 Grenoble, France.}

\address[L.E. Hientzsch]{Institute for Analysis, Karlsruhe Institute of Technology (KIT), D-76128 Karlsruhe, Germany
}

\address[C. Lacave]{Univ. Savoie Mont Blanc, CNRS, LAMA, ISTerre, 73000 Chamb\'ery, France}

\begin{document}
\begin{abstract} 
The evolution of highly concentrated vorticity around rings in the three-dimensional axisymmetric Euler equations is studied in a regime for which the leapfrogging dynamics predicted by Helmholtz is expected to occur. We provide in this paper the first result deriving this phenomenon for a general class of initial data in the suitable regime. The singular interaction of rings requires significant improvements of weak and strong localization estimates obtained in prior works. Our method is based on the combination of a new variational argument and a recently introduced double iterative procedure.
\end{abstract}

\maketitle

\section{Introduction}
The goal of this paper is to rigorously study the evolution of interacting vortex rings in $3D$ incompressible ideal fluids.
Vortex rings are experimentally and numerically well observed coherent structures propagating in an ideal fluid. Mathematically, they correspond to axisymmetric solutions to the 3D incompressible Euler equations concentrated around annuli. Due to its self-induced motion, linked to \emph{local induction approximation} (LIA) \cite{daRios}, a single vortex ring translates at leading order along its symmetry axis. Specifically, it was computed by Kelvin \cite{Kelvin, Lamb} that a vortex ring of major radius $r$, thickness $\eps$, and circulation $\gamma$ has a self-induced velocity at leading order equal to
\begin{equation}\label{eq:ring_v}
 v= \frac{\gamma}{4\pi r} |\ln \eps| e_z.
\end{equation}

It is well-known since the work of Helmholtz \cite{Helmholtz1} in 1858, see also \cite{Dyson, Dyson2, Hicks, Lamb, Borisov}, that two coaxial vortex rings can successively pass through each other in a periodic motion. More precisely and according to \cite{Helmholtz1}, two given vortex rings translating along the same axis may exhibit the so-called \emph{leapfrogging motion} in which `\emph{the foremost widens and travels more slowly, the pursuer shrinks and travels faster till finally, if their velocities are not too different, it overtakes the first and penetrates it. Then the same game goes on in the opposite order, so that the rings pass through each other alternately'}, see Figure~\ref{fig:leapfrog}. 

\begin{figure}[ht]
 \centering
\begin{tikzpicture}[
 scale=1.4,
 >={Stealth[scale=1.2]},
 z axis/.style={thick,->,shorten >=2pt},
 red tube/.style={red!80, line width=2.5pt},
 blue tube/.style={blue!80, line width=2.5pt},
 red reflection/.style={red!40, line width=1pt},
 blue reflection/.style={blue!40, line width=1pt}
]

\begin{scope}[xshift=0cm]
 \draw[z axis] (0,0) -- (2,0) node[right] {$z$};
 
\draw[blue tube] (0.3,0) ++(5:0.4 and 0.96) arc (5:355:0.4 and 0.96);
\draw[blue reflection] (0.3,0.96) arc (90:270:0.4 and 0.96);
 
\draw[red tube] (1.3,0) ++(5:0.4 and 0.96) arc (5:355:0.4 and 0.96); \draw[red reflection] (1.3,0.96) arc (90:270:0.4 and 0.96);
 
 \draw[->, blue, thick] (0.3,0.96) -- (0.7,0.56);
 
 \draw[->, red, thick] (1.3,0.96) -- (1.7,1.36);
\end{scope}

\begin{scope}[xshift=3cm] 
 \draw[z axis] (0,0) -- (2,0) node[right] {$z$};
 
\draw[blue tube] (0.9,0) ++(5:0.3 and 0.72) arc (5:355:0.3 and 0.72);
\draw[blue reflection] (0.9,0.72) arc (90:270:0.3 and 0.72);
 
\draw[red tube] (0.9,0) ++(4:0.5 and 1.2) arc (4:356:0.5 and 1.2); 
\draw[red reflection] (0.9,1.2) arc (90:270:0.5 and 1.2);
 
 \draw[->, red, thick] (0.9,1.2) -- (1.3,1.2);
 \draw[->, blue, thick] (0.9,0.72) -- (1.7,0.72);
\end{scope}

\begin{scope}[xshift=6cm] 
 \draw[z axis] (0,0) -- (2,0) node[right] {$z$};
 
\draw[red tube] (0.3,0) ++(5:0.4 and 0.96) arc (5:355:0.4 and 0.96);
\draw[red reflection] (0.3,0.96) arc (90:270:0.4 and 0.96);
 
\draw[blue tube] (1.3,0) ++(5:0.4 and 0.96) arc (5:355:0.4 and 0.96); \draw[blue reflection] (1.3,0.96) arc (90:270:0.4 and 0.96);
 
 \draw[->, red, thick] (0.3,0.96) -- (0.7,0.56);
 
 \draw[->, blue, thick] (1.3,0.96) -- (1.7,1.36);
\end{scope}

\begin{scope}[xshift=9cm] 
 \draw[z axis] (0,0) -- (2,0) node[right] {$z$};
 
\draw[red tube] (0.9,0) ++(5:0.3 and 0.72) arc (5:355:0.3 and 0.72);
\draw[red reflection] (0.9,0.72) arc (90:270:0.3 and 0.72);
 
\draw[blue tube] (0.9,0) ++(4:0.5 and 1.2) arc (4:356:0.5 and 1.2); 
\draw[blue reflection] (0.9,1.2) arc (90:270:0.5 and 1.2);
 
 \draw[->, blue, thick] (0.9,1.2) -- (1.3,1.2);
 \draw[->, red, thick] (0.9,0.72) -- (1.7,0.72);
\end{scope}
\end{tikzpicture}
 \caption{Leapfrogging phenomenon}\label{fig:leapfrog}
\end{figure}

For such phenomena to occur, suitable assumptions on the initial locations, radii $r$ and thickness $\eps$ and intensity $\gamma$ of the coaxial vortex rings need to be made. Specifically, coaxial vortex rings need to be located at very small mutual distance leading to sufficiently strong interactions inducing a change of their radii at a suitable scale. This in turn causes a relevant variation of their self-induced velocity~\eqref{eq:ring_v}. Here, we say that leapfrogging dynamics occurs when the variation of the self-induced motion \eqref{eq:ring_v} is of the same order as the mutual interaction of the vortex rings, namely $\sqrt{|\ln\eps|}$. We refer the reader to Appendix~\ref{appendix:regimes} for a detailed discussion of the suitable scaling regime and motion law. 

Up to a few crossings, such motion can be experimentally \cite{Yamada} and numerically \cite{AlvarezNing, lim, RileyStevens} observed, supporting this conjecture for finite times.

The suitable framework to set up the mathematical problem is given by the axisymmetric Euler equations without swirl in vorticity form in cylindrical coordinates $(z,r,\theta)$. After dimension reduction due to the axial symmetry and for the potential vorticity $\omega_\eps$, see Section~\ref{sec:WP} below for details, we are led to the system
\begin{equation}\label{eq:Axi2D_rescaled}
\begin{cases}
 \ds \partial_t \, \omega^\eps + \frac{1}{|\ln \eps|} u^\eps\cdot \nabla \omega^\eps = 0 & \text{ in } \R_+^* \times \R^2_+, \vspace{1mm}\\
 \div (ru^\eps) =0 & \text{ in } \R_+ \times \R^2_+ ,\vspace{1mm} \\
 \curl (u^\eps) = r\omega^\eps & \text{ in } \R_+ \times \R^2_+ , \vspace{1mm} \\
 \omega^\eps(0,\cdot) = \omega_0^\eps & \text{ in } \R^2_+,
\end{cases}
\end{equation}
where $\R_+^2$ denotes the space $\{(z,r) \in \R \times \R_+\}$.

In the present setting, all vortex rings obey at leading order the same self-induced motion \eqref{eq:ring_v}. This suggests to scale the time variable by a factor $1/|\ln \eps|$ yielding \eqref{eq:Axi2D_rescaled}. The mutual distance of vortex rings is then chosen to be of order $1/\sqrt{|\ln \eps|}$ leading to a \emph{relative motion} of velocity $v' = \mathcal{O}\left(\sqrt{|\ln \eps|}\right)$. Indeed, in a time of order $1$ for this system, the vortex rings are expected to exhibit a complete rotation around each other while traveling a vertical distance of order $1$, see also Appendix~\ref{appendix:regimes}.

In a slightly informal way, the leapfrogging conjecture can be formulated as follows. 
\begin{conj}\label{conj:LeapFrog}
Let $N>1$ and $0<\eps\ll 1$. Given initial vorticity $\omega_0^\eps$ decomposed in $N$ blobs, each of circulation $\gamma>0$ and sharply `concentrated' around $N$ points $X_{i,0}^\eps$ (in terms of $\eps$) located at mutual distance of order $1/\sqrt{|\ln \eps|}$, the unique solution $\omega^\eps(t,\cdot)$ to \eqref{eq:Axi2D_rescaled} remains `concentrated' around $N$ points $X_i^\eps(t)$.\\
Moreover, the trajectories of these points are governed by the system
\begin{equation}\label{def:X_i}
 \begin{cases}
 \displaystyle \der{}{t} X_i^\eps(t) = \frac{\gamma}{4\pi X_{i,r}^\eps(t)}e_z + \frac{\gamma}{2\pi|\ln \eps|} \sum_{j \neq i} \frac{\big(X_i^\eps(t)-X_j^\eps(t)\big)^\perp}{\big|X_i^\eps(t)-X_j^\eps(t)\big|^2},\\
 X_i^\eps(0) = X_{i,0}^\eps.
 \end{cases}
\end{equation}
\end{conj}
 The meaning of `concentration' will be precised below.
Concerning \eqref{def:X_i}, the motion of the vortex centers is driven by the vertical self-induced velocity corresponding to the first term and the second term that represents the well-known mutual interaction in the 2D point-vortex dynamics. Due to the chosen fast time scale, both are reported up to a factor $1/|\ln \eps|$. The choice of placing the vortex centers at mutual distance of order ${1}/{\sqrt{|\ln\eps|}}$ yields that the relative self-induced motion, i.e. the second order term in \eqref{eq:ring_v}, is of the same order than the motion due to the interaction, see Section~\ref{sec:LP}. Note that the velocity induced by the mutual interactions is inverse proportional to the inter-vortex distance. In the case $N=2$, System~\eqref{def:X_i} describes the evolution of two leapfrogging vortex rings.

The key difficulties in attempting to prove this conjecture are twofold. The present regime is critical in the sense that the aforementioned effects are of the same strength and the vortex centers are situated very close to each other. This makes the problem at hand conceivably harder than the ones in which the motion is dominated by either of these effects. Commonly, showing the persistence of the localization property on suitable time-scales turns out to be the major difficulty. Subsequently, the asymptotic motion law for small $\eps>0$ needs to be derived. Experimental and numerical evidence supports the conjecture up to a few number of crossings. However, in these experiments various instability mechanisms are displayed and the periodic motion does not occur for arbitrarily large times \cite{lim, Maxworthy, RileyStevens, WidnallSullivan, Yamada}. 
The \emph{vortex filamentation phenomenon} appears to be particularly relevant, it causes a vortex ring to leave a small tail of vorticity behind which may potentially collide with following vortex rings, see for this mechanism in various scenarios \cite{Childress, ChildressGilbertValiant, Dritschel, Moffatt, PullinMoore}.

The main result of this paper validates the above conjecture: we prove that an arbitrary but finite number $N$ of initially sharply concentrated vortex rings remains concentrated for a small but positive time independent of $\eps>0$. The initial concentration assumptions allow for general data to be considered, and in particular, do not require any specific vorticity profile.

\subsection{Main result}
We introduce the precise mathematical setup for the problem and state our main result. For $x \in \R^2_+$, we denote by $x = (x_z,x_r)$ with $x_z \in \R$ and $x_r \in \R_+$, and by $(e_z, e_r)$ the canonical basis on $\R^2$.
We consider a fixed reference point 
\begin{equation*}
 X^*=(z^\ast, r^\ast).
\end{equation*}
For $\eps>0$, the initial data $\omega_0^\eps$ will be given by a superposition of $N$ vorticities of the same circulation $\gamma$ and sharply concentrated around points located at distance of order $1/\sqrt{|\ln\eps|}$ from the reference point
\begin{equation}\label{distance-regime}
X^\eps_{i,0}:= X^* + \frac{Y_{i,0}}{\sqrt{|\ln\eps|}},
\end{equation}
where $Y_{i,0}\in \R^2$ is fixed independently of $\eps$. Specifically, let $\gamma \neq 0$ and assume that
$\omega_0^\eps \in L^\infty (\R^2_+)$ such that

\begin{equation}\label{hyp:omega0}
 \left\{\begin{aligned}
& \omega_0^\eps = \sum_{i=1}^N \omega_{i,0}^\eps, \\
& \supp \omega_{i,0}^\eps \subset B\big( X_{i,0}^\eps,\eps \big), \\
& \omega_{i,0}^\eps \text{ has a definite sign and } \int_{\R^2_+} x_r \omega_{i,0}^\eps(x) \dd x = \gamma, \\
& |\omega_0^\eps| \le \frac{M_0}{\eps^2},
\end{aligned}\right.
\end{equation} 
for some $M_0>0$ independent of $\eps$, see Figure~\ref{fig:initial_data} for an example.
\begin{figure}
    \centering
    \begin{tikzpicture}
    % 1. Draw the axes
    \draw[thick, ->] (-1,0) -- (5.5,0) node[right] {$z$};
    \draw[thick, ->] (0,-0.5) -- (0,4.5) node[above] {$r$};
    
    % 2. Define Coordinates
    \coordinate (Xast) at (2.5,2.5);
    \coordinate (X1) at (1.5,1.5);
    \coordinate (X2) at (3,3);

    % 4. Draw the squared box with black boundary centered at X*
    % Center (2.5, 2.5), side length 1.0
    \draw[ dashed, thick] (0.9,0.9) rectangle (4,4);
    
    % 5. Add the points and labels
    % X* is blue
    \fill[blue] (Xast) circle (2pt) node[right, blue] {$X_{\ast}$};
    
    % X_{0,1} with blue boundary circle (radius 0.25)
    \fill (X1) circle (1.5pt) node[below right] {$X_{0,1}^{\varepsilon}$};
    \draw[red, thick] (X1) circle (0.25);
% Radius indicator for the first circle
    \draw[red, ->] (X1) -- ++(135:0.25) node[above right, scale=0.8] {$\mathcal{O}(\varepsilon)$};
    
    % X_{0,2} with blue boundary circle (radius 0.25)
    \fill (X2) circle (1.5pt) node[above right] {$X_{0,2}^{\varepsilon}$};
    \draw[red, thick] (X2) circle (0.25);
    
    % 6. Side length indicator as a "superscript" of the upper side
    \draw[<->, dotted, thick] (1, 4.2) -- (4, 4.2) 
        node[midway, above] {$\mathcal{O}(|\ln\varepsilon|^{-\frac{1}{2}})$};
    
    % 7. Origin label
    \node[below left] at (0,0) {$0$};

    % Support
    \draw[red] (3.1,2.95) -- (4.3,2.5) node[right]{$\supp \omega_{2,0}^\eps$};

\end{tikzpicture}
    \caption{Example of initial data for two vortex rings }
    \label{fig:initial_data}
\end{figure}
Given that the potential vorticity $\omega^\eps$ satisfies the transport equation \eqref{eq:Axi2D_rescaled} and the existence of a unique associated flow, see Remark~\ref{rmk: decomposition vorticity}, it also holds $\omega^\eps(t)=\sum_{i=1}^N\omega_{i}^\eps(t)$ for $t\geq 0$. Further, it follows from the transport equation that $\omega^\eps(t)$ conserves the sign definiteness as well as the $L^\infty$-norm. By the divergence free condition of $ru^\eps$, \eqref{eq:Axi2D_rescaled} also implies the conservation of the mass of each blob:
 \[
 \int_{\R^2_+} x_r \omega_{i}^\eps(t,x) \dd x = \gamma,\quad \forall t\geq 0,
 \]
see Section~\ref{sec-WP} and \eqref{eq:consgamma} for more details.

We are now in position to precisely formulate our main result.
\begin{theorem}\label{theo:main} 
Let $X^*= (z^*, r^*) \in \R^2_+$, $(Y_{i,0})_{1 \le i \le N}$ be $N$ points in $\R^2$ such that $Y_{i,0,r} \neq Y_{j,0,r}$ for all $i\neq j$. Let $\gamma \neq 0$, $M_0>0$, and $\omega_0^\eps \in L^\infty (\R^2_+)$ satisfying \eqref{hyp:omega0} for every $\eps > 0$. 

There exists $T_0>0$, $\eps_0>0$ and $C_0>0$, depending only on $X^*$, $(Y_{i,0})_{1 \le i \le N}$ and $M_0$, such that for every $\eps\in (0,\eps_0]$, there exists a unique solution $(t\mapsto X_i^\eps)_{1\le i \le N}$ of \eqref{def:X_i} on $[0,T_0]$, and the unique weak solution $(u^\eps,\omega^\eps)$ of \eqref{eq:Axi2D_rescaled} on $[0,T_0]$ (in the sense of Proposition ~\ref{prop:WP} below) can be decomposed as 
\[
\omega^\eps=\sum_{i=1}^N \omega_i^\eps, \quad \omega_i^\eps \in L^\infty(\R_+\times \R^2_+)
\]
with the following properties.
\begin{enumerate}[(i)]
 \item A weak localization property:
 \begin{equation*}
 \sup_{t \in [0,T_0]}\left| \gamma - \int_{B\left( X_i^\eps(t) , C_0\frac{\ln |\ln \eps|}{|\ln \eps|}\right)} x_r\omega_i^\eps(t,x)\dd x\right| \le \gamma\frac{\ln |\ln \eps|}{|\ln \eps|}.
 \end{equation*}
 
 \item A strong localization property: 
 \begin{equation*}
 \supp \omega_i^\eps(t,\cdot) \subset \left\{ (z,r) \in \R^2_+ \, , \, \big| r - X_{i,0,r}^\eps\big| \le \frac14\frac{\min_{j\neq k} |Y_{j,0,r} - Y_{k,0,r}|}{\sqrt{|\ln\eps|}} \, , \, \big| z - z^*\big| \le 1 \right\}, \text{ for all}\ t\in [0,T_0].
 \end{equation*}
\end{enumerate}
\end{theorem}

We emphasize that the conclusions of Theorem~\ref{theo:main} hold for \emph{any} initial data $\omega_0^\eps \in L^\infty(\R^2_+)$ satisfying assumptions~\eqref{hyp:omega0}, which ranges from smooth compactly supported blobs to vortex patches, even possibly multiply connected. No further assumptions on the regularity or geometric properties of the data are required. 

The weak localization property, namely item $(i)$ in Theorem~\ref{theo:main} states that most of vorticity's mass is contained in a ball centered in $X_i^\eps$ obeying \eqref{def:X_i} and provides precise rates of convergence. We largely improve the previously known weak localization estimates \cite{Mar3} and similarly in \cite{HLM24,DonatiLacaveMiot} in different geometric settings and references therein. Proposition~\ref{prop:Turkington} below provides even more precise estimates than Theorem~\ref{theo:main}. Specifically, for any $\eta>0$ it gives the existence of a domain of size $\eps^{1-\eta}$ containing most of the potential vorticity. However, interestingly enough, our method fails to provide a precise location for this domain, see Remark~\ref{rem:Turkington} for a detailed discussion. 

The strong localization property $(ii)$ in Theorem~\ref{theo:main} is pivotal for deriving the asymptotic motion and constitutes a substantial improvement compared to previously known results in \cite{Mar3, DonatiLacaveMiot, HLM24}. These previous results state that the potential vorticity remains supported in a strip of radial size $1/|\ln \eps|^{\kappa}$, for $\kappa < 1/4$. Here, we prove a localization of radial size $1/\sqrt{|\ln \eps|}$. This improvement is crucial to reach the critical size in order to observe a motion on spatial scales of order $1/\sqrt{|\ln \eps|}$ on which leapfrogging occurs. Due to the anisotropic nature of the problem leading to \eqref{eq:ring_v} in $x_z$-direction, the rings will cover a large distance, leading to confinement at order 1 in the $x_z$-direction. The radial confinement can also simply be written as
\begin{equation*}
 \big| r - X_{i,0,r}^\eps\big| \le \frac14\min_{j\neq k} |X_{j,r}(0) - X_{k,r}(0)|
\end{equation*}
by definition~\eqref{def:X_i} of the $X_i(t)$.

The localization property and its link to the small-time assumption for our main result are discussed in view of possible future extensions more in detail in Remark~\ref{rem:main} and a survey of previous results is provided in Section~\ref{sec:literature}.

\begin{remark}
Our method provides further information on the localization of vorticity. The center of potential vorticity of each vortex blob in the $2D$-reduction, namely the expression
\begin{equation}\label{def:beps}
 b^\eps_i(t) = \int_{\R^2_+} x \frac{ x_r \omega^\eps_i(t,x)}{\gamma} \dd x
\end{equation}
can be localized as follows. We obtain in Proposition~\ref{prop:full_dynamic} as a consequence of the aforementioned localization properties that
 \begin{equation*}
 \sup_{t \in [0,T_0]} \left|b_i^\eps(t) - X_i^\eps(t)\right| \le C_0\frac{\ln |\ln \eps|}{|\ln \eps|}.
 \end{equation*}
\end{remark}

\subsection{The limiting system}\label{sec:LP}

Several remarks concerning System~\eqref{def:X_i} are in order. While it is Hamiltonian, it is not known whether solutions are global due to a potential collapse of vortex rings at $r=0$. We prove that there exists $T_X>0$ only depending on $(Y_{i,0})_{1\le i \le N}$ defined in \eqref{distance-regime} but not on $\eps$ such that System~\eqref{def:X_i} has a unique solution on the time interval $[0,T_X]$, see Section~\ref{sec:limiting_traj}. We recall that all vortex rings are chosen to be of circulation $\gamma$ and hence it is convenient to consider the problem in a vertically moving frame in terms of $(Y_{i}^\eps)_{1\le i \le N}$ that are defined by the relation
 \begin{equation}\label{eq:LP-bis}
X_i^\eps(t) = X^* + \frac{\gamma}{4\pi r^*} te_z + \frac{1}{\sqrt{|\ln \eps|}} Y_i^\eps(t).
 \end{equation}
Further, we show that with that definition, the trajectories $t\mapsto Y_i^\eps(t)$ converge to limiting trajectories $t\mapsto \tilde{Y_i}(t)$ as $\eps \to 0$ which satisfy
\begin{equation}\label{eq:LP}
 \left\{\begin{aligned}
 & \frac{\dd \tilde{Y}_i}{\dd t} =\frac{\gamma}{2\pi} \sum_{j\neq i} \frac{(\tilde{Y}_i-\tilde{Y}_j)^\perp}{|\tilde{Y}_i-\tilde{Y}_j|^2} - \frac{\gamma}{4\pi (r^*)^2}\tilde{Y}_{i,r} e_z \\
 & \tilde{Y}_i(0)=Y_{i,0},
 \end{aligned}\right.
 \end{equation}
 In this system, the relative self-induced motion and the interactions are readily seen to be of the same order. Note the negative sign of the second term in \eqref{eq:LP} linked to the fluctuations in the relative self-induced vertical motion. The larger the radius $\tilde{Y}_{i,r}$ the larger is the deceleration compared to the uniform vertical translation, see Appendix~\ref{appendix:regimes} for a detailed discussion of the limit dynamics. In the literature, such as \cite{JerrardSmets,DavilaDelPinoMussoWei,Meyer}, it is \emph{this} dynamics that is considered as the motion of the rings in the appropriate scaling.

In this paper, we prefer to work with the trajectories $(X_i^\eps)_{1\le i \le N}$ and $(Y_i^\eps)_{1\le i \le N}$ as it is much more natural in the computations. We prove in Section~\ref{sec:motion} that subsequently, one could replace the trajectories $t\mapsto X_i^\eps(t)$ by 
 \begin{equation}\label{def:tildeX_i}
 \tilde{X}_i^\eps(t) := X^* + \frac{\gamma}{4\pi r^*} te_z + \frac{1}{\sqrt{|\ln \eps|}} \tilde{Y}_i(t)
 \end{equation}
 to match the notations of \cite{JerrardSmets,DavilaDelPinoMussoWei,Meyer}. All the results, in particular Theorem~\ref{theo:main} remain valid with that substitution.

 On a related note, we point out that the limiting dynamics \eqref{eq:LP} differs qualitatively and quantitatively from \eqref{def:X_i}. While, we can rule out collisions of points in \eqref{eq:LP} in finite time, see Appendix~\ref{app:nocollision}, global existence for solutions to \eqref{def:X_i} remains open. For further differences in the behavior of solutions for $0<\eps\ll 1$ and $\eps=0$ including numerical evidence, see Appendix~\ref{app:different dynamics}. Moreover, we recover the same limit dynamics \eqref{eq:LP} as for the leapfrogging in \cite{DavilaDelPinoMussoWei, GarciaHassainiaHmidi} and for the Gross-Pitaevskii equation in \cite{JerrardSmets} even though the dynamics for $\eps>0$ in \cite{JerrardSmets} appears to differ from the present one. Finally, we emphasize that neither the dynamics for $\eps=0$ nor $\eps>0$ considered here coincides with the ones studied in \cite{Mar24}. Specifically, the self-induced motion (not the relative one) is chosen to be of the same strength than the displacement due to interactions. In the respective limit dynamics \eqref{eq:LPMAR-appB} for $\eps=0$, the displacement linked to the vertical self-induced motion does not depend on the specific radius $Y_{i,r}$. We refer to Appendix~\ref{appendix:regimes} for details and figures in the case $N=2$.

\subsection{Survey of previous results and outlook}\label{sec:literature}

The interest in the motion of vortex rings and their hydrodynamic stability dates back to the work of Lord Kelvin \cite{Kelvin} and Helmholtz \cite{Helmholtz1}. Rigorous derivations of the motion of a single vortex ring in an inviscid ideal fluid are provided in \cite{Fraenkel, Benedetto2000}. The respective motion law to \eqref{eq:ring_v} for a single vortex ring was rigorously established for weakly viscous fluids in \cite{BrunelliMar, BuCaMarVisc, GallaySverak} where the recent paper \cite{GallaySverak} provides a very accurate description for large times. 

More generally, assuming that the vorticity is concentrated around a curve, the inviscid motion is expected to be asymptotically governed by the binormal curvature flow equation, as was formally computed by Da Rios \cite{daRios}, see \cite{JerrardSeis, FontelosVega} for partial results in this direction. The binormal curvature flow exhibits a rich panel of dynamics that have been intensively studied. In particular a series of papers by Banica and Vega is devoted to self-similar solutions or infinite polygonal lines solutions (we refer to \cite{BanicaVegaReview} for a review of the results). Finally, Jerrard and Smets \cite{JerrardSmetsbinormal} introduced a new formulation of the binormal curvature flow equation allowing for weaker solutions.

The issue of the persistence in time of such vortex filaments by the Euler flow is the so-called \emph{vortex filament conjecture}. While the conjecture remains open in its generality, results are available for specific geometric and dimension reduced settings. The rigorous asymptotic dynamics of straight filaments (the motion of Dirac masses of vorticity in the 2D Euler equations) was obtained in \cite{MarPul93}. In the aforementioned case of annular vortex rings evolving according to \eqref{eq:transport}, positive answers have been given with various assumptions. We refer to the work by Butt\`a, Cavallaro and Marchioro \cite{Mar3, Mar24} and references therein. 
This is in particular due to the fact that the problem allows for a dimension reduction leading to the 2D problem~\eqref{eq:Axi2D_rescaled}.

We also mention that the desingularization problem for the vortex rings has been addressed recently by Cao, Wan and Zhan \cite{Cao19} and by Davila, Del Pino, Musso and Wei \cite{DavilaDelPinoMussoWei}.

We emphasize that System~\eqref{eq:Axi2D_rescaled} notably differs from the 2D Euler equations. Indeed, the dynamics of 2D Euler \emph{point-vortices} do not include any self-induced motion of vortices. Even though Problem~\eqref{eq:Axi2D_rescaled} is reduced to two dimensions, as in the study of helical filaments in \cite{DonatiLacaveMiot} or point-vortices for the lake equations \cite{HLM24}, the Biot-Savart law (the relation giving $u$ in terms of $\omega$) is much more complex compared to the planar case. Additional techniques are required and the results developed for the 2D Euler equations (see for instance \cite{MarPul93}) cannot be applied. Roughly speaking, the motion of vortex rings is driven by the self-induced motion \eqref{eq:ring_v} and by their mutual interactions. Depending on the chosen asymptotic regime, the former maybe dominant, leading to a translation of the vortex rings along the symmetry axis \cite{Mar2, Mar3}. In case the latter dominates, the dynamics is reminiscent to the 2D dynamics of point vortices \cite{ButtaCavMar24}. In particular, it was shown in \cite{Mar3} that the interactions of vortex rings of different radii (independent of $\eps$) are bounded independently of $\eps$ and can be neglected compared to their self-interaction (inducing the velocity given by \eqref{eq:ring_v}) which is of order $|\ln \eps|$. The result of Theorem~\ref{theo:main} also improves those of \cite{Mar3} in the case of a single vortex ring.

As discussed, the \emph{leapfrogging} phenomenon first described by Helmholtz \cite{Helmholtz1} concerns the regime where both effect are balanced and hence requires mutual distances of order $1/\sqrt{|\ln \eps|}$ to each other, and a detailed explanation of the suitable scaling regime is provided in Appendix~\ref{sec:compRegimes}. The small mutual distance and strong interactions make the scenario hardly accessible by previously mentioned techniques.

To the best of our knowledge, the only rigorous results concerning the leapfrogging dynamics for the 3D Euler equations are \cite{DavilaDelPinoMussoWei, GarciaHassainiaHmidi, Mar24} and \cite{JerrardSmets} for the 3D Gross-Pitaevskii equation. While \cite{DavilaDelPinoMussoWei, GarciaHassainiaHmidi, JerrardSmets} consider the same scaling regime, the regime studied in \cite{Mar24} differs from the present one. The results from \cite{Mar24} cannot be adapted to our situation, as also pointed out by the authors in \cite[page 791]{Mar24}, see also Appendix~\ref{sec:compRegimes} for a comparison. In \cite{DavilaDelPinoMussoWei}, the authors construct initial data $\omega_0^\eps$ of a very specific class to obtain a smooth solution exhibiting leapfrogging vortex rings for finite times. After completion of our work, the result \cite{GarciaHassainiaHmidi} providing a global time-periodic solution of a leap-frogging pair of vortex rings has appeared. We emphasize the fact that the conclusions of Theorem~\ref{theo:main} hold for \emph{any} initial data $\omega_0^\eps \in L^\infty(\R^2_+)$ satisfying relations~\eqref{hyp:omega0}, which ranges from smooth compactly supported blobs to vortex patches, even possibly multiply connected. Opposite to \cite{DavilaDelPinoMussoWei, GarciaHassainiaHmidi}, our proof does not rely on sophisticated desingularization techniques.

While our result is limited to a fraction of a periodic motion but treats very general initial data, \cite{DavilaDelPinoMussoWei, JerrardSmets} provide solutions exhibiting several crossings. The method in \cite{JerrardSmets} for the Gross-Pitaevskii equation is based on Jacobian estimates not available for \eqref{eq:3DEuler}. In \cite{Mar24}, several crossings are observed upon choosing large radii but in a different scaling regime and shorter time scales.

In view of these results, additional assumptions on the general class of initial data \eqref{hyp:omega0} considered here might be necessary. Indeed, it is conceivable that the possible occurrence of filamentation of the vortex ring in the vertical direction leads to instabilities preventing an improvement of the $\mathcal{O}(1)$ localization in the $z$ direction for general initial data. Namely, a \emph{tail} in the wake of the ring may appear, a phenomenon that was described for instance in \cite{Moffatt,Pozrikidis} and recently proved in \cite{Choijeong2021, ChoiJeongHvortex, GuoJeongZhao, MarinRoulley}. This fact may explain similar restrictions in \cite{Mar3,Mar2,DonatiLacaveMiot,HLM24}. In \cite{DavilaDelPinoMussoWei}, this filamentation is not present since the authors work with specific initial data and a stable vorticity profile is constructed. In \cite{Mar24}, this filamentation does also not occur suggesting that the tail does not appear on the short time scale considered.
Consequently, it appears to be of great relevance to investigate the appearance of vortex filamentation and related potentially destabilizing effects.

\subsection{Strategy of the proof}
The 3D axisymmetric Euler equations without swirl allow for a dimension reduction in terms of a transport type equation for the vorticity coupled with an anelastic constraint for the velocity field \eqref{eq:Axi2D_rescaled}. The self-induced motion of a vortex ring is closely linked to this anelastic constraint and also captured by the local induction approximation (LIA), see also \cite{HLM24, DonatiLacaveMiot}. While this constitutes a main difference to the 2D setting, it is well known that the Biot-Savart law in the present setting shares some similarities with the 2D one even though it admits no explicit representation formula. Indeed, we revisit the well-posedness of \eqref{eq:Axi2D_rescaled} in Section~\ref{sec:WP} and provide a suitable expansion of the Biot-Savart law. This allows for a decomposition of the velocity field displaying the self-induced motion and mutual interactions of vortex rings \eqref{def:u_K,u_L,u_R}. Based on this decomposition and a boot-strap hypothesis on the strong localization up to some time $T_{\eps}>0$, we provide suitable estimates on the strength of the interactions of vortices. A common feature in previous papers \cite{Mar3,HLM24,DonatiLacaveMiot} is that vorticity which is initially localized in balls of radius $\eps$ is only known to be localized in balls of radius $|\ln\eps|^{-\kappa}$ for some $0<\kappa<1/4$ for positive times. Given the mutual distance of $|\ln\eps|^{-\frac12}$ this turns out to be insufficient for our purpose and needs to be significantly improved by upgrading both weak and strong localization properties.

The next step consists in showing the weak localization property. Commonly and in the aforementioned results, this property is inferred for the $i$-th vortex by exploiting the conservation of $L^p$-norms of the potential vorticity, the total energy as well as suitable bounds on the local energy of the $i$-th vortex and higher order moments such as the moment of inertia. In particular, the use of the moment of inertia $I_{\eps}(t)=\int|x-b_\eps(t)|^2x_r\omega_\eps(t)\dd x$ causes a loss of precision in the localization for positive times due to a repeated application of H\"older-inequalities. To overcome this issue, we revisit the variational approach by Turkington \cite{Tur87} developed for a single vortex evolving according to the 2D Euler flow. Given a precise expansion of the local energy and the conservation of the mass of the $i$-th vortex, it allows one to precisely localize the vorticity in weak sense. Several key difficulties arise: the techniques need to be adapted to the present setting with an anelastic constraint and several vortices displaying singular interactions. As these interactions may potentially lead to large variations in the energy of the $i$-th vortex, a precise estimate on the respective energy as required by Turkington is difficult to obtain. We overcome these difficulties among others by a multi-layered boot-strap argument. Specifically, assuming additionally a bootstrap hypothesis on the weak localization up to some time $0<\tilde{T_\eps}\leq T_\eps$, we are able to provide precise estimates on the center of vorticity and the energy of the $i$-th vortex.
With this information at hand, we show a substantially more precise weak localization property by adapting Turkington's argument. In particular, it follows $\tilde{T_\eps}=T_\eps$. We believe the extension of Turkington's approach to the setting of an anelastic constraint and strongly interacting vortices to be of independent interest.

A first consequence of this property is the derivation of the motion law up to time $T_\eps$ in Section~\ref{sec:motion}. Finally, we prove the desired strong localization property in Section~\ref{sec:strong_loc}. While no self-induced motion linked to the LIA occurs in the purely 2D setting, here it amounts to a fast vertical translation \eqref{eq:ring_v}. This anisotropic behavior together with the possible appearance of vortex filamentation suggests that a strong localization on scales of order $\mathcal{O}(\sqrt{|\ln\eps|}^{-1})$ for generic solutions can only be expected in the radial direction. Indeed, a small filament emanating from a vortex rings is expected to travel at speeds of order $\mathcal{O}(1)$ while the core of vortex rings travels at speed of order $\mathcal{O}(|\ln\eps|)$ in $x_z$-direction. Concerning the strong localization in $x_r$-direction, we require a localization of the support of order $\mathcal{O}(|\ln\eps|^{-\frac12})$.
In the literature, the proof of the strong localization is typically based on an iterative method first introduced by Marchioro and Pulvirenti in \cite{MarPul93, MP94} for $2D$ Euler flows. Up to today, the approach has reached an apparent complexity by being gradually upgraded and refined in a series of papers to various dimension-reduced settings and interacting vortices. Notably, in \cite{Mar3} and subsequently in \cite{HLM24, DonatiLacaveMiot} the anisotropic character of the problem has been dealt with by introducing directional localization estimates. Most recently, a double iterative procedure has been introduced in \cite{Mar24} to show an improved strong localization result sufficient to deal with a leapfrogging type dynamics, though on a shorter time scale, see Appendix~\ref{appendix:regimes} for a comparison. Our aim is twofold. First, we adapt this double iterative procedure to prove the strong localization at the required precision up to time $T_0$ independent of $\eps$, where we crucially rely on the improved weak localization estimates of Section~\ref{sec:weak_loc}.
Second and as a service to the reader not familiar with the refined iteration procedure, we present the method in a self-contained and accessible manner.

\begin{remark}\label{rem:main}
While our result crucially improves the strong radial localization compared to the literature enabling us to derive the leapfrogging motion for small times independent of $\eps$, the localization in the $x_z$ direction of order $\mathcal{O}(1)$ constitutes the main obstruction to observing a full period of the leapfrogging motion. Indeed, it does not allow one to show that the vortex blobs' supports remain disjoint for more than a fraction of the period. Specifically, vortex rings may exhibit a crossing but not the exchange of roles, see Figure~\ref{fig:crossing}. The possible occurrence of vortex filamentation for generic data appears to prevent a more precise strong localization over a full period. Rigorously proving the occurrence of filamentation in the present regime remains an open problem. The recent preprint \cite{GuoJeongZhao} discusses filamentation phenomena for a single vortex ring and less sharp initial concentration of vorticity.

\begin{figure}
\begin{minipage}[c]{0.48\textwidth}
\centering
 \includegraphics[height=4.4cm]{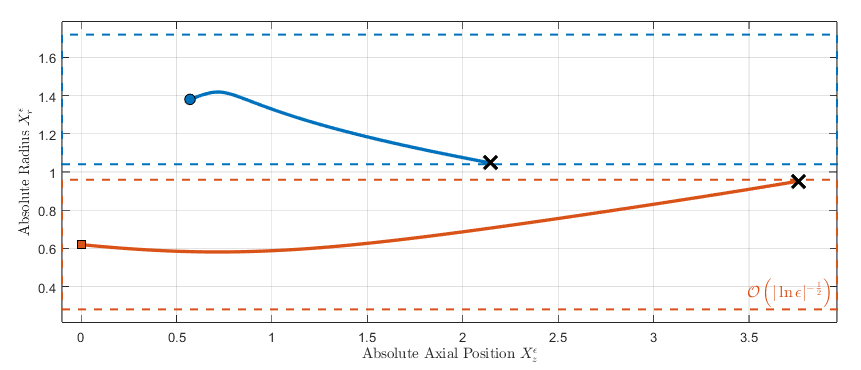}
\end{minipage}
\hfill
\begin{minipage}[c]{0.48\textwidth}
    \centering
 \includegraphics[height=4cm]{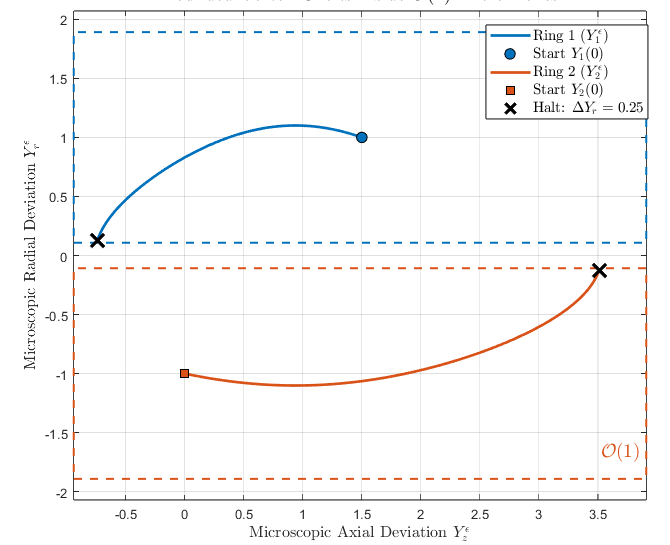}
 \end{minipage}
    \caption{The result of Theorem~\ref{theo:main} yielding confinement in strips. The pursuing vortex ring shrinks and overtakes the larger one. Then it starts enlarging and slowing down, while the other ring shrinks and accelerates. Theorem~\ref{theo:main} stops before both rings can have the same radius and exchange roles. The distance covered by the rings during that process is of order 1. The simulation solves equations \eqref{def:X_i} with $N=2$, $\eps=10^{-3}$, $Y_{0,1}=(1.5,1)^t$, $Y_{0,2}=(0,-1)^t$. On the left in the original frame of reference and on the right for $Y_i^{\eps}$ in the co-moving frame.}
    \label{fig:crossing}
\end{figure}

We expect our method to potentially reach a slightly larger fraction of a full period by revisiting the proof and optimizing constants. However, we do not pursue this issue here for the sake of a concise presentation and as it does not overcome the natural limitation due to the concentration scale in $x_z$.

While this suggests that one may not expect a strong localization property in both directions to hold at suitable spatial and temporal scales to describe a full period of motion, it is conceivable that the weak localization persists on larger time scales. These observations further motivate the study of the stability of a vortex ring under collision with a tail of another vortex ring. Different from our current method, this would require the weak localization to be obtained without assuming \emph{a priori} strong localization in a bootstrap argument. Note that all previous results regarding the control of general concentrated vortices in inviscid flows require the separation of vortex blobs \cite{Mar3, ButtaCavMar24, Mar24, DonatiLacaveMiot, HLM24,MarPul93, Meyer}. While we expect our new variational approach for multiple vortex rings developed in Section~\ref{sec:weak_loc} to be of independent interest for that purpose, this challenging mathematical problem is expected to require additional new tools. 
\end{remark}

\subsection{Organization of the paper}

The plan of the paper is the following. 
In Section~\ref{sec:WP}, we recall well-posedness results for \eqref{eq:Axi2D_rescaled} and provide a suitable expansion for the Biot-Savart law. 
In Section~\ref{sec:reduc}, we introduce useful notations and gather properties on the interactions of vortex rings. 
In Section~\ref{sec:weak_loc}, we obtain the weak localization property by new bootstrap argument and assuming the strong localization, in an unspecified region of space.
In Section~\ref{sec:motion}, we obtain information on the motion of the center of mass and obtain the weak localization Theorem~\ref{theo:main}-$(i)$.
In Section~\ref{sec:strong_loc}, we conclude the proof of Theorem~\ref{theo:main} by proving the strong localization, successively in the $r$ and $z$ direction. 
Appendix~\ref{app:lemma} provides two technical lemmas concerning rearrangement inequalities and a variant of Gronwall's inequality. Appendices~\ref{appendix:regimes}-\ref{app:limit dyn} discusses various asymptotic dynamics and provides a comparison to the regimes considered in \cite{Mar24, DavilaDelPinoMussoWei, JerrardSmets}.

\section{Well-posedness and Biot-Savart law}\label{sec:WP}

This section recalls and establishes well-posedness results for the axisymmetric 3D Euler equations without swirl. 

\subsection{2D reduction}
The evolution of an incompressible and inviscid three-dimensional fluid is governed by the standard 3D Euler equations
\begin{equation}\label{eq:3DEuler}
 \begin{cases}
 \partial_t \mathbf{U} + (\mathbf{U} \cdot \nabla) \mathbf{U} = -\nabla P, \\
 \div \mathbf{U} = 0,
 \end{cases}
\end{equation}
where for $(t,x) \in [0,T] \times \R^3$, the quantity $\mathbf{U}(t,x)$ denotes the velocity of the fluid at time $t$ at the position $x$. Letting $\mathbf{\Omega} = \nabla \wedge \mathbf{U}$ denote the vorticity of the fluid, \eqref{eq:3DEuler} can be expressed as
\begin{equation}\label{eq:3Dvorticity}
 \begin{cases}
 \partial_t \mathbf{\Omega} + (\mathbf{U}\cdot \nabla) \mathbf{\Omega} = (\mathbf{\Omega} \cdot \nabla) \mathbf{U} \vspace{2mm}\\
 \displaystyle \mathbf{U}(t,x) = -\frac{1}{4\pi}\int_{\R^3} \frac{x - y}{|x-y|^3} \wedge \mathbf{\Omega}(t,y) \dd y,
 \end{cases}
\end{equation}
where $\mathbf{U}$ is uniquely determined by $\mathbf{\Omega}$ for all $t\in [0,T]$ by means of the Biot-Savart law, i.e. the second identity in~\eqref{eq:3Dvorticity}. In the present paper, we restrict our attention to solutions that are axisymmetric vector fields without swirl. Thus, it is conveniently studied in cylindrical coordinates in $\R^3$ which we denote by $(z,r,\theta)$ with basis vectors $(e_z,e_r,e_\theta)$. A vector field $\mathbf{U}: \R^3\rightarrow \R^3$ is called \emph{axisymmetric without swirl} if of the form 
\begin{equation*}
 \mathbf{U}(z,r,\theta)=U_z(z,r)e_z+U_r(z,r)e_r,
\end{equation*}
i.e. $\mathbf{U}$ is independent of the toroidal variable $\theta$ (axisymmetric) and $\mathbf{U}\cdot e_{\theta}=0$ (without swirl). This implies that the vorticity $\mathbf{\Omega}=\nabla\wedge \mathbf{U}$ is a purely toroidal vector field:
\begin{equation*}
 \mathbf{\Omega}=\Omega e_\theta, \quad \text{with} \,\, \Omega=\partial_z U_r-\partial_r U_z\in \R.
\end{equation*}
Here and below we adopt the notation $x\in \R_{+}^2=\R\times (0,\infty)$ with $x=(x_z,x_r)$. Note the change of convention w.r.t. the usual cylindrical coordinates. This allows one to reduce \eqref{eq:3DEuler} to a 2D problem by defining the scalar $\omega: \R_+^2\rightarrow \R$ and vector field $u:\R_{+}^2\rightarrow \R^2$ respectively as
\begin{equation*}
 \omega=\frac{\Omega}{x_r}\in \R, \qquad u=(U_z,U_r)\in \R^2. 
\end{equation*}
Note that the div-curl problem $\div \mathbf{U}=0$ and $\nabla \wedge \mathbf{U}=\mathbf{\Omega}$ reduces to
\begin{equation*}
 \partial_{r}u_{r}+\frac{1}{x_r} u_{r}+\partial_{z} u_{z}=0, \quad x_r\omega=\partial_z u_{r}-\partial_{r} u_{z}. 
\end{equation*}
In particular, the elliptic problem posed on $\R_+^2=\R\times (0,\infty)$ also writes
\begin{equation}\label{eq:elliptic}
 \div_{(z,r)}(x_r u)=0, \quad \curl_{(z,r)}(u)=x_r\omega.
\end{equation}
Under suitable boundary and regularity conditions, the two-dimensional velocity field $u$ is hence uniquely determined by $\omega$, see Section~\ref{sec:BiotSavart} for details. On the half plane $\R_+^2$, the potential vorticity $\omega$ satisfies the transport type equation
\begin{equation}
 \label{eq:transport} \partial_t \omega +u\cdot \nabla \omega=0,
\end{equation}
that also reads as continuity equation
\begin{equation}\label{eq:continuity}
 \partial_t(x_r\omega)+\div(x_ru\omega)=0,
\end{equation}
given that $\div(x_ru)=0$ from \eqref{eq:elliptic}.

\subsection{Well-posedness results}\label{sec-WP}

Global existence and uniqueness of weak axisymmetric without swirl solutions of \eqref{eq:3DEuler} has been proved in \cite{UhkovskiiYudovich} and \cite[Theorem 3.3]{SaintRaymond}, \cite{AbidiHmidiKeraani, Danchin} as well as references therein. In \cite[Theorem 5.4]{cetrone-serafini}, global existence and uniqueness of the weak vorticity solution of \eqref{eq:transport} is proved when the initial velocity is axisymmetric without swirl, the initial vorticity satisfies $\Omega_0\in L^1\cap L^\infty(\R_+^2)$ and its support does not overlap the symmetry axis. More precisely, this solution satisfies 
$$\frac{\Omega(t,X(t,0,x))}{X_r(t,x)}=\frac{\Omega_0(x)}{x_r}$$ where the flow map $X:[0,\infty)\times \R^3\to \R^3$ is defined by 
 \begin{equation*}
 \frac{\dd X(t,t_0,x)}{\dd t}=U\left(t,X(t,t_0,x)\right), \quad X(t_0,t_0,x)=x.
 \end{equation*}
 From now on we shall write $X(t,x)=X(t,0,x)$ for simplicity.
 
Since the velocity is axisymmetric without swirl, the flow map $X$ may be seen as a map on $\R_+ \times \R^2_+$ to $\R^2_+$ namely $X=X_z(t,x_z,x_r)e_z+X_r(t,x_z,x_r)e_z$, where 
\begin{equation}\label{eq:lagrangian-axi}\begin{split}
 \frac{\dd X_z(t,x_z,x_r)}{\dd t}&=u_z\left(t,X(t,x_z,x_r)\right), \quad X_z(0,x_z,x_r)=x_z,\\
 \frac{\dd X_r(t,x_z,x_r)}{\dd t}&=u_r\left(t,X(t,x_z,x_r))\right), \quad X_r(0,x_z,x_r)=x_r.\end{split}
 \end{equation}
In the solution constructed by \cite[Theorem 5.4]{cetrone-serafini}, the velocity is reconstructed from the vorticity with the Biot-Savart Kernel: $U(t,x)=\int H(x-y)\wedge \mathbf{\Omega}(t,y)\,\dd y$ with $H(x)=-x/(4\pi |x|^3)$ for $x,y\in \R^3$ and with $\mathbf{\Omega}$ depending only on $(z,r)$ and purely toroidal. Writting this 3D Biot-Savart law in cylindrical coordinates, we have an explicit representation\footnote{We may note that a $\cos \theta$ is missing in \cite{Mar3,Mar2}.}:
 \begin{equation}\label{eq:b-s-axi}
 \begin{split}
 u_z(t,x)&=-\frac{1}{4\pi}\int \, \dd y_z \int_0^\infty y_r^2\,\dd y_r \int_{-\pi}^\pi \, \dd \theta \frac{\omega(t,y)(x_r\cos\theta-y_r)}{[|x-y|^2+2x_ry_r(1-\cos\theta)]^{3/2}},\\
 u_r(t,x)&=\frac{1}{4\pi}\int \, \dd y_z \int_0^\infty y_r^2\,\dd y_r \int_{-\pi}^\pi \, \dd\theta \frac{\omega(t,y)(x_z-y_z)\cos\theta}{[|x-y|^2+2x_ry_r(1-\cos\theta)]^{3/2}},
 \end{split}
 \end{equation} 
 where $x=(x_z,x_r)$ and $y=(y_z,y_r)$ are vectors in $\R^2_+$.
 Moreover, it is classical that $u$ is bounded (see \eqref{ineq:bound-u} below) and log-lipschitz (by adapting, for instance, the planar case done in \cite[Appendix 2.3]{MP94}). Finally, it is established in \cite[Corollary 5.6]{cetrone-serafini} that the solution is a weak solution of \eqref{eq:transport} in the sense of distributions (see (5.3) in \cite{cetrone-serafini}). 
We summarize these results for
\begin{equation}\label{eq:Axi2D}
\begin{cases}
 \partial_t \, \omega + u\cdot \nabla \omega = 0 & \text{ in }\R_+^* \times \R^2_+ , \vspace{1mm}\\
 \div (x_r u) = 0 & \text{ in } \R_+ \times \R^2_+ ,\vspace{1mm} \\
 \curl (u)=x_r \omega & \text{ in } \R_+ \times \R^2_+, \vspace{1mm} \\
 \omega(0,\cdot) = \omega_0 & \text{ in } \R^2_+.
\end{cases}
\end{equation} 
on the half plane $\R_+^2$ in the next proposition.

\begin{proposition}[\cite{cetrone-serafini, GallaySverak}]\label{prop:WP}
 Let $\Omega_0\in L^1\cap L^{\infty}(\R_+^2)$ and assume that $\inf_{x\in \supp(\Omega_0)}x_r=\delta_0>0$. Let $\omega_0=\Omega_0/x_r$ so that $\omega_0\in L^1\cap L^{\infty}(\R_+^2)$. For all $T>0$ there exists a unique pair $(\omega,u)$ with $x_r\omega\in L^{\infty}([0,T];L^1\cap L^{\infty} (\R_+^2))$, $x_r^2\omega\in L^{\infty}([0,T];L^1(\R_+^2))$, $\omega\in L^{\infty}([0,T]; L^{\infty}(\R_+^2))$, and $u\in L^\infty([0,T]\times \R_+^2)$, solution to \eqref{eq:transport} in the following sense
 \begin{enumerate}[(i)]
 \item $u=(u_z,u_r)$ is given in terms of $\omega$ by the axisymmetric Biot-Savart law \eqref{eq:b-s-axi}.
 \item The vorticity equation is satisfied in the sense of distribution, i.e. for all $\Phi\in C^{\infty}_c([0,T]\times \R_+^2)$ and for any $t\in [0,T]$
 \begin{equation*}
 \int_{\R_{+}^2}\Phi(t,x) x_r\omega(t,x)\dd x-\int_{\R_+^2}\Phi(0,x)x_r\omega_0(x)\dd x=\int_0^t\int_{\R_+^2}x_r\omega(s,x)\left(\partial_t\Phi+u\cdot\nabla \Phi\right)(s,x)\dd x\dd s.
 \end{equation*}
 \item We have $\omega(t,X(t,x))=\omega_0(x)$, where $X$ is the flow of $u$ defined in \eqref{eq:lagrangian-axi}\footnote{Such a flow exists in the classical sense by log-Lipschitz regularity of $u$.}.
 \end{enumerate}
 Moreover, we have the estimates
\begin{equation}\label{ineq:bound-u}
 \|u\|_{L^\infty(\R_{+}^2)}\leq C \|x_r^2\omega\|_{L^1(\R_+^2)}^{1/2}\|\omega\|_{L^\infty(\R_+^2)}^{1/2}
\end{equation}
and
\begin{equation*}
 \left|u(t,x)-u(t,y)\right|\leq C(\|x_r\omega\|_{L^1\cap L^{\infty}([0,T]\times\R_+^2)})|x-y|\left(1+\left|\ln|x-y|\right|\right).
\end{equation*}
\end{proposition}

A direct consequence of the log-Lipschitz estimate of the velocity field is the following:

\begin{corollary}\label{coro:flow}
 Let $\omega_0\in L^1\cap L^\infty (\R_+^2)$ such that $x_r\omega_0\in L^1\cap L^\infty (\R_+^2)$ and $\delta_0:=\inf_{x\in \supp(\omega_0)}x_r>0$. Let $T>0$ and $(\omega,u)$ be the unique weak solution of \eqref{eq:Axi2D} on $[0,T]$ provided by Proposition~\ref{prop:WP}. Then the following holds true.
 \begin{enumerate}[(i)]
 \item There exists $\delta_T>0$ depending only on $\|\omega_0\|_{L^\infty}$, $\delta_0$ and $T$ such that 
 \begin{equation*}
 \supp\omega(t,\cdot)\subset \{x\,|\, x_r\geq \delta_T\}, \quad \forall t\in [0,T].
 \end{equation*} 
 \item If $\omega_0$ has compact support in $\R_+^2$, denoting by $R_0=\sup \{|x|\,|\, x\in \supp(\omega_0)\}$, there exists a compact set $\cal{K}_T\subset \R_+^2$ depending only on $\|x_r\omega_0\|_{L^\infty}$, $\delta_0, R_0$ and $T$ such that 
 \begin{equation*}
 \supp\omega(t,\cdot)\subset \cal{K}_T, \quad \forall t\in [0,T].
 \end{equation*}
 \item For all $t\in [0,T]$ and $p\in [1,\infty)$ it holds $\|x_r^{\frac{1}{p}}\omega(t,\cdot)\|_{L^p}=\|x_r^{\frac{1}{p}}\omega_0\|_{L^p}$.
\end{enumerate}
\end{corollary}

\begin{proof}
We provide a quick proof for sake of completeness although this is probably a standard fact. Since $\omega$ is constant along trajectories it suffices to provide estimates on the Lagrangian flow. Let $x\in \supp(\omega_0)$. Writing $X(t)=X(t,x)$ and observing that $u_r(t,x_z,0)=0$ in view of \eqref{eq:b-s-axi} we have
\begin{align*}
 \left| \frac{\dd X_r(t)}{\dd t}\right|&=\left|u_r\left(t,X(t)\right)\right|
 =\left|u_r\left(t,X_z(t),X_r(t)\right)-u_r\left(t,X_z(t),0\right) \right|\\
 &\leq C |X_r(t)|(1+|\ln |X_r(t)||)
\end{align*}
where we have used the log-Lipschitz property for $u$. Hence this proves $(i)$ by means of a Gronwall-type argument (see e.g. the arguments to obtain (6.26) from (6.23) in \cite{MP94}). Property $(ii)$ is a straightforward consequence of the $L^\infty$ bound on the velocity. Finally, $(iii)$ follows from the continuity equation \eqref{eq:continuity}.
\end{proof}

\begin{remark}\label{rmk: decomposition vorticity}
Finally, one obtains the scaled system \eqref{eq:Axi2D_rescaled} from \eqref{eq:Axi2D} by considering the logarithmic time scale $s=|\ln\eps|t\in [0,T]$ and $\omega^{\eps}(s,x)=\omega\left(\frac{s}{|\ln\eps|},x\right)$. 

Further, for initial vorticity as in \eqref{hyp:omega0} consisting of $N$ disjoint blobs, the existence of a unique Lagrangian flow allows one to decompose the unique weak solution as 
\begin{equation}
 \omega^\eps(t,x)=\sum_{i=1}^N\omega_i^\eps(t,x), \qquad \omega_i^{\eps}(t,X(t,\cdot))=\omega_{i,0}^\eps.
 \label{dec:flow}
\end{equation}
Moreover, Part \textit{(i)} of Corollary~\ref{coro:flow} states that there exists a positive time such that a vorticity that is initially localized away from the symmetry axis does not reach the axis and remains localized. Concerning the localization in the radial direction, a quantitative estimate of this type is stated in \cite[Section 2.3]{Choijeong2021}, it is shown that the radius of the support
\begin{equation*}
 R_\eps(t)=\sup\{x_r : x\in \supp(\omega^\eps(t,\cdot))\}
\end{equation*}
grows at most as 
\begin{equation*}
 R_\eps(t)\leq C_\eps\Big(1+\frac{t}{|\ln\eps|}\Big)^2,
\end{equation*}
where $C_\eps=C(\|\omega_0^\eps\|_{L^1(\R^3)}, \|\omega_0^\eps\|_{L^{\infty}(\R^3)})>0$. Unfortunately, these estimates of the radius and of the distance to the axis depend on $\|\omega_0^\eps\|_{L^{\infty}(\R^3)}$ which blows up as $\eps^{-2}$ when $\eps\to 0$.
Hence, we will not be able to use this kind of estimate and we will later define a time $T_\eps$ such that the support of $\omega^\eps$ is included in a strip $\R\times[r^*/2,2 r^*]$. One of the main objectives is to get a strong localization property to prove that $T_\eps \geq T_0$ for all $\eps\in (0,\eps_0]$, provided that $T_0$ and $\eps_0$ are chosen small enough.
\end{remark}

\subsection{Expansion of the Biot-Savart law}\label{sec:BiotSavart}
In the 3D axisymmetric framework, the Biot-Savart law admits an explicit representation formula \eqref{eq:b-s-axi}. A suitable expansion of the Biot-Savart kernel used in sequel is introduced. We notice that the 2D div-curl problem \eqref{eq:elliptic} including the weighted incompressibility condition yields that there exists a stream function $\Psi: \R_+^2\rightarrow \R$ such that 
\begin{equation}\label{eq:2Du}
u=\begin{pmatrix}u_z\\u_r\end{pmatrix} = \frac1{x_r} \nabla^\perp \Psi = \frac1{x_r}\begin{pmatrix}-\partial_r \Psi\\ \partial_z \Psi\end{pmatrix}.
\end{equation}
Following \cite{FengSverak, GallayConfl} the boundary conditions $\Psi(x_z,0)=0$ are imposed. We are led to solve 
\begin{equation}\label{eq:elliptic2D}
 \begin{cases}
 \displaystyle \div \left(\frac{1}{x_r}\nabla \Psi(x)\right) = x_r\omega(x) \\
 \Psi(x_z,0)=0.
 \end{cases}
\end{equation}
In relation with \eqref{eq:b-s-axi}, there is an explicit formula of the Green kernel, namely it is given in \cite[Section 2]{GallayConfl} that 
\begin{equation}\label{eq:stream}
 \Psi(x)=\int_{\R_+^2}\mathcal{G}(x,y)y_r\omega(y)\dd y,
 \end{equation}
where 
\begin{equation}\label{eq:F}
 \mathcal{G}(x,y)= -\frac{1}{2\pi}\sqrt{x_ry_r}F\left(\frac{|x-y|^2}{x_ry_r}\right) \quad \text{with} \quad F(s)=\int_{0}^{\frac{\pi}{2}}\frac{\cos(2\theta)}{\sqrt{\sin^2(\theta)+\frac{s}{4}}}\dd \theta, \quad s>0.
\end{equation}
This function $F$ was intensively studied in the literature. We bring together in the next lemma some estimates we need, and we refer to \cite[Lemma 2.1]{GallayConfl} and \cite[Lemma 2.7]{FengSverak} for more details.

\begin{lemma}\label{lem:F}
For any $R_0>0$ given, there exists $C_0>0$ such that the following estimates holds for any $s\in (0,R_0]$:
\[
\max\Bigg(\Big|F(s)+\frac{1}{2}\ln(s)\Big|, \frac1{|\ln s|}\Big|F'(s) + \frac1{2s}\Big|, |s^2F''(s)| \Bigg)\leq C_0.
\]
\end{lemma}

The explicit formula of the leading term of $F$ allows for a decomposition of the Green function $\mathcal{G}$. 

\begin{corollary}\label{cor:expansion Green function}
 We have the following decomposition
 \begin{equation*}
 \mathcal{G}(x,y)=\frac{1}{2\pi}\sqrt{x_ry_r}\ln|x-y|+S(x,y) \quad \text{with} \,\, S\in W^{1,\infty}_{\loc}(\R^2_+\times \R^2_+).
 \end{equation*}
\end{corollary}

\begin{proof}
 We define 
 \begin{equation*}
 S(x,y):=\mathcal{G}(x,y)-\frac{1}{2\pi}\sqrt{x_ry_r}\ln|x-y|=-\frac{1}{2\pi}\sqrt{x_ry_r}\Bigg(F\left(\frac{|x-y|^2}{x_ry_r}\right)+\ln|x-y|\Bigg).
 \end{equation*}
 Let $\cal{K}$ a compact subset of $\R^2_+$, then there exists $R_0>0$ such that for any $x,y\in \cal{K}$ we have $\frac{|x-y|^2}{x_ry_r}\leq R_0$.
 Hence, by Lemma~\ref{lem:F}, one has for all $x,y\in \cal{K}$
 \begin{equation*}
 \Bigg|F\left(\frac{|x-y|^2}{x_ry_r}\right)+\ln|x-y| \Bigg|=\Big|\ln(\sqrt{x_ry_r})+\mathcal{O}(1)\Big|
 \leq C_{\cal{K}}. 
 \end{equation*}
 For the derivative, we compute
 \begin{align*}
 \nabla_x F\left(\frac{|x-y|^2}{x_ry_r}\right)&=\Big(\frac{2(x-y)}{x_ry_r} - \frac{|x-y|^2 e_r}{x_r^2y_r}\Big)F'\left(\frac{|x-y|^2}{x_ry_r}\right)\\
 &=-\frac{x-y}{|x-y|^2}+\mathcal{O}\left(\frac{1}{\sqrt{x_r y_r}}\right) +\mathcal{O}\left(\frac{1}{x_r}\right),
 \end{align*}
 which implies that
 \[
 \Bigg|\nabla_x F\left(\frac{|x-y|^2}{x_ry_r}\right)+\nabla_x\ln|x-y| \Bigg| \leq C_\cal{K}
 \]
 for all $x,y\in \cal{K}$. As we have the same estimates for $\nabla_y$, this ends the proof.
 \end{proof}

This expansion of the Green's function together with \eqref{eq:2Du} allows one to derive a decomposition of the velocity field $u$. For further reference, we denote the most singular part of the stream function by
\begin{equation}\label{def:psi}
 \psi^\eps(t,x)= \psi[\omega^\eps](t,x) := \int_{\R^2_+} \frac{\sqrt{x_ry_r}}{2\pi}\ln|x-y|y_r \omega^\eps(t,y)\dd y.
\end{equation}
As will be detailed later, this quantity plays a central role in determining the energy and the local induction for the propagation of vortex rings.

\begin{lemma}\label{lem:decomp_u}
Let $\cal{K}$ a compact subset of $\R^2_+$ and $\omega \in L^{\infty}(\R_+^2)$ with $\supp(\omega)\subset \cal{K}$. There exists $C_{\cal{K}}>0$ depending only on $\cal{K}$ such that the velocity field $u$ defined by \eqref{eq:2Du} satisfies 
\begin{equation*}
 u := u_K + u_L + u_R
\end{equation*}
with
\begin{equation}\label{def:u_K,u_L,u_R}
 \begin{split}
 u_K(x) & = \frac{1}{x_r}\int_{\R^2_+} \frac{\sqrt{x_ry_r}}{2\pi} K(x-y) y_r\omega(y)\dd y, \quad \text{with } \, K(x,y)=\nabla^{\perp}\ln|x-y|=\frac{(x-y)^{
 \perp}}{|x-y|^2}, \\
 u_L(x) & = - \frac{e_z}{4\pi x_r} \int_{\R^2_+} \sqrt{\frac{y_r}{x_r}}\ln|x-y| y_r\omega(y)\dd y = -\frac{e_z}{2x_r^2}\psi[\omega](x), \\
 u_R(x) & = \frac{1}{x_r}\int_{\R^2_+}\nabla_x^\perp S(x,y)y_r\omega(y)\dd y,
 \end{split}
\end{equation}
where $\|u_R\|_{L^{\infty}(\cal{K})}\leq C_\cal{K} \int y_r\omega(y)\dd y$.
\end{lemma} 

\begin{proof}
By the decomposition of Corollary~\ref{cor:expansion Green function}, we have
\begin{equation*}
 \nabla^{\perp}_x \mathcal{G}(x,y)=\frac{\sqrt{x_ry_r}}{2\pi}\frac{(x-y)^{\perp}}{|x-y|^2}-\frac{1}{4\pi}\sqrt{\frac{{y_r}}{{x_r}}}\ln|x-y|e_z+\nabla^{\perp}_x S(x,y),
\end{equation*}
hence the decomposition of $u$ comes from \eqref{eq:2Du} and \eqref{eq:stream}.
The desired estimate of $u_R$ follows the estimate of $S$ in Corollary~\ref{cor:expansion Green function}.
\end{proof}

Note that the most singular term $u_K$ in the decomposition stems from the analogue of the 2D Biot-Savart law. It corresponds to the standard spinning around the filament and does not contribute to the motion of the vortex rings. The second term given by $u_L$ is linked to the radial correction of the 2D Biot-Savart law in the present axisymmetric setting. It points in the vertical $e_z$-direction and its speed amounts to $\frac{-1}{2x_r^{2}}\psi^{\eps}(t,x)$ where the latter is closely linked to local energy, see Lemma~\ref{lem:encadrement_psi} below. In the absence of strong interactions, it accounts for the leading order displacement of the vortex ring according to the LIA or binormal curvature flow, see \eqref{eq:ring_v}. Finally, the expansion of the Biot-Savart kernel allows one to treat $u_R$ as term of lower order.

\begin{remark}
In dimension reduced settings like the present one, it is well-known that the $3D$-elliptic problem translates to a weighted $2D$ elliptic problem like \eqref{eq:elliptic2D}. For instance, this also occurs for the (degenerate) lake equations and the $3D$ Euler equations under helicoidal symmetry. A common feature is that the unique solution to the respective weighted elliptic problems is given by 
\begin{equation*}
 \Psi(x) = \int_{\R^2_+}\mathcal{G}(x,y)y_r\omega(y)\dd y
\end{equation*}
where the Green function shares similarities with the classical one on $\R^2$. For the lake equations on bounded domains, it was proven in \cite{DekeyserVanS, AlTakiLacave} that
\begin{equation*}
 \mathcal{G}_{\mathrm{lake}}(x,y)=\frac{\sqrt{b(x)b(y)}}{2\pi}\ln |x-y|+S(x,y) \quad \text{with} \quad S\in W^{1,\infty}_{\loc}(\Omega^2) ,
\end{equation*}
is the Green kernel associated to the following elliptic problem
\begin{equation*}
 \begin{cases}
 \displaystyle \div \left(\frac{1}{b(x)}\nabla \Psi(x)\right) = b(x)\omega(x) \\
 \Psi\vert_{\partial\Omega}=0.
 \end{cases}
\end{equation*}
Note that the 3D axisymmetric Euler equations without swirl can formally be recovered as special case of the lake equations by setting $b(x_z,x_r)=x_r$ and considering a lake given by the half-plane. A similar result holds for the $3D$ Euler equations under helicoidal symmetry \cite{DonatiLacaveMiot}, but with a more complicated expansion due to the fact that the anelastic constraint is not a scalar function $1/b(x)$, but a matrix. A key advantage of such an expansion of the respective Green function lies in the decomposition of the associated velocity field, see Lemma~\ref{lem:decomp_u}, that highlights the link to the LIA, namely the contribution of $u_L$, see \cite{DekeyserVanS, HLM24, DonatiLacaveMiot}. 
\end{remark}

\section{A vortex ring under the influence of the other rings} \label{sec:reduc}

Our strategy for proving Theorem~\ref{theo:main} consists in focusing on the dynamics of a single blob of vorticity, by considering the effect of the other blobs as an exterior field. To proceed, we place ourselves in the framework of Theorem~\ref{theo:main}. Specifically, let $X^*= (z^*, r^*) \in \R^2_+$ and $(Y_{i,0})_{1 \le i \le N}$ be $N$ points in $\R^2$ such that $Y_{i,0,r} \neq Y_{j,0,r}$ for all $i\neq j$ and let $X_{i,0}^{\eps}$ given by \eqref{distance-regime}. Finally, let $M_0>0$ be fixed and $\omega_0^\eps \in L^\infty (\R^2_+)$ satisfying \eqref{hyp:omega0} for every $\eps > 0$. 
Further, we denote by $(u^\eps,\omega^\eps)$ the unique weak solution of \eqref{eq:Axi2D_rescaled} with initial data $\omega_0^\eps$ on $[0,T_0)$ where $T_0$ is arbitrary large (in the sense of Proposition~\ref{prop:WP}). The vorticity can be decomposed as $\omega^\eps(t,\cdot) = \sum_{i=1}^N \omega_i^\eps(t,\cdot)$ where $\omega_i^\eps$ corresponds to the transport of $\omega_{i,0}^\eps$, see Remark~\ref{rmk: decomposition vorticity}. For each $i \in \{1,\ldots,N\}$, we introduce the exterior field $F^\eps_i : [0,+\infty)\times \R^2_+ \to \R^2$ given by
\begin{equation}\label{def:F}
 F_i^\eps(t,x) = \frac1{x_r}\sum_{j \neq i} \int_{\R^2_+} \nabla_x^\perp\mathcal{G}(x,y) y_r\omega_j^\eps(t,y)\dd y,
\end{equation}
so that the $i$-th blob $\omega_i^\eps$ satisfies the following equations:
\begin{equation}\label{eq:Axi2D_rescaled_F}
\begin{cases}
 \ds \partial_t \, \omega_i^\eps + \frac{1}{|\ln \eps|}(u_i^\eps + F_i^\eps)\cdot \nabla \omega_i^\eps = 0 & \text{ in } \R_+^* \times \R^2_+, \vspace{1mm}\\
 u_i^\eps = \frac{1}{x_r}\nabla^\perp \Psi_i^\eps & \text{ in } \R_+ \times \R^2_+ ,\vspace{1mm} \\
 \div \big( \frac{1}{x_r} \nabla \Psi_i^\eps \big) = x_r \omega_i^\eps & \text{ in } \R_+ \times \R^2_+ ,\quad \Psi_i^\epsilon=0\text{ on } \R_+ \times \partial \R^2_+ , \vspace{1mm} \\
 \omega_i^\eps(0,\cdot) = \omega_{i,0}^\eps & \text{ in } \R^2_+.
\end{cases}
\end{equation}

As $\div (x_r(u_i^\eps +F_i^\eps))=0$, it is clear that the total mass of vorticity is conserved:
\begin{equation}\label{eq:consgamma}
 \int_{\R^2_+} x_r \omega_i^\eps(t,x)\dd y =\int_{\R^2_+} x_r \omega_{i,0}^\eps(t,x)\dd y=\gamma,\qquad \forall t\in [0,T_0].
\end{equation}

For any $r_0,h_0,d > 0$ and $\eps>0$ we define the rectangles $\mathcal{A}_{\eps}^{r_0}(d,h_0)$ corresponding to hollow cylinders in the 3D setting as
\begin{equation}\label{def:Aeta}
 \mathcal{A}_{\eps}^{r_0}(d,h_0) = \left\{ (z,r) \in \R^2_+, |r-r_0 | \leq \frac{d}{\sqrt{|\ln\eps}|} \text{ and } |z-z^*|\leq h_0 \right\}.
\end{equation}
The height $h_0>0$ will be chosen later, and we set
\begin{equation}\label{def:min-dist}
 d_0 := \frac{1}{4}\min_{j\neq k} |Y_{j,0,r}-Y_{k,0,r} |.
\end{equation}
From the assumption \eqref{hyp:omega0} on $\omega_0^\eps$, we have for every $i \in \{1,\ldots,N\}$,
\[
\supp \omega_{i,0}^\eps \subset B(X^\eps_{i,0},\eps) \subset \mathcal{A}_{\eps}^{X^\eps_{i,0,r}}\left(\frac{d_0}{4},\frac{h_0}{4}\right)
\]
as soon as $\eps \le \eps_0$ for $\eps_0>0$ depending on $h_0$ and $d_0$.

Let $T_0>0$ to be chosen later. We define
 \begin{equation}\label{def:Teps}
 T_\eps := \sup \Big\{ t\in [0,T_0],\, \forall s\in [0,t],\ \forall i \in \{1,\ldots,N\},\ \supp \omega_i^\eps(s,\cdot) \subset \mathcal{A}_{\eps}^{X^\eps_{i,0,r}}\left(d_0,h_0\right) \Big\},
 \end{equation}
which verifies, by continuity of the trajectories, $T_\eps > 0$ for every $\eps \in (0,\eps_0]$. Note that in particular, up to choosing $\eps_0$ small enough, we have
\begin{multline}\label{est:supp}
 \supp \omega_i^\eps(t,\cdot) \subset\mathcal{A}_{\eps}^{X^\eps_{i,0,r}}\left(d_0,h_0\right) \subset [z^*-h_0,z^*+h_0]\times[r^*/2,2r^*],\\ \quad \forall t\in [0,T_\eps],\ i\in \{1,\ldots, N\}, \ \eps\in (0,\eps_0].
\end{multline}
Next, we prove that the influence on $\omega_i^\eps$ of the other vortex rings is not too singular. More precisely we have the following lemma.

\begin{lemma}\label{lem:hyp_F}
 Provided $h_0\leq 1$, there exists $\eps_0>0$ depending only on $X^*$, $(Y_{i,0})$ and $M_0$ such that for every $\eps \in (0,\eps_0]$, every $i \in \{1,\ldots,N\}$, every $x,y \in\mathcal{A}_{\eps}^{X^\eps_{i,0,r}}\left(d_0,h_0\right)$ and every $t\in [0,T_\eps]$,
\begin{equation*}
 \left\{\begin{aligned}
 & \div (x_r F_i^\eps(t,x)) = 0,\\
 & | F_i^\eps(t,x) | \le C_F\sqrt{|\ln\eps|} \\
& \big| F_i^\eps(t,x) - F_i^\eps(t,y) \big| \le C'_F |\ln\eps| |x-y| ,
\end{aligned}\right.
\end{equation*}
with
\begin{equation*}
 C_F = \frac{\gamma C^*}{d_0} \qquad \text{ and } \qquad C_F' = \frac{\gamma C^*}{d_0^2},
\end{equation*}
and where $C^*$ depends only on $X^*$.
\end{lemma}

\begin{proof}
First, we note that $ \div (x_r F_i^\eps(t,x)) = 0$ by definition of $F_i^{\eps}$ \eqref{def:F}. 

By \eqref{est:supp}, we infer the existence of $R_0>0$ such that $\frac{|x-y|^2}{x_ry_r}\leq R_0$ for all $x,y\in [z^*-1,z^*+1]\times[r^*/2,2r^*]$. Thanks to the expression of $\cal G$ in terms of $F$ \eqref{eq:F}, we deduce from Lemma~\ref{lem:F} that
\[
|\nabla \cal G(x,y) |\leq \frac{C}{|x-y|} ,\quad |\nabla^2 \cal G(x,y) |\leq \frac{C}{|x-y|^2},
\]
for some $C>0$ depending only on $X^*$.

Next we notice that for all $i\neq j$ and $(x,y)\in \mathcal{A}_{\eps}^{X^\eps_{i,0,r}}\left(d_0,h_0\right)\times \mathcal{A}_{\eps}^{X^\eps_{j,0,r}}\left(d_0,h_0\right) $, we have
\begin{equation*}
|x-y|\geq |x_r-y_r|\geq |X^\eps_{i,0,r}-X^\eps_{j,0,r}| - | x_r -X^\eps_{i,0,r}| - | y_r -X^\eps_{j,0,r}| \geq \frac{4d_0}{\sqrt{|\ln \eps|}} - \frac{2d_0}{\sqrt{|\ln \eps|}}.
\end{equation*}
 By the definition of $F_i$ \eqref{def:F} and the conservation of the mass \eqref{eq:consgamma}, this ends the proof.
\end{proof}

The remainder of the paper is then devoted to proving that a blob of vorticity $\omega_i^\eps$ which is transported by $u_i^\eps$, the vector field generated by $\omega_i^\eps$ through the Biot-Savart law, see Lemma~\ref{lem:decomp_u}, and by an exterior field $F_i^\eps$, which satisfies Lemma~\ref{lem:hyp_F}, remains concentrated in the sense of Theorem~\ref{theo:main} on $[0,T_\eps]$. More precisely, by choosing $T_0>0$ and $\eps_0>0$ small enough, we will show a stronger localization than the one given in the definition of $T_\eps$, i.e.
\begin{equation*}
 \supp \omega_i^\eps(t,\cdot) \subset \mathcal{A}_{\eps}^{X^\eps_{i,0,r}}\left(\frac{d_0}2,\frac{h_0}2\right), \qquad \forall t\in [0,T_\eps], \ \eps \in (0,\eps_0], \ i\in \{1,\ldots,N\}.
\end{equation*}
By definition of $T_\eps$ \eqref{def:Teps}, this will imply that $T_\eps=T_0$ and that the strong localization remains valid on $[0,T_0]$. \\
In this respect, it is crucial to verify that all the estimates below are uniform with respect to $i$. Since the proofs will focus only on one patch for a fixed $i\in \{1,\dots,N\}$, we will often drop the index $i$ for the sake of readability.

In the sequel, every constant $C,C_0,\eps_0$ may depend on $X^*$, $(Y_{i,0})$ and $M_0$ and we omit to mention this dependence. At the opposite, we will carefully track the dependence on $h_0$ and $T_0$. Finally, we consider from now on that $\gamma >0$, and the opposite sign will be reached simply by reversing the time. By \eqref{hyp:omega0}, this means that we are focusing only on non-negative vortices.

\section{Weak localization}\label{sec:weak_loc}

In this section, we focus on the weak localization which corresponds to Item~\textit{(i)} in Theorem~\ref{theo:main}. By the conservation of the total mass of vorticity \eqref{eq:consgamma}, we recall that $\int_{\R^2} x_r\omega^\eps_t(t,x)\, dx =\gamma$ for all $t\geq 0$. The key proposition here is the analog of \cite[Lemma 2.1]{Tur87}, also used in a different context in \cite{Guo_Zhao}, which states that almost all the vorticity is contained in small disks.

\begin{proposition}[Weak concentration]\label{prop:Turkington}
Given $T_{\eps}>0$ as in \eqref{def:Teps}, provided $T_0\leq 1/(24 C'_F)$ and $h_0\leq \min(r^*/240,1/2)$ the following holds true: for any positive function $\eps\mapsto \rho(\eps)$ such that $\rho(\eps)\sqrt{|\ln\varepsilon|}$ tends to zero as $\varepsilon\to 0$, there exists $\varepsilon_{\rho}>0$, depending only on $\rho$ and $h_0$, such that for any $\varepsilon\in (0,\varepsilon_{\rho}]$, every $i\in\{1,\ldots,N\}$ and all $t\in [0,T_\eps]$ there is an open set $\Omega_{\rho,i}^\eps(t)$ satisfying
\[ \frac1\gamma \int_{\Omega_{\rho,i}^\eps(t)} x_r \omega^\eps_i(t,x) \dd x \geq 1-\frac{h_0}{r^*}\rho(\eps)\ln|\ln\eps| \quad \text{and}\quad
 \mathrm{diam}\,\Omega_{\rho,i}^\eps(t) \leq \eps e^{120/\rho(\eps)}.
\] 
\end{proposition}

We do not claim that the constants are optimal, but we write them here explicitly in order to focus on the dependence with respect to $h_0$ and $T_0$.

Even if the proof was originally inspired by \cite{Tur87} in the context of the two dimensional Euler equations, we need to add several new arguments due to the anelastic constraint and the stronger singularity of our problem. Indeed, the vortex patches are moving by the binormal flow, highlighted by the presence of $u_L^\eps$ in Lemma~\ref{lem:decomp_u}, and are here close to each other, with a distance of order $1/\sqrt{|\ln\eps|}$. The strategy in this section is to first assume that weak localization holds up to some time $\tilde{T}_\eps$, which is less precise than the desired estimates. This will allow to obtain sharp energy estimates and to deduce a sharper weak localization property on $[0,\tilde{T}_\eps]$, which will imply that $\tilde{T}_\eps=T_\eps$ and that all the estimates proved up to $\tilde{T}_\eps$ hold actually true up to $T_\eps$. Among these estimates and as they will be used later, we state a lemma encompassed in this bootstrap argument, which will be established in the proof of Proposition~\ref{prop:Turkington}. This lemma involves the center of vorticity of each blob, which was defined at \eqref{def:beps}.

\begin{lemma}\label{lem:maj|x-b|}
Provided $T_0\leq 1/(24 C'_F)$ and $h_0\leq \min(r^*/240,1/2)$, there exists $\eps_0$ depending on $h_0$ such that for every $\eps \in (0,\eps_0]$, $i\in\{1,\ldots,N\}$ and $t \in [0,T_\eps]$,
 \begin{equation}\label{eq:maj|x-b|}
 \frac1\gamma \int_{\R^2_+} |x-b^\eps_i(t)| x_r\omega^\eps_i(t,x)\dd x \le 5h_0 \frac{\ln|\ln\eps|}{|\ln \eps|}.
 \end{equation}
\end{lemma}

Finally, at the end of this section, we will deduce from Proposition~\ref{prop:Turkington} the concentration property around the center of vorticity.

\begin{corollary}\label{cor:weak_loc} Given $T_{\eps}>0$ as in \eqref{def:Teps}, provided $T_0\leq 1/(24 C'_F)$ and $h_0\leq \min(r^*/240,1/2)$ the following holds true: there exists $\eps_0 > 0$ depending on $h_0$, such that for every $\eps \in (0,\eps_0]$ and $i\in\{1,\ldots,N\}$, 
 \begin{equation*}
 \sup_{t \in [0,T_\eps]}\left| \gamma - \int_{B\left( b^\eps_i(t) ,3h_0\frac{ \ln |\ln\eps|}{|\ln \eps|} \right)} x_r\omega^\eps_i(t,x)\dd x\right| \le \gamma \frac{\ln |\ln \eps|}{|\ln \eps|}.
 \end{equation*}
\end{corollary}

All these statements are proved within this section. 

\subsection{Weak localization bootstrapping time}

We begin by introducing a time $\tilde{T}_\eps$ such that the solution is weakly concentrated in a less precise way than in Proposition~\ref{prop:Turkington}. Namely, we set $\tilde{T}_\eps>0$ the largest time $\tilde{T}_\eps \le T_\eps$ such that for every $t \in [0,\tilde{T}_\eps]$ and every $i \in \{1,\ldots,N\}$, there exists a domain $\Omega_i^\eps(t)$ such that 
\begin{equation}\label{def:tildeTeps}
 \frac1\gamma \int_{\R^2_+\setminus\Omega_i^\eps(t)} x_r \omega_i^\eps(t,x) \dd x \leq \frac{\ln |\ln \eps|}{|\ln \eps|} \quad \text{and}\quad \diam \Omega_i^\eps(t) \leq \frac{1}{|\ln \eps|}.
\end{equation}
To check that $\tilde{T}_\eps>0$, we may set $\Omega_i^\eps(t)=X^\eps(t,\supp(\omega_{i,0}^\eps))$, so that
$$\int_{\R^2_+\setminus\Omega_i^\eps(t)} x_r \omega_i^\eps(t,x) \dd x=\int_{\R^2_+\setminus \supp(\omega_{i,0}^\eps)} x_r \omega_{i,0}(x) \dd x=0,$$ and observe that by virtue of \eqref{ineq:bound-u} and \eqref{dec:flow} we have $$\text{diam}(\Omega_i^\eps(t))\leq \text{diam}(\supp(\omega_{i,0}^\eps))+ C_\eps t\leq 2\eps +C_\eps t,$$ with $C_\eps$ possibly diverging as $\eps \to 0$, so that $\text{diam}(\Omega_i^\eps(t))<1/|\ln \eps|$ for $t>0$ sufficiently small. 

The different localization properties related to the definition of $T_\eps$ \eqref{def:Teps} and $\tilde{T}_\eps$ are illustrated in Figure~\ref{fig:all_times}.

\begin{figure}
 \centering

\begin{minipage}[c]{0.48\textwidth}
\centering
\begin{tikzpicture}
\draw[->] (-2,0) -- (5,0) node[below] {$z$};
\draw[->] (0,-0.5) -- (0,4) node[left] {$r$};

\draw[red, thick] (-1.5,1) -- (4,1);
\draw[red, thick] (-1.5,3) -- (4,3);
\draw[red, thick] (-1.5,1) -- (-1.5,3);
\draw[red, thick] (4,1) -- (4,3);

\begin{scope}
\clip (-1.5,1) rectangle (4,3);
\foreach \x in {-4,-3.75,...,6} {
 \draw[blue!20] (-2,\x) -- (5,\x+5);
}
\end{scope}

\filldraw (1.25,2) circle (2pt) node[below] {$X^\eps(0)$};
\draw[->,smooth] (1.25,2) .. controls (2.6,2.1) and (3,2.4) .. (3.25,2.55);
\filldraw (3.3,2.6) circle (2pt) node[left] {$X^\eps(t)$};

\draw[<->, dashed] (-1.5,3.3) -- (4,3.3) node[midway, above] {$2h_0$};
\draw[<->, dashed] (-1.8,1) -- (-1.8,3) node[midway, left] {$\frac{2d_0}{\sqrt{|\ln \epsilon|}}$};

\draw[red!80] (3,0.5) node {$t \le T_\eps$};
\draw (1.5,-1) node {Strong localization};
\end{tikzpicture}
\end{minipage}
\hfill
\begin{minipage}[c]{0.48\textwidth}
\centering
\begin{tikzpicture}
\filldraw (2,4) circle (2pt); 
\draw[->,smooth] (2,4) .. controls (3.2,4.2) and (4,4.8) .. (4.5,5.1);
\filldraw (4.55,5.15) circle (2pt);
\draw (6.3,5.15) node[right]{$\Omega^\eps(t)$};

\draw[blue!80, thick,dashed] (5.05,4.8) circle (1.3);
\draw[<->,dashed] (5.05,4.8) -- (5.05,6.1) node[above] {$\frac{1}{2|\ln \eps|}$};

\draw[red!80, thick] (2,4) circle (0.8);
\draw[<->, dashed] (2,4) -- (2,4.8) node[above] {$\eps$};

\draw[blue!80] (5.7,3.2) node {$t \le \tilde{T}_\eps$};
\draw[red!80] (3.1,3.2) node {$t =0$};

\draw[red!80,thick] (1.1,2.5) -- (7,2.5);
\draw[red!80,thick] (1.1,7) -- (7,7);

\draw[<->, dashed] (0.8,2.5) -- (0.8,7) node[midway, left] {$\frac{2d_0}{\sqrt{|\ln \epsilon|}}$};

\draw (4,1.8) node {Weak localization};
\end{tikzpicture}
\end{minipage}

 \caption{Strong and weak localization: on the left, the support of $\omega^\eps$ remains within a rectangle of height $2h_0$ and width $2d_0/\sqrt{|\ln \eps}$ for $t \le T_\eps$. On the right, most of the mass is concentrated within a ball of diameter $1/|\ln \eps|$ for $t \le \tilde{T}_\eps$, while initially all the mass was contained in a ball of radius $\eps$.}
 \label{fig:all_times}
\end{figure}

The first consequence of this definition of $\tilde{T}_\eps$ is that Lemma~\ref{lem:maj|x-b|} holds true up to $\tilde{T}_\eps$, namely we have the following:

\begin{lemma}\label{lem:maj|x-b|_vfaible}
 Provided $h_0\leq1$, there exists $\eps_0>0$ depending on $h_0$ such that for every $\eps \in (0,\eps_0]$, $i\in\{1,\ldots,N\}$, the inequality~\eqref{eq:maj|x-b|} holds true on the time interval $[0,\tilde{T}_\eps]$. In addition, it holds for all $[0,\tilde{T}_\eps]$ that
 \[
 \Omega_i^\eps(t) \subset B\left(b_i^\eps(t),\frac{3h_0\ln |\ln \eps|}{|\ln \eps|}\right).
 \]
 where $\Omega_i^\eps(t)$ is defined in \eqref{def:tildeTeps}.
\end{lemma}

\begin{proof}
Here and below, we frequently drop the index $i$ in the notation of $\omega_i^\eps$ and $\Omega_i^\eps$ during the proofs as they are performed for a fixed $i\in \{0,\ldots,N\}$ and we derive uniform estimates with respect to this parameter.

Let $x_0\in \Omega^\eps(t)\cap \supp(\omega^\eps(t,\cdot))$, hence $\Omega^\eps(t) \subset B\Big(x_0,\frac{1}{|\ln \eps|}\Big)$. By definition of $b^\eps$ \eqref{def:beps} and by the conservation of the total mass of the vorticity \eqref{eq:consgamma}, we write
\[
\gamma|b^\eps(t)-x_0| = \left| \int_{\R^2_+} (x-x_0) x_r \omega^\eps(t,x)\dd x\right| \le \int_{\R^2_+} |x-x_0| x_r \omega^\eps(t,x)\dd x .
\]
 We split the integration domain, then using that $\diam (\supp \omega^\eps(t)) \le 2h_0$ (see the definition of $T_\eps$ \eqref{def:Teps}) and the definition of $\tilde T_\eps$ \eqref{def:tildeTeps} we get that for $\eps$ small enough
 \begin{align*}
 \int_{\R^2_+} |x-x_0| x_r \omega^\eps(t,x)\dd x 
 & \le \int_{\Omega^\eps(t)} |x-x_0| x_r \omega^\eps(t,x)\dd x + \int_{\R^2_+ \setminus \Omega^\eps(t)} |x-x_0| x_r \omega^\eps(t,x)\dd x \\
 & \le \frac{\gamma}{|\ln \eps|} + 2h_0\gamma \frac{\ln| \ln \eps|}{|\ln \eps|},
 \end{align*}
 and thus
 \begin{align*}
 \int_{\R^2_+} |x-b^\eps(t)|x_r \omega^\eps(t,x)\dd x & \le \gamma|x_0-b^\eps(t)| + \int_{\R^2_+} |x-x_0|x_r \omega^\eps(t,x)\dd x \\
 & \le 2\gamma\left(\frac{1}{|\ln \eps|} + 2h_0 \frac{\ln| \ln \eps|}{|\ln \eps|}\right),
 \end{align*}
 which implies \eqref{eq:maj|x-b|} by choosing $\eps_0$ small enough such that $1\leq \frac{h_0}2 \ln| \ln \eps| $. 

 Another consequence with this choice of $\eps_0$ is that for any $x\in \Omega^\eps(t)$, 
 \[
 |x-b^\eps(t)|\leq |x-x_0| + |x_0-b^\eps(t)|\leq \frac{2}{|\ln \eps|} + 2h_0 \frac{\ln| \ln \eps|}{|\ln \eps|} \leq 3h_0 \frac{\ln| \ln \eps|}{|\ln \eps|}.
 \]
\end{proof}

\begin{remark}\label{rm:pourlem4.5}
 A direct consequence of Lemma~\ref{lem:maj|x-b|_vfaible} and definition \eqref{def:tildeTeps} of $\tilde{T}_\eps$ is that, for every $t \in [0,\tilde{T}_\eps]$ and every $i \in \{1,\ldots,N\}$, we have
 \begin{equation*}
 \int_{|x-b^\eps(t)| > \frac{3h_0\ln |\ln \eps|}{|\ln \eps|}} x_r \omega^\eps(t,x)\dd x \le \gamma\frac{\ln |\ln\eps|}{|\ln \eps|}.
\end{equation*}
\end{remark}

\subsection{Expansion of the energy}

Recalling the definition of $\psi_i^\eps$ given in \eqref{def:psi}, we first bound it by above and from below.

\begin{lemma}\label{lem:encadrement_psi}
Provided $h_0\leq 1/2$, there exists $\eps_0>0$ depending on $h_0$ such that for every $\eps \in (0,\eps_0]$, $i\in\{1,\ldots,N\}$, $t \in [0,\tilde{T}_\eps]$ and $x \in\supp \omega^\eps_i(t)$, it holds
\begin{equation*}
 0 \le -\psi^\eps_i(t,x) \le \frac{\gamma}{2\pi} |\ln\eps| \sqrt{x_r b_{i,r}^\eps(t)} + \frac{2\gamma}{\pi} h_0\ln |\ln\eps|.
\end{equation*}
 In particular, we get that for $\eps$ small enough,
 \begin{equation*}
 |\psi_i^\eps(t,x)| \le \frac{2\gamma}{\pi} |\ln \eps| r^*.
 \end{equation*}
\end{lemma}

\begin{proof}
 Since $t \le T_\eps$, we have $\diam \big(\supp \omega^\eps(t)\big) \le 2 h_0$, so provided that $h_0 < 1/2$, then
 \begin{equation*}
 \psi^\eps(t,x) \le 0\quad \text{if }x\in \supp \omega^\eps(t),
 \end{equation*} 
where we recall that we are focusing only on non-negative vortices.
 
 Next, we observe that
 \begin{align*}
 \psi^\eps(t,x) = &\frac{\sqrt{x_r}}{2\pi} \int \sqrt{b^\eps_r(t)}\ln|x-y| y_r \omega^\eps(t,y)\dd y \\
 & + \frac{\sqrt{x_r}}{2\pi} \int \left(\sqrt{y_r} - \sqrt{b^\eps_r(t)}\right)\ln|x-y| y_r \omega^\eps(t,y)\dd y \\
 := & A_1 + A_2.
 \end{align*}
 
 To estimate $A_1$, we apply the rearrangement of the mass (see Lemma~\ref{lem:rearrangement}) to $g(s)=-\ln s \mathds{1}_{(0,1)}$ and $M=(2r^*)M_0/\eps^2$, to get
 \begin{equation}\label{eq:rearrangement1}
 \begin{aligned}
 -\int\ln|x-y|y_r \omega^\eps(t,y)\dd y &\le -2\pi M \int_0^{\sqrt{\frac{\gamma}{\pi M }}} s \ln s \dd s = -2\pi M \Big(\frac12 s^2\ln s -\frac14 s^2 \Big)\Big|_{s=\sqrt{\frac{\gamma}{\pi M }}}\\
 &\leq \gamma |\ln \eps| + \mathcal{O}(1)
 \end{aligned}
\end{equation}
 and thus
 \begin{equation*}
 -A_1 \leq \frac{\gamma}{2\pi} \sqrt{x_r b_r^\eps(t)} |\ln \eps| + \mathcal{O}(1).
 \end{equation*}
 
 To estimate $A_2$, we recall that for $\eps$ small enough, every $x \in \supp \omega^\eps(t)$ satisfies $r^*/2\leq x_r \leq 2r^*$ which implies that $r^*/2\leq b_r^\eps(t) \leq 2r^*$, then 
 \begin{equation*}
 |A_2| = \left|\frac{\sqrt{x_r}}{2\pi} \int \frac{y_r - b^\eps_r(t)}{\sqrt{y_r} + \sqrt{b^\eps_r(t)}}\ln|x-y| y_r \omega^\eps(t,y)\dd y\right| 
 \le \frac{\sqrt{2r^*}}{2\pi} \int \frac{|y_r - b^\eps_r(t)| }{2\sqrt{r^*/2}}\big|\ln |x-y|\big|y_r \omega^\eps(t,y)\dd y.
 \end{equation*}
 We now split the integral 
 \begin{align*}
 |A_2| \le& \frac{1}{2\pi} \int_{|y-b^\eps(t)|\le \frac{3h_0 \ln |\ln \eps|}{|\ln \eps|}} |y_r - b^\eps_r(t)| \ln |x-y||y_r \omega^\eps(t,y)\dd y \\
 & + \frac{1}{2\pi} \int_{|y-b^\eps(t)|> \frac{3h_0 \ln |\ln \eps|}{|\ln \eps|}} |y_r - b^\eps_r(t)| \big| \ln |x-y|\big|y_r \omega^\eps(t,y)\dd y \\
 :=& A_{21} + A_{22}.
 \end{align*}
 We simply bound the first term 
 \begin{equation*}
 A_{21} \le \frac{3h_0}{2\pi} \frac{\ln |\ln \eps|}{|\ln \eps|} (\gamma|\ln \eps| + \mathcal{O}(1)) = \frac{3h_0\gamma}{2\pi}\ln |\ln \eps|+o(\ln |\ln \eps|).
 \end{equation*}
 For the second term, we use again the mass rearrangement lemma but this time on a function with a mass less than $\gamma\frac{\ln |\ln \eps|}{|\ln \eps|}$ (see Remark~\ref{rm:pourlem4.5}), to get that
 \begin{align*}
 A_{22} & \le \frac{1}{2\pi}\frac{2d_0}{\sqrt{|\ln\eps|}} \int_{|y-b^\eps(t)|> \frac{3h_0\ln |\ln \eps|}{|\ln \eps|}} -\ln |x-y|y_r \omega^\eps(t,y)\dd y \\
 & \le \frac{d_0}{\pi\sqrt{|\ln\eps|}} \frac{\gamma\ln |\ln \eps|}{|\ln\eps|} (|\ln\eps|+\cal O(1)) = o(\ln |\ln \eps|).
 \end{align*}
 This means that 
 \begin{equation*}
 |A_2| \le \frac{3h_0\gamma}{2\pi}\ln |\ln \eps|+o(\ln |\ln \eps|),
 \end{equation*}
 and thus that
 \begin{equation*}
 -\psi^\eps_i(t,x) \le \frac{\gamma}{2\pi} |\ln\eps| \sqrt{x_r b_{i,r}^\eps(t)} + \frac{2\gamma}{\pi} h_0\ln |\ln\eps|.
 \end{equation*}
\end{proof}

We now define the total energy related to the $i$-th vortex as
\begin{equation}\label{def:energy}
 \mathcal{E}^\eps_i(t) := -\iint_{\R^2_+\times\R^2_+}\mathcal{G}(x,y)x_r\omega^\eps_i(t,x)y_r\omega^\eps_i(t,y)\dd x \dd y.
\end{equation}
Due to the initial concentration, namely \eqref{hyp:omega0}, we have the following.

\begin{lemma}\label{lem:energy0}
 We have
 \begin{equation*}
 \mathcal{E}_i^\eps(0) = \frac{\gamma^2}{2\pi}|\ln \eps| b_{i,r}^\eps(0) + \mathcal{O}(1).
 \end{equation*}
\end{lemma}

\begin{proof}
Using the decomposition of the Green kernel given in Corollary~\ref{cor:expansion Green function}, we write
\begin{align*}
 \mathcal{E}^\eps(0)=& -\iint_{\R^2_+\times\R^2_+} \frac{\sqrt{x_r y_r}}{2\pi}\ln |x-y|x_r\omega^\eps_0(t,x)y_r\omega^\eps_0(t,y)\dd x \dd y + \mathcal{O}(1)\\
 =& -\frac{b_r^\eps(0) + \mathcal{O}(\eps)}{2\pi} \iint_{\R^2_+\times\R^2_+} \ln |x-y|x_r\omega^\eps_0(t,x)y_r\omega^\eps_0(t,y)\dd x \dd y + \mathcal{O}(1),
\end{align*}
where we have used that $\supp \omega^\eps_0\subset B(X^\eps_0,\eps)$.
Next we use the mass rearrangement as in \eqref{eq:rearrangement1} to state that on the one hand we have
\[
\iint_{\R^2_+\times\R^2_+} \ln |x-y|x_r\omega^\eps_0(t,x)y_r\omega^\eps_0(t,y)\dd x \dd y \geq -\gamma^2 |\ln\eps | + \mathcal{O}(1),
\]
whereas on the other hand, we use that $|x-y|\leq 2\eps$ for $x,y\in \supp \omega^\eps_0$ to get
\[
\iint_{\R^2_+\times\R^2_+} \ln |x-y|x_r\omega^\eps_0(t,x)y_r\omega^\eps_0(t,y)\dd x \dd y \leq \ln (2\eps)\gamma^2 =-\gamma^2 |\ln\eps | + \mathcal{O}(1)
\]
which ends the proof.
\end{proof}

We aim to precisely express $\mathcal{E}_i^\eps(t)$. To this end, we compute its time derivative.

\begin{proposition}\label{prop:der_energy}
Provided $h_0\leq 1/2$ and $T_0\leq 1/C_F'$, there exists $\eps_0>0$ depending on $h_0$ such that for every $\eps \in (0,\eps_0]$, $i\in\{1,\ldots,N\}$ and $t \in [0,\tilde{T}_\eps]$, we have that
 \begin{equation*}
 \left|\der{}{t} \mathcal{E}^\eps_i(t) - \frac{\gamma^2}{2\pi}F_{i,r}^\eps(t,b^\eps_i(t))\right| \le \frac{21 C'_F \gamma^2}{\pi} h_0 \ln |\ln \eps|,
 \end{equation*}
 where $C'_F$ is the lipschitz constant of $(F_i^\eps)_i$ (see Lemma~\ref{lem:hyp_F}).
\end{proposition}

\begin{proof}
By the same computations as the ones in the proofs of \cite[Lemma 4.10]{DonatiLacaveMiot} and \cite[Proposition 5.6]{HLM24}, using the symmetry of $\mathcal{G}$ and the fact that $\omega^\eps$ satisfies~\eqref{eq:Axi2D_rescaled_F} in the sense of Proposition~\ref{prop:WP}, we get
\begin{equation*}
 -\der{}{t} \mathcal{E}^\eps(t) = \frac{2}{|\ln \eps|} \iint (u^\eps+F^\eps)(t,x) \cdot \nabla_x \mathcal{G}(x,y) x_r\omega^\eps(t,x) y_r \omega^\eps(t,y) \dd x \dd y.
\end{equation*}
We now notice that
\begin{equation*}
 u^\eps(t,x) \cdot \int \nabla_x \mathcal{G}(x,y)y_r \omega^\eps(t,y) \dd y = u^\eps(t,x) \cdot \big(-x_r u^\eps(t,x) \big)^\perp = 0.
\end{equation*}
Therefore,
 \begin{align*}
 - \der{}{t} \mathcal{E}^\eps(t) 
 =& \frac{2}{|\ln \eps|}\iint F^\eps(t,x)\cdot\nabla_x \cal G(x,y) x_r\omega^\eps(t,x) y_r\omega^\eps(t,y)\dd x \dd y \\
 =& \frac{1}{\pi|\ln \eps|}\iint F^\eps(t,x) \cdot \sqrt{x_ry_r}\frac{x-y}{|x-y|^2}x_r\omega^\eps(t,x) y_r\omega^\eps(t,y)\dd x \dd y \\
 & + \frac{1}{2\pi|\ln \eps|}\iint F^\eps(t,x) \cdot e_r \sqrt\frac{y_r}{x_r}\ln|x-y|x_r\omega^\eps(t,x) y_r\omega^\eps(t,y)\dd x \dd y\\
 & + \frac{2}{|\ln \eps|}\iint F^\eps(t,x)\cdot\nabla_x S(x,y) x_r\omega^\eps(t,x) y_r\omega^\eps(t,y)\dd x \dd y \\
 :=& A_1+A_2+A_3.
\end{align*}
Since $\nabla S\in L^{\infty}([z^\ast-h_0,z^\ast+h_0]\times [r^\ast/2, 2r^\ast])$ by Corollary~\ref{cor:expansion Green function}, recalling \eqref{est:supp} and Lemma~\ref{lem:hyp_F} for $F_\eps$, we obtain $A_3 = \mathcal{O}\big(\frac{1}{\sqrt{|\ln \eps}|}\big)$.

For $A_1$ we use a standard symmetrization technique: we have by exchanging $x$ and $y$
\begin{align*}
 A_1&= \frac{1}{2\pi|\ln \eps|}\iint \left(F^\eps(t,x)-F^\eps(t,y)\right) \cdot \sqrt{x_ry_r}\frac{x-y}{|x-y|^2}x_r\omega^\eps(t,x) y_r\omega^\eps(t,y)\dd x \dd y \\
 |A_1| &\leq \frac{C_{F'}}{2\pi}\iint \sqrt{x_ry_r}x_r\omega^\eps(t,x) y_r\omega^\eps(t,y)\dd x \dd y\leq \frac{r^\ast C_{F'}\gamma^2}{\pi},
\end{align*}
where we have used Lemma~\ref{lem:hyp_F}.
 
We turn to $A_2$. We start by a naive estimate: bounding $F^\eps(t,x)$ by $C_F\sqrt{|\ln \eps|}$, we get that
\begin{equation*}
 |A_2| \le \frac{C_F}{\sqrt{|\ln \eps|}}\int \frac{1}{x_r} |\psi^\eps(t,x)|x_r\omega^\eps(t,x) \dd x
 \le \frac{2C_F}{r^*\sqrt{|\ln \eps|}}\int |\psi^\eps(t,x)|x_r\omega^\eps(t,x) \dd x
\end{equation*} 
then using Lemma~\ref{lem:encadrement_psi}, we get that
\begin{equation*}
 A_2 = \mathcal{O}(\sqrt{|\ln \eps|}),
\end{equation*}
and thus
\begin{equation}\label{est:badenergy}
 \der{}{t} \mathcal{E}^\eps(t) = \mathcal{O}(\sqrt{|\ln \eps|}).
\end{equation}
This of course is much less precise than what we are aiming to prove, but we use this information to improve the estimate of $A_2$.

Using Lemma~\ref{lem:energy0}, we get from \eqref{est:badenergy} that for every $t \le T_\eps\leq T_0\leq 1/C_F'$,
\begin{equation}\label{eq:E_naive}
 \mathcal{E}(t) = \frac{\gamma^2 b_r^\eps(0)}{2\pi}|\ln \eps| + \mathcal{O}(\sqrt{|\ln \eps|}).
\end{equation}
We now express $A_2$ as
\begin{align*}
 A_2 =& \frac{1}{2\pi |\ln \eps|} F_r^\eps(t,b^\eps(t)) \iint \sqrt\frac{y_r}{x_r} \ln |x-y |x_r\omega^\eps(t,x) y_r\omega^\eps(t,y)\dd x \dd y \\
 & + \frac{1}{2\pi |\ln \eps|} \iint \big( F_r^\eps(t,x)-F_r^\eps(t,b^\eps(t))\big)\sqrt\frac{y_r}{x_r}\ln|x-y|x_r\omega^\eps(t,x) y_r\omega^\eps(t,y)\dd x \dd y \\
 =& \frac{1}{ |\ln \eps| b_r^\eps(t)} F_r^\eps(t,b^\eps(t)) ( -\mathcal{E}^\eps(t) + \mathcal{O}(1)) \\
 & +\frac{1}{2\pi|\ln \eps|} F_r^\eps(t,b^\eps(t)) \iint \left(\frac{1}{x_r}-\frac{1}{b_r^\eps(t)} \right)\sqrt{x_r y_r}\ln|x-y|x_r\omega^\eps(t,x) y_r\omega^\eps(t,y)\dd x \dd y \\
 & + \frac{1}{2\pi |\ln \eps|} \iint \big( F_r^\eps(t,x)-F_r^\eps(t,b^\eps(t))\big)\sqrt\frac{y_r}{x_r}\ln|x-y|x_r\omega^\eps(t,x) y_r\omega^\eps(t,y)\dd x \dd y \\
 :=& A_{21} + A_{22} + A_{23}.
\end{align*}
Using successively the Lipschitz bound on $F^\eps$, the fact that $x_r > r^*/2$, Lemma~\ref{lem:encadrement_psi} then Lemma~\ref{lem:maj|x-b|_vfaible} we get that
\begin{align*}
 |A_{23}|& \le \frac{2C'_F}{r^*} \int |x-b^\eps(t)| |\psi^\eps(t,x)| x_r \omega^\eps(t,x)\dd x \\
 & \le \frac{2C'_F}{r^*} \frac{2\gamma}{\pi}r^*|\ln \eps| \int |x-b^\eps(t)| x_r \omega^\eps(t,x)\dd x \\
 & \le \frac{ 4 C'_F \gamma^2}{\pi} 5h_0\ln |\ln\eps|.
\end{align*}
On the other hand, using the $L^\infty$ bound on $F^\eps$, the same arguments as above yield
\begin{equation*}
 |A_{22}| =\mathcal{O}(1).
\end{equation*}
We now turn to $A_{21}$. Using the $L^\infty$ bound of $F^\eps$ and relation~\eqref{eq:E_naive}, we get that
\begin{align*}
 \left|A_{21} + \frac{\gamma^2}{2\pi} F_r^\eps(t,b^\eps(t))\right| & = \frac{|F_r^\eps(t,b^\eps(t))|}{|\ln \eps|b_r^\eps(t)} \left|-\mathcal{E}^\eps(t) + \frac{\gamma^2}{2\pi}|\ln \eps|b_r^\eps(t) +\mathcal{O}(1)\right| \\
 & \le \frac{2C_F \gamma^2}{r^* 2\pi\sqrt{|\ln \eps|}}\Big(|\ln \eps| |b_r^\eps(0)- b_r^\eps(t)|+\mathcal{O}(1)\Big)\\
 &= \mathcal{O}(1),
\end{align*}
where we have used that $b^\eps(0), b^\eps(t) \in \mathcal{A}_{\eps}^{X^\eps_{0,r}}\left(d_0,h_0\right)$ (see \eqref{def:Aeta} for the definition of this set), hence $|b_r^\eps(0)- b_r^\eps(t)|\leq 2d_0/\sqrt{|\ln\eps|}$. 

In conclusion, 
\begin{equation*}
 \left|\der{}{t} \mathcal{E}^\eps(t) - \frac{\gamma^2}{2\pi}F_r^\eps(t,b^\eps(t))\right| \le \frac{20 \gamma^2}{\pi} C'_F h_0\ln |\ln \eps| + \mathcal{O}(1),
\end{equation*}
which gives the desired result.
\end{proof}

Next, we relate the radial velocity of the center of mass with the field generated by the other vortex rings.

\begin{lemma}\label{lem:der_br}
Provided $h_0\leq 1/2$, there exists $\eps_0>0$ depending on $h_0$ such that for every $\eps \in (0,\eps_0]$, $i\in\{1,\ldots,N\}$ and $t \in [0,\tilde{T}_\eps]$, we have that
\begin{equation*}
 \left|\der{}{t} b_{i,r}^\eps(t) - \frac{1}{|\ln\eps|} F_{i,r}^\eps(t,b_i^\eps(t))\right| \le 6 C'_F h_0\frac{\ln |\ln \eps|}{|\ln \eps|}.
\end{equation*}
\end{lemma}

\begin{proof}
As for the computation of $\der{\cal E}{t}$ in the beginning of the proof of Proposition~\ref{prop:der_energy}, we compute from the definition of $b^\eps$ \eqref{def:beps} and the equations verified by $\omega^\eps$ \eqref{eq:Axi2D_rescaled_F} that
 \begin{align*}
 \der{}{t} b_r^\eps(t) 
 &=\frac{e_r}{\gamma |\ln\eps|} \cdot \int_{\R^2_+} \big(u^\eps(t,x) + F^\eps(t,x)\big) x_r\omega^\eps(t,x)\dd x .
\end{align*}
Using the decomposition of $u^\eps$ \eqref{def:u_K,u_L,u_R} and noticing that $u_L^\eps \cdot e_r = 0$, we get that
\begin{align*}
 \der{}{t} b_r^\eps(t) = &\frac{e_r}{\gamma |\ln\eps|}\cdot\int_{\R^2_+} u_K^\eps(t,x) x_r\omega^\eps(t,x)\dd x \\
 & +\frac{e_r}{\gamma |\ln\eps|}\cdot\int_{\R^2_+} u_R^\eps(t,x) x_r\omega^\eps(t,x)\dd x \\
 & +\frac{e_r}{\gamma |\ln\eps|}\cdot \int F^\eps(t,x) x_r\omega(t,x)\dd x \\
 :=& A_1 + A_2 + A_3.
\end{align*}
Since $u_R^\eps$ is uniformly bounded with respect to $\eps$ by Lemma~\ref{lem:decomp_u}, we first observe that $A_2 = \mathcal{O}\left(\frac{1}{|\ln \eps|}\right)$. Next, we compute using the expression of $u_K^\eps$ \eqref{def:u_K,u_L,u_R} that
\begin{equation*}
 A_1 = \frac{e_r}{2\pi \gamma|\ln\eps|}\cdot\iint \frac{\sqrt{y_r}}{\sqrt{x_r}} \frac{(x-y)^\perp}{|x-y|^2}x_r\omega^\eps(t,x)y_r\omega^\eps(t,y)\dd x \dd y
\end{equation*}
and symmetrize to get that $A_1 = \mathcal{O}\left(\frac{1}{|\ln \eps|}\right)$.

For the term $A_3$ we compute that
\begin{equation*}
 A_3 = \frac{1}{|\ln\eps|} F_r^\eps(t,b^\eps(t)) + \frac{1}{\gamma|\ln\eps|} \int \big(F_r^\eps(t,x)-F^\eps_r(t,b^\eps(t))\big) x_r\omega(t,x)\dd x 
\end{equation*}
then using the Lipschitz hypothesis on $F^\eps$ (see Lemma~\ref{lem:hyp_F}) with Lemma~\ref{lem:maj|x-b|_vfaible}, we get that 
\begin{equation*}
 \left|A_3 - \frac{1}{|\ln\eps|} F_r^\eps(t,b^\eps(t))\right| \le \frac{C'_F}\gamma\int |x-b^\eps(t)|x_r \omega^\eps(t,x)\dd x \le 5 C'_F h_0\frac{\ln |\ln \eps|}{|\ln \eps|}.
\end{equation*}
In conclusion,
\begin{equation*}
 \left|\der{}{t} b_r^\eps(t) - \frac{1}{|\ln\eps|} F_r^\eps(t,b^\eps(t))\right| \le 5 C'_F h_0\frac{\ln |\ln \eps|}{|\ln \eps|} +\mathcal{O}\left(\frac{1}{|\ln \eps|}\right),
\end{equation*}
which ends the proof.
\end{proof}

The conclusion of the section is the following.

\begin{corollary}\label{coro:energy}
Provided $h_0\leq 1/2$ and $T_0\leq 1/C_F'$, there exists $\eps_0>0$ depending on $h_0$ such that for every $\eps \in (0,\eps_0]$, $i\in\{1,\ldots,N\}$ and $t \in [0,\tilde{T}_\eps]$, we have that
 \begin{equation*}
 \left|\mathcal{E}_i^\eps(t)- \frac{\gamma^2}{2\pi}|\ln \eps| b_{i,r}^\eps(t)\right| \le \frac{(24 T_0 C'_F + \frac14) \gamma^2}{\pi} h_0\ln|\ln\eps|.
 \end{equation*}
\end{corollary}

\begin{proof}
Gathering the results of Proposition~\ref{prop:der_energy} and Lemma~\ref{lem:der_br}, we have that for every $t \le \tilde{T}_\eps$,
\begin{equation*}
 \left|\der{}{t} \left( \mathcal{E}^\eps(t) - \frac{\gamma^2}{2\pi}|\ln \eps| b^\eps_r(t) \right)\right| \le \frac{21C'_F \gamma^2}{\pi} h_0 \ln |\ln \eps| +\frac{3 C'_F \gamma^2}{\pi} h_0\ln |\ln \eps|.
\end{equation*}
 There simply remains to integrate in time this relation using Lemma~\ref{lem:energy0} to get that
 \begin{equation*}
 \left|\mathcal{E}^\eps(t)- \frac{\gamma^2}{2\pi}|\ln \eps| b_r^\eps(t)\right| \le \mathcal{O}(1) + \frac{24 C'_F \gamma^2 }{\pi} T_0 h_0\ln|\ln\eps|.
 \end{equation*}
\end{proof}

\subsection{Turkington's-type Lemma}

In this section we prove that Proposition~\ref{prop:Turkington} holds for any $t \le \tilde{T}_\eps$. More precisely, we prove the following.

\begin{proposition}\label{prop:Turkington_avant_bootstrap} Let $\eps\mapsto \rho(\eps)$ a positive function such that $\rho(\eps)\sqrt{|\ln\varepsilon|}\to 0$ as $\varepsilon\to 0$.
Provided $h_0\leq 1/2$ and $T_0\leq 1/C_F'$, there exist $\varepsilon_{\rho}>0$, depending only on $\rho$ and $h_0$, such that for any $\varepsilon\in (0,\varepsilon_{\rho}]$, every $i\in \{1,\ldots,N\}$ and all $t\in [0,\tilde{T}_\eps]$, there is an open set $\widetilde{\Omega}_{\rho,i}^\eps(t)$ satisfying
\begin{equation*}
 \int_{\R^2_+\setminus\widetilde{\Omega}_{\rho,i}^\eps(t)} x_r \omega^\eps_i(t,x) \dd x \le \gamma \frac{h_0}{r^*}\rho(\eps)\ln|\ln\eps| \quad \text{and}\quad
 \mathrm{diam}\,\widetilde{\Omega}_{\rho,i}^\eps(t) \leq \eps e^{B/\rho(\eps)}.
\end{equation*}
with
\begin{equation*}
 B = 20( 5+24 C'_F T_0)
\end{equation*}
which crucially does not depend on $\rho$.
\end{proposition}

\begin{proof}
For $A>0$ which will be determined later, let us define the open set
 \[
\widetilde{\Omega}_\rho^\eps(t) := \Big\{x\in \supp \omega^\eps(t,\cdot),\ - \psi^\eps(t,x) \ge \frac{\gamma}{2\pi} |\ln \eps| b_r^\eps(t) -\frac{A}{\rho(\eps)} \Big\}.
 \]

We have from Corollary~\ref{coro:energy} that
\[
- \frac{(24 T_0 C'_F+\frac14) \gamma^2}{\pi} h_0\ln|\ln\eps| \leq \mathcal{E}^\eps(t) - \frac{\gamma^2}{2\pi}|\ln\eps| b_r^\eps(t).
\] 
As
\[
\mathcal{E}^\eps(t) =- \int_{\R^2_+} \psi^\eps(t,x) x_r\omega^\eps(t,x) \dd x + \mathcal{O}(1),
\]
we obtain that 
\begin{align*}
 0 &\leq \mathcal{E}^\eps(t) - \frac{\gamma^2}{2\pi}|\ln\eps| b_r^\eps(t)+\frac{(24 T_0 C'_F+\frac14) \gamma^2}{\pi} h_0\ln|\ln\eps| \\
 &\leq \int_{\R^2_+} \Big( -\psi^\eps(t,x)-\frac{\gamma}{2\pi}|\ln\eps| b_r^\eps(t)\Big)x_r\omega^\eps(t,x) \dd x+\frac{(24 T_0 C'_F+\frac14) \gamma^2}{\pi} h_0\ln|\ln\eps|+ \mathcal{O}(1),
\end{align*}
hence
\begin{multline*}
 \int_{\R^2_+\setminus \widetilde{\Omega}_\rho^\eps(t)} \Big( \psi^\eps(t,x)+\frac{\gamma}{2\pi}|\ln\eps| b_r^\eps(t) \Big)x_r\omega^\eps(t,x) \dd x\\
\leq \int_{\widetilde{\Omega}_\rho^\eps(t)} \Big( -\psi^\eps(t,x)-\frac{\gamma}{2\pi}|\ln\eps| b_r^\eps(t) \Big)x_r\omega^\eps(t,x) \dd x +\frac{(24 T_0 C'_F+\frac12) \gamma^2}{\pi} h_0\ln|\ln\eps| .
\end{multline*}

Outside $\widetilde{\Omega}_\rho^\eps(t)$, we have by definition of this set:
\[
 \frac{A}{\rho(\eps)} \leq \psi^\eps(t,x) + \frac{\gamma}{2\pi}|\ln \eps|b_r^\eps(t).
\]
It follows that
\begin{align*}
\frac{A}{\rho(\eps)} \int_{\R^2_+\setminus \widetilde{\Omega}_\rho^\eps(t)} x_r\omega^\eps(t,x) & \dd x\\
\leq & \int_{\R^2_+\setminus \widetilde{\Omega}_\rho^\eps(t)}\Big( \psi^\eps(t,x)+\frac{\gamma}{2\pi}|\ln\eps| b_r^\eps(t) \Big)x_r\omega^\eps(t,x) \dd x \\
\leq& \int_{\widetilde{\Omega}_\rho^\eps(t)} \Big( -\psi^\eps(t,x)-\frac{\gamma}{2\pi}|\ln\eps| b_r^\eps(t) \Big)x_r\omega^\eps(t,x) \dd x +\frac{(24 T_0 C'_F+\frac12) \gamma^2}{\pi} h_0\ln|\ln\eps|\\
\leq & \int_{\widetilde{\Omega}_\rho^\eps(t)} \Big( -\psi^\eps(t,x) - \frac{\gamma}{2\pi} \sqrt{x_r b_r^\eps(t)}|\ln\eps| \Big)x_r\omega^\eps(t,x) \dd x\\
&+\frac{\gamma}{2\pi}|\ln \eps| \sqrt{b_r^\eps(t)}\int_{\widetilde{\Omega}_\rho^\eps(t)} \Big( \sqrt{x_r} - \sqrt{b_r^\eps(t)}\Big)x_r\omega^\eps(t,x) \dd x \\
& +\frac{(24 T_0 C'_F+\frac12) \gamma^2}{\pi} h_0\ln|\ln\eps|\\
\le & \frac{2\gamma^2}{\pi} h_0 \ln|\ln \eps| + \frac{\gamma}{2\pi}|\ln \eps|\frac{\sqrt{2r^*}}{2\sqrt{r^*/2}} \int_{\widetilde{\Omega}_\rho^\eps(t)}| x_r - b_r^\eps(t)|x_r\omega^\eps(t,x) \dd x \\
& + \frac{(24 T_0 C'_F+\frac12) \gamma^2}{\pi} h_0\ln|\ln\eps|\\
\leq & \frac{\gamma^2 }{\pi} h_0 ( 5+24C'_FT_0) \ln|\ln\eps|,
\end{align*}
where we used Lemma~\ref{lem:encadrement_psi} and Lemma~\ref{lem:maj|x-b|_vfaible}.
Choosing $A = \frac{1}{\pi}\gamma r^* ( 5+24C'_F T_0)$ gives the first estimate in the lemma for any $\eps \in (0,\eps_0]$.

Inside $\widetilde{\Omega}_\rho^\eps(t)$, we have
\begin{align*}
0\le &-\psi^\eps(t,x) - \frac{\gamma}{2\pi}|\ln\eps| b_r^\eps(t) + \frac{A}{\rho(\eps)} = -\frac{1}{2\pi}\int_{\R^2_+} \big(\sqrt{x_ry_r}\ln|x-y| - (\ln\eps )b_r^\eps(t) \big)y_r \omega^\eps(t,y) \dd y+ \frac{A}{\rho(\eps)} \\
\le &-\frac{b_r^\eps(t)}{2\pi} \int_{\R^2_+} \ln\left(\frac{|x-y|}{\eps}\right) y_r\omega^\eps(t,y) \dd y
-\frac{1}{2\pi}\int_{\R^2_+} (\sqrt{x_ry_r} -b_r^\eps(t) )\ln|x-y| y_r\omega^\eps(t,y) \dd y+ \frac{A}{\rho(\eps)} \\
\le &-\frac{b_r^\eps(t)}{2\pi} \int_{\R^2_+} \ln\left(\frac{|x-y|}{\eps}\right) y_r\omega^\eps(t,y) \dd y +C\sqrt{|\ln\eps|} + \frac{A}{\rho(\eps)},
\end{align*}
where we have used $\sqrt{x_ry_r} -b_{r}^\eps(t) = \mathcal{O}(1/\sqrt{|\ln \varepsilon|})$ for any $x,y\in \supp \omega^\eps(t)$ and \eqref{eq:rearrangement1}.
Therefore, for $R\geq 1$
\begin{align*}
\int_{\R^2_+\setminus B(x,R\varepsilon)} \ln\left(\frac{|x-y|}{\eps}\right) y_r\omega^\eps(t,y) \dd y 
\leq&
-\int_{\R^2_+\cap B(x,R\varepsilon)} \ln\left(\frac{|x-y|}{\eps}\right) y_r\omega^\eps(t,y) \dd y +C\sqrt{|\ln\eps|} + \frac{2\pi A}{\rho(\eps)b_r^\eps(t)}\\
\leq&
-\int_{\R^2_+\cap B(x,\varepsilon)} \ln\left(\frac{|x-y|}{\eps}\right) y_r\omega^\eps(t,y) \dd y +C\sqrt{|\ln\eps|} + \frac{4\pi A}{\rho(\eps)r^*}
\end{align*}
because $-\ln\frac{|x-y|}{\eps} \leq 0 $ outside $B(x,\varepsilon)$.
By the mass rearrangement Lemma~\ref{lem:rearrangement} applied to $g(s)=-\ln\frac{s}{\eps}\mathds{1}_{(0,\eps)}$ and $M=2r^*M_0/\eps^2$, we get
 \begin{align*}
 -\int_{\R^2_+\cap B(x,\varepsilon)} \ln\left(\frac{|x-y|}{\eps}\right) y_r \omega^\eps(t,y)\dd y &\le -2\pi M \int_0^{\min(\eps,\sqrt{\frac{\gamma}{\pi M }})} s \ln\frac{s}{\eps} \dd s \\
 &\leq -2\pi M \Big(\frac12 s^2\ln \frac{s}{\eps} -\frac14 s^2 \Big)\Big|_{s=\min(\eps,\sqrt{\frac{\gamma}{\pi M }})}
 = \mathcal{O}(1).
 \end{align*}
 This implies that
 \[
 \ln R\int_{\R^2_+\setminus B(x,R\varepsilon)} y_r\omega^\eps(t,y) \dd y \leq C\sqrt{|\ln\eps|} + \frac{4\pi A}{\rho(\eps)r^*}.
 \]
Choosing $B = A \frac{20\pi }{\gamma r^*} $ and $R=\frac12e^{B/\rho(\eps)}>1$, we have for $\eps$ small enough (depending on $\rho$ and $A$, hence on $T_0$) that
\[
\int_{\R^2_+\setminus B(x,R\varepsilon)} y_r\omega^\eps(t,y) \dd y \leq \frac{4\pi}{r^*}\frac{C\rho(\eps) \sqrt{|\ln\eps|} + A}{B-\rho(\eps)\ln 2}\leq \frac{5\pi}{r^*}\frac{A}{B} = \frac14 \gamma
\]
because $\rho(\eps) \sqrt{|\ln\eps|}\to 0$ as $\eps\to 0$. As this inequality holds true for any $x\in \widetilde{\Omega}_\rho^\eps(t)$, it follows that
\[
\mathrm{diam}\, \widetilde{\Omega}_\rho^\eps(t) \leq 2R\varepsilon = \eps e^{B/\rho(\eps)}.
\]
Indeed, if we have $x,\tilde x \in\widetilde{\Omega}_\rho^\eps(t)$ such that $|x-\tilde x|>2R\varepsilon$, it would mean that
\[ 
\int_{B(x,R\varepsilon)\cup B(\tilde x,R\varepsilon)} y_r\omega^\eps(t,y) \dd y = \int_{B(x,R\varepsilon)} y_r\omega^\eps(t,y) \dd y + \int_{ B(\tilde x,R\varepsilon)} y_r\omega^\eps(t,y) \dd y \geq 2\frac34 \gamma
\]
which is impossible by the assumption $\int y_r\omega^\eps(y) \dd y=\gamma$ and $\omega^\eps\geq 0$.
\end{proof}

\subsection{Bootstrap argument and ends of the proof of the main results of this section}\label{sec:bootstrapweak}

In this paragraph we finally establish that $\tilde{T}_\eps = T_\eps$, by proving that for every $t \le \tilde{T}_\eps$ and $i\in \{1,\ldots,N\}$, there exists a domain $\widetilde{\Omega}_i^\eps(t)$ such that
\begin{equation*}
 \diam \widetilde{\Omega}_i^\eps(t) \le \frac{1}{2|\ln \eps|},
\end{equation*} and
\begin{equation*}
 \frac1\gamma\int_{\R^2_+\setminus\widetilde{\Omega}_i^\eps(t)} x_r \omega^\eps_i(t,x)\dd x \le \frac{3}{4}\frac{\ln |\ln \eps|}{|\ln \eps|}.
\end{equation*}
By continuity of the trajectories, this will prevent the inequality $\tilde{T}_\eps < T_\eps$ from being possible.

By Proposition~\ref{prop:Turkington_avant_bootstrap} applied to $\rho(\eps) = \frac{B}{|\ln \eps|-2\ln |\ln \eps|}$, there exists a set $\widetilde{\Omega}_i^\eps(t)$ such that
\begin{equation*}
\mathrm{diam}\,\widetilde{\Omega}_i^\eps(t) \leq \eps e^{B/\rho(\eps)} = \frac{1}{|\ln \eps|^2} \le \frac{1}{2|\ln \eps|}
\end{equation*}
for $\eps_0$ small enough, and
\begin{equation*}
 \frac{1}{\gamma}\int_{\R^2_+\setminus\widetilde{\Omega}_i^\eps(t)} x_r \omega^\eps(t,x) \dd x \le \frac{h_0}{r^*}\rho(\eps) \ln |\ln \eps| = \frac{h_0}{r^*}\frac{B\ln|\ln\eps|}{|\ln \eps|- 2\ln |\ln \eps|} \le \frac{3}{4}\frac{\ln |\ln \eps|}{|\ln \eps|},
\end{equation*}
provided that
\begin{equation}\label{ine:cond-h0}
\frac{Bh_0}{r^*} = 20(5+24 C_F' T_0)\frac{h_0}{r*} \leq \frac 12 .
\end{equation}
The inequality \eqref{ine:cond-h0} entails the constraint on $h_0$ and $T_0$ in the main results stated at the beginning of this section. Assuming for instance $T_0 \leq 1/(24 C_F')$ and $h_0 \leq r^*/240$, all the results of this section read as follows: there exists $\eps_0>0$ depending only on $h_0$ such that for any $\eps\in (0,\eps_0]$ and every $i\in \{1,\ldots,N\}$, all the inequalities stated hold true for all $t\in [0,T_\eps]$. In particular, Proposition~\ref{prop:Turkington} and respectively Lemma~\ref{lem:maj|x-b|} follow directly from Proposition~\ref{prop:Turkington_avant_bootstrap} and resp. Lemma~\ref{lem:maj|x-b|_vfaible}. Finally, Corollary~\ref{cor:weak_loc} is a direct consequence of Lemma~\ref{lem:maj|x-b|_vfaible} combined with the definition of $\Omega_i^\eps$ \eqref{def:tildeTeps}, see Remark~\ref{rm:pourlem4.5}.

\begin{remark}\label{rem:h0}
 The constant $h_0$ was chosen small enough so that some constants are absorbed, see \eqref{ine:cond-h0}, which explained why we have tracked in all this section the dependence on $h_0$. The possibility to choose $h_0$ small will not be used below, therefore, $h_0$ is fixed as $h_0=\min(r^*/240,1/2)$ and we will not specify the dependence of the constants on $h_0$ in the sequel.
\end{remark}

\begin{remark}\label{rem:Turkington}
Before going to the next part of the paper, let us remark that Proposition~\ref{prop:Turkington} can give us a much sharper concentration result than Corollary~\ref{cor:weak_loc}. Indeed, by taking for instance $\rho(\eps) = \frac{120}{\eta |\ln \eps|}$, for any given $\eta > 0$, we get the existence of a domain $\hat \Omega^\eps_i(t)$ such that
\[ \frac1\gamma \int_{\hat \Omega_i^\eps(t)} x_r \omega^\eps_i(t,x) \dd x \geq 1-\frac{h_0}{r^*}\frac{120}{\eta}\frac{\ln|\ln\eps|}{|\ln \eps|} \quad \text{and}\quad
 \mathrm{diam}\,\hat\Omega^\eps(t) \leq \eps^{1-\eta},
\] 
which means that most of the mass can be found in an area of radius $\eps^{1-\eta}$. However, this comes at the price of a relatively rough information on the location of $\hat\Omega^\eps_i(t)$. Specifically, with the localization estimates at hand and the mass being located away from $\hat\Omega^\eps_i(t)$ of order $\mathcal{O}\left(\frac{\ln|\ln\eps|}{|\ln \eps|}\right)$, we can only show that \begin{equation*}
 \dist(b^\eps(t), \hat{\Omega}^\eps(t)) = \mathcal{O}\left(\frac{\ln|\ln\eps|}{|\ln \eps|}\right),
\end{equation*}
by adapting the proof of Lemma~\ref{lem:maj|x-b|_vfaible} and no better than that. Therefore, the center of mass can be quite far from the actual location of the vortex filament compared to the size of filament's core. This is due to the possible filamentation of the vortex filament as discussed in the introduction.
 \end{remark}

\section{Equation of the motion}\label{sec:motion}

Lemma~\ref{lem:der_br} provides the equation of the motion of the radial component of $b^\eps_i(t)$ that holds up to $T_\eps$ as $\tilde{T}_\eps=T_\eps$ being defined in \eqref{def:tildeTeps}, \eqref{def:Teps} respectively. The goal of this section is to derive the full equation for both components and then establish the link with $(X_i^\eps)_{1\le i \le N}$ as in \eqref{def:X_i}.

\subsection{Motion of the center of mass}

We prove in the following lemma that the motion of the center of mass of the vortex ring is governed by the binormal flow together with the influence of the other vortex rings.

\begin{lemma}\label{lem:eq_b}
Provided $T_0\leq 1/(24 C'_F)$, there exists $\eps_0 > 0$ and $C>0$ independent of $T_0$, such that for every $\eps \in (0,\eps_0]$, $i\in\{1,\ldots,N\}$ and $t \in [0,T_\eps]$, we have that
\begin{equation*}
 \left|\der{}{t} b_{i}^\eps(t) - \frac{1}{|\ln\eps|} F_{i}^\eps(t,b_i^\eps(t))- \gamma\frac{e_z}{4\pi b_{i,r}^\eps(t)} \right| \le C \frac{\ln |\ln \eps|}{|\ln \eps|}.
\end{equation*}
\end{lemma}

\begin{proof}
As in the beginning of the proof of Lemma~\ref{lem:der_br}, we start by writing
\begin{equation*}
 \der{}{t} b^\eps(t) = \frac{1}{\gamma|\ln\eps|}\int_{\R^2_+} \big(u^\eps(t,x) + F^\eps(t,x)\big) x_r\omega^\eps(t,x)\dd x.
\end{equation*}
A few computations were already done in Lemma~\ref{lem:der_br}, and we obtain that
\begin{align*}
 \der{}{t} b^\eps(t) =& \frac{1}{\gamma|\ln\eps|}\int_{\R^2_+} u_K^\eps(t,x) x_r\omega^\eps(t,x)\dd x \\
 & + \frac{1}{\gamma|\ln\eps|}\int_{\R^2_+} u_L^\eps(t,x) x_r\omega^\eps(t,x)\dd x \\
 & + \frac{1}{\gamma|\ln\eps|} \int u_R^\eps(t,x) x_r\omega(t,x)\dd x\\
 & + \frac{1}{\gamma|\ln\eps|}\int_{\R^2_+} F^\eps(t,x) x_r\omega^\eps(t,x)\dd x \\
 := & A_1 + A_2 + A_3 + A_4,
\end{align*}
with $A_1 = \mathcal{O}\left(\frac{1}{|\ln \eps|}\right)$ by symmetrizing, and $A_3 = \mathcal{O}\left(\frac{1}{|\ln \eps|}\right)$.

For the term $A_4$ we compute that
\begin{equation*}
 A_4 = \frac{1}{|\ln\eps|} F^\eps(t,b^\eps(t)) + \frac{1}{\gamma|\ln\eps|} \int \big(F^\eps(t,x)-F^\eps(t,b^\eps(t))\big) x_r\omega(t,x)\dd x 
\end{equation*}
then using from Lemma~\ref{lem:hyp_F} that $F^\eps$ is Lipschitz and Lemma~\ref{lem:maj|x-b|}
\begin{equation*}
 \left|A_4 - \frac{1}{|\ln\eps|} F^\eps(t,b^\eps(t))\right| \le \frac{C'_F}\gamma \int |x-b^\eps(t)| x_r \omega(t,x) \le 5 C'_Fh_0\frac{\ln |\ln \eps|}{|\ln \eps|}.
\end{equation*}
We now turn to $A_2$. Using the definition of $u_L$ \eqref{def:u_K,u_L,u_R}, we have
\begin{align*}
 A_2 = &-\frac{e_z}{2\gamma|\ln\eps|}\int \frac{1}{x_r^2}\psi^\eps(t,x) x_r\omega^\eps(t,x)\dd x \\
 =& - \frac{e_z}{2\gamma|\ln\eps||b_r^\eps(t)|^2}\int_{\R^2_+} \psi^\eps(t,x) x_r\omega^\eps(t,x)\dd x \\
 & - \frac{e_z}{2\gamma|\ln\eps|} \int_{\R^2_+} \left(\frac{1}{x_r^2}-\frac{1}{|b_r^\eps(t)|^2} \right)\psi^\eps(t,x)x_r\omega^\eps(t,x)\dd x.
\end{align*}
By Lemma~\ref{lem:encadrement_psi} and Lemma~\ref{lem:maj|x-b|}, we have
\begin{align*}
 \int_{\R^2_+} \left|\frac{1}{x_r^2}-\frac{1}{|b_r^\eps(t)|^2} \right||\psi^\eps(t,x)|x_r\omega^\eps(t,x)\dd x & \le C |\ln\eps| \int_{\R^2_+} \big|x_r - b_r^\eps(t)| x_r \omega^\eps(t,x)\dd x\\
 & \le C \ln |\ln \eps|.
\end{align*}
Next we use the definition of the energy \eqref{def:energy}, the expression of $\psi^\eps$ \eqref{def:psi}, the decomposition of the kernel provided in Corollary~\ref{cor:expansion Green function}, with the uniform bound of $S$ in $L^{\infty}_{\loc}$ and the support of $\omega^\eps(t)$ satisfying \eqref{est:supp}, to write
\begin{equation*}
 -\int_{\R^2_+} \psi^\eps(t,x) x_r \omega^\eps(t,x)\dd x = \mathcal{E}^\eps(t) + \mathcal{O}(1),
\end{equation*}
then, by using Corollary~\ref{coro:energy}, we obtain that
\begin{equation*}
 \left|A_2 -
 \gamma\frac{e_z}{4\pi b_r^\eps(t)}\right| \le C(1+ T_0) \frac{\ln|\ln\eps|}{\ln \eps}\le C \frac{\ln|\ln\eps|}{\ln \eps},
\end{equation*}
provided $T_0\leq 1/(24 C'_F)$.
\end{proof}

Finally, we show that the velocity $F_{i}^\eps$ created by the other vortices behaves as the interaction in the classical point vortex system (see \cite{MarPul93,MP94}).

\begin{lemma}\label{lem:pour_gronwall}
 Provided $T_0\leq 1/(24 C'_F)$, there exists $\eps_0 > 0$ and $C>0$ independent of $T_0$, such that for every $\eps \in (0,\eps_0]$, $i\in\{1,\ldots,N\}$ and $t \in [0,T_\eps]$, we have 
 \begin{equation*}
 \left| F_i^\eps(t,b_i^\eps(t)) - \frac{\gamma}{2\pi} \sum_{j \neq i} \frac{\big(b_i^\eps(t)-b_j^\eps(t)\big)^\perp}{\big|b_i^\eps(t)-b_j^\eps(t)\big|^2}\right| \le C \ln |\ln \eps|.
 \end{equation*}
\end{lemma}

\begin{proof}
 Going back to the definition of the $F_i^\eps$ \eqref{def:F} and the decomposition of $\mathcal{G}$, see Corollary~\ref{cor:expansion Green function}, we have that
 \begin{align*}
 F_i^\eps(t,b_i^\eps(t)) =& \frac{1}{b_{i,r}^\eps(t)}\sum_{j \neq i} \int_{\R^2_+} \nabla^\perp_x \mathcal{G}(b_i^\eps(t),y) y_r \omega_j^\eps(t,y)\dd y \\
 = &\frac{1}{2\pi b_{i,r}^\eps(t)}\sum_{j \neq i} \int_{\R^2_+} \sqrt{y_r b_{i,r}^\eps(t)}\frac{\big(b_i^\eps(t)-y\big)^\perp}{\big|b_i^\eps(t)-y\big|^2} y_r \omega_j^\eps(t,y)\dd y \\
 & - \frac{e_{z}}{4\pi b_{i,r}^\eps(t)}\sum_{j \neq i} \int_{\R^2_+} \sqrt\frac{y_r}{b_{i,r}^\eps(t)}\ln |b_i^\eps(t)-y| y_r \omega_j^\eps(t,y)\dd y \\
 & + \frac{1}{b_{i,r}^\eps(t)}\sum_{j \neq i} \int_{\R^2_+} \nabla^\perp_{x} S(b_i^\eps(t),y) y_r \omega_j^\eps(t,y)\dd y \\
 :=& A_1 + A_2 + A_3.
 \end{align*}
 By the localization of the vorticity supports entailed in the definition of $T_{\eps}$ \eqref{def:Teps}, we know that $|b_{i,r}^\eps(t)-X_{i,0,r}|\leq d_0/\sqrt{|\ln \eps|}$. Hence for $j\neq i$ and $y\in \supp \omega_j^\eps(t)$ we have by definition of $d_0$ in \eqref{def:min-dist} that
 $$|b_i^\eps(t)-y| \ge |X_{i,0,r}-X_{j,0,r}|-|X_{i,0,r}-b_{i,r}^\eps(t)|-|y-X_{j,0,r}|\geq \frac{2d_0}{\sqrt{|\ln \eps|}}.$$
 It follows that $A_2 = \mathcal{O}(\ln |\ln \eps|)$. On the other hand, we infer from Corollary~\ref{cor:expansion Green function} that $A_3 = \mathcal{O}(1)$. 
 
 We now turn to $A_1$. By definition of $T_\eps$ \eqref{def:Teps}, we have that for $y\in \supp \omega_{j}^\eps$ and $j\neq i$
 \[
\frac{|\sqrt{b_{i,r}^\eps}-\sqrt{y_{r}}|}{|b_i^\eps(t)-y|} =\frac{|b_{i,r}^\eps-y_{r}|}{(\sqrt{b_{i,r}^\eps}+\sqrt{y_{r}})|b_i^\eps(t)-y|}\leq \frac1{r^*}
 \]
 hence
\[
 A_1 
 = \frac{1}{2\pi}\sum_{j \neq i} \int_{\R^2_+} \frac{\big(b_i^\eps(t)-y\big)^\perp}{\big|b_i^\eps(t)-y\big|^2} y_r \omega_j^\eps(t,y)\dd y + \mathcal{O}(1).
\]
Moreover, since $|b_i^\eps-b_j^\eps|\geq |X_{i,0,r}-X_{j,0,r}|-2d_0/\sqrt{|\ln \eps|}\geq 2d_0/\sqrt{|\ln \eps|}$, 
using that $a/|a|^2-b/|b|^2=(a-b)/(|a||b|)$, we deduce from Lemma~\ref{lem:maj|x-b|} that
 \begin{align*}
 \left|A_1 - \frac{1}{2\pi}\sum_{j \neq i} \int_{\R^2_+} \frac{\big(b_i^\eps-b_j^\eps\big)^\perp}{\big|b_i^\eps-b_j^\eps\big|^2} y_r \omega_j^\eps(t,y)\dd y \right| &\le \frac1{2\pi} \sum_{j \neq i}\int_{\R^2_+} \frac{|b_j^\eps-y|}{|b_i^\eps-b_j^\eps||b_i^\eps-y|} y_r \omega_j^\eps(t,y) \dd y + \mathcal{O}(1)\\
 &\le C h_{0}|\ln \eps| \frac{\ln |\ln \eps|}{|\ln \eps|} \leq C\ln |\ln \eps|.
 \end{align*}
 Gathering all the previous estimates concludes the proof.
 \end{proof}

\subsection{Limiting trajectories}\label{sec:limiting_traj}

We prove in this section that the center of vorticities $(b_{i}^\eps)_{1\leq i\leq N}$, defined in \eqref{def:beps}, are located close to $(X_{1\leq i\leq N}^\eps)_{i}$ and $(\tilde X_{i}^\eps)_{1\leq i\leq N}$, defined in \eqref{def:X_i} and \eqref{def:tildeX_i} respectively. 

First, we prove that $(X_{i}^\eps)_{i}$ and $(\tilde X_{i}^\eps)_{i}$ are located in $\cal{A}_\eps^{X^\eps_{i,0,r}}(d_0,h_0)$, see \eqref{def:Aeta} for the definition of this set, and remain suitably close. 

\begin{lemma}\label{lem:X-tildeX}
 For $h_0$ fixed, there exist positive numbers $T_X$, $\eps_0$ and $C$, such that for every $t\in [0,T_X]$, $\eps\in (0,\eps_0]$ and $i\in \{1,\ldots, N\}$, we have
 \[
 X_i^\eps(t),\tilde X_i^\eps(t) \in \cal{A}_\eps^{X^\eps_{i,0,r}}(d_0,h_0) \quad \text{and} \quad |X_i^\eps(t) - \tilde{X}_i(t)| \le \frac{C}{|\ln \eps|}.
 \]
\end{lemma}

\begin{proof}
First, we consider $\tilde{Y}_i$ defined in \eqref{eq:LP}. Given the independence of $\eps$ of the respective system, it is clear by continuity of the trajectories that there exists $T_X\in (0,1]$ such that
\[
\tilde X_i^\eps(t) \in \cal{A}_\eps^{X^\eps_{i,0,r}}\Big(\frac{d_0}2,\frac{h_0}2\Big), \quad \forall t\in [0,T_X].
\]

Next, we remark that by definition $Y_i^\eps$ is given by
\begin{equation*}
 Y_i^\eps(t) = \sqrt{|\ln \eps|} \left( X_i^\eps(t) - X^* - \frac{\gamma}{4\pi r^*} t e_z\right),
\end{equation*}
which satisfies that $Y_i^\eps(0) = Y_{i,0}$ and
\begin{align*}
 \der{}{t}Y_i^\eps(t) & = \frac{\gamma}{2\pi} \sum_{j\neq i} \frac{(Y_i^\eps-Y_j^\eps)^\perp}{|Y_i^\eps-Y_j^\eps|^2} + \sqrt{|\ln \eps|}\frac{\gamma}{4\pi r^*} \left( \frac{1}{1+\frac{Y_{i,r}^\eps}{r^*\sqrt{|\ln \eps|}}}- 1 \right) e_z.
\end{align*}
By continuity of the trajectories, we state that the exists $\widehat T_\eps\in (0,T_X]$ such that $X_i^\eps(t) \in \cal{A}_\eps^{X^\eps_{i,0,r}}\Big(d_0,h_0\Big)$ for all $t\in [0,\widehat T_\eps]$.

Noticing that
\begin{equation*}
 -\frac{\gamma}{4\pi (r^*)^2}\tilde{Y}_{i,r} = \sqrt{|\ln \eps|}\frac{\gamma}{4\pi r^*}\left( \frac{1}{1+\frac{\tilde{Y}_{i,r}}{r^*\sqrt{|\ln \eps|}}}- 1 \right) + \mathcal{O}\left( \frac{1}{\sqrt{|\ln \eps|}}\right),
\end{equation*}
we conclude that for all $t\in [0,\widehat T_\eps]$
\begin{align*}
 & \der{}{t} \cal Y^\eps = \mathcal{F}_1^\eps(\cal Y^\eps) \\
 & \left| \der{}{t} \tilde{\cal Y} - \mathcal{F}_1^\eps(\tilde{\cal Y}) \right| = \mathcal{O}\left( \frac{1}{\sqrt{|\ln \eps|}}\right),
\end{align*}
where
 \begin{equation*}
 \mathcal{F}^\eps_{1} : \begin{cases}
 (\R^2_+)^N \to (\R^2)^N \\
 \displaystyle \mathcal{Y} = (Y_1,\ldots,Y_N) \mapsto \left( \frac{\gamma}{2\pi} \sum_{j\neq i} \frac{(Y_i-Y_j)^\perp}{|Y_i-Y_j|^2} + \sqrt{|\ln \eps|}\frac{\gamma}{4\pi r^*} \left( \frac{1}{1+\frac{Y_{i,r}}{r^*\sqrt{|\ln \eps|}}}- 1 \right) e_z\right)_{1 \le i \le N}.
 \end{cases}
 \end{equation*}
The map $\mathcal{F}^\eps_{1}$ is Lipschitz continuous when $\cal X^\eps,\tilde{\cal X}^\eps \in\prod_{i} \mathcal{A}_{\eps}^{X^\eps_{i,0,r}}\left(d_0,h_0\right)$ with constant $K = \mathcal{O}(1)$, so by Gronwall lemma (see for instance the variant in Lemma~\ref{lem:gronwall}), we get that for every $t \le \widehat{T}_\eps$,
 \begin{equation*}
 \left|\mathcal{Y}^\eps(t) - \tilde{\mathcal{Y}}(t)\right| \le \left( \frac{CT_X}{\sqrt{|\ln \eps|}} + 0 \right) e^{KT_X},
 \end{equation*}
for all $t\in [0,\widehat T_\eps]$.
 
In conclusion there exists a constant $C$ depending on $(Y_{i ,0})_i$ but not on $\eps$ such that for any $t \in [0,\widehat T_\eps]$,
\begin{equation*}
 |Y_i^\eps(t) - \tilde{Y}_i(t)| \le \frac{C}{\sqrt{|\ln \eps|}},
\end{equation*}
and thus
\begin{equation*}
 |X_i^\eps(t) - \tilde{X}_i(t)| \le \frac{C}{|\ln \eps|}.
\end{equation*}
Choosing $\eps_0$ small enough, this allows us to state that 
\[
X_i^\eps(t) \in \cal{A}_\eps^{X^\eps_{i,0,r}}\Big(\frac{3d_0}4,\frac{3h_0}4\Big),
\]
and hence no blow-up happens on $[0,\hat{T}_\eps]$, and by continuity of the trajectories, we get that $\widehat T_\eps\ge T_X$, which ends the proof.
\end{proof}

We are now in position to prove that $b_{i}^\eps$ is close to $X_{i}^\eps$.

\begin{proposition}\label{prop:full_dynamic}
 Provided $T_0\leq \min(1/(24 C'_F),T_X)$, there exist $\eps_0 > 0$ and $C>0$ independent of $T_0$, such that for every $\eps \in (0,\eps_0]$, $i\in\{1,\ldots,N\}$ and $t \in [0,T_\eps]$, we have that
\begin{equation*}
 |b_i^\eps(t) - X_i^\eps(t)| \le C\frac{\ln |\ln \eps|}{|\ln \eps|}.
\end{equation*}
\end{proposition}

\begin{proof}
 Let $\mathcal{X}^\eps(t) = (X_1^\eps(t),\ldots,X_N^\eps(t))$, $\mathcal{B}^\eps(t) = (b_1^\eps(t),\ldots,b_N^\eps(t))$. 
 
 Let
 \begin{equation*}
 \mathcal{F}^\eps_{2} : \begin{cases}
 (\R^2_+)^N \to (\R^2)^N \\
 \displaystyle \mathcal{X} = (X_1,\ldots,X_N) \mapsto \left( \frac{\gamma}{4\pi X_{i,r}} e_z + \frac{\gamma}{2\pi |\ln \eps|}\sum_{j \neq i} \frac{\big(X_i-X_j\big)^\perp}{\big|X_i-X_j\big|^2}\right)_{1 \le i \le N}.
 \end{cases}
 \end{equation*}
 Then, we have by definition of the $X_i^\eps$ that
 \begin{equation*}
 \der{}{t} \mathcal{X}^\eps(t) = \mathcal{F}^\eps_{2}(\mathcal{X}^\eps(t))
 \end{equation*}
 and by Lemma~\ref{lem:eq_b} and Lemma~\ref{lem:pour_gronwall} that
 \begin{equation*}
 \left|\der{}{t} \mathcal{B}^\eps(t) - \mathcal{F}^\eps_{2}(\mathcal{B}^\eps(t))\right| \le C \frac{\ln |\ln \eps|}{|\ln \eps|}.
 \end{equation*}
 By Lemma~\ref{lem:X-tildeX}, the map $\mathcal{F}^\eps_{2}$ is Lipschitz continuous on $\prod_{i} \mathcal{A}_{\eps}^{X^\eps_{i,0,r}}\left(d_0,h_0\right)$ with constant $K = \mathcal{O}(1)$, so by Gronwall lemma (see for instance the variant in Lemma~\ref{lem:gronwall}), we get that for every $t \le T_\eps$,
 \begin{equation*}
 \left|\mathcal{B}^\eps(t) - \mathcal{X}^\eps(t)\right| \le \left(C \, T_0 \frac{\ln |\ln \eps|}{|\ln \eps|} + \eps \right) e^{KT_0},
 \end{equation*}
which ends the proof.
\end{proof}

From this result, we are now in position to prove Theorem~\ref{theo:main}-$(i)$ on the time interval $t \in [0,T_\eps]$.
\begin{corollary}\label{cor:weak-X}
 Provided $T_0\leq \min(1/(24 C'_F),T_X)$, there exist $\eps_0 > 0$ and $C_{0}>0$ independent of $T_0$, such that for every $\eps \in (0,\eps_0]$, $i\in\{1,\ldots,N\}$, we have that
 \begin{equation*}
 \sup_{t \in [0,T_\eps]}\left| \gamma - \int_{B\left( X_i^\eps(t) , C_0\frac{\ln |\ln \eps|}{|\ln \eps|}\right)} x_r\omega_i^\eps(t,x)\dd x\right| \le \gamma\frac{\ln |\ln \eps|}{|\ln \eps|}.
 \end{equation*}
\end{corollary}

\begin{proof}
 By Proposition~\ref{prop:full_dynamic} and as $h_0\leq\frac12$ (see Remark~\ref{rem:h0}), we have that 
 \begin{equation*}
 B\left(b_i^\eps(t), 3h_{0}\frac{\ln |\ln \eps|}{|\ln \eps|} \right) \subset B\left(X_i^\eps(t), (C+3h_{0}) \frac{\ln |\ln \eps|}{|\ln \eps|} \right)\subset B\left(X_i^\eps(t), \Big(C+\frac32\Big) \frac{\ln |\ln \eps|}{|\ln \eps|} \right),
 \end{equation*}
 hence
 \begin{equation*}
B\left(X_i^\eps(t), \Big(C+\frac32\Big) \frac{\ln |\ln \eps|}{|\ln \eps|} \right)^c \subset B\left(b_i^\eps(t), 3h_{0}\frac{\ln |\ln \eps|}{|\ln \eps|} \right)^c,
 \end{equation*}
 and
 \[
 \int_{B\left( X_i^\eps(t) , (C+\frac32) \frac{\ln |\ln \eps|}{|\ln \eps|}\right)^c} x_r\omega_i^\eps(t,x)\dd x \leq \int_{B\left( b_i^\eps(t) , 3h_{0} \frac{\ln |\ln \eps|}{|\ln \eps|}\right)^c} x_r\omega_i^\eps(t,x)\dd x.
 \]
 So using Corollary~\ref{cor:weak_loc}, we have the desired result by taking $C_0 = C+3/2$.
\end{proof}

\begin{remark}
 It is clear from the previous proof that we have the same result replacing $X_i^\eps$ by $\tilde X_i^\eps$, which proves the claim made in the introduction.
\end{remark}

\section{Strong localization}\label{sec:strong_loc}

All the above results are valid up to $T_{\eps}$, defined in \eqref{def:Teps}, in particular Part~\textit{(i)} of Theorem~\ref{theo:main}. The goal of this section is to show Part \textit{(ii)} of Theorem~\ref{theo:main}, but also to state that $T_{\eps}=T_{0}$, namely we will prove the following.

\begin{proposition}\label{prop:locforte}
For $h_0$ fixed, there exists $T_0 > 0$ and $\eps_0>0$ such that for every $\eps \in (0,\eps_0]$, every $t \le T_\eps$ and $i\in \{1,\ldots,N\}$, we have that
\begin{equation*}
 \supp \omega_i^\eps(t,\cdot) \subset \mathcal{A}_{\eps}^{X^\eps_{i,0,r}}\left(\frac{d_0}{2},\frac{h_0}{2}\right)
\end{equation*}
and thus, by the continuity of the trajectories of the fluid particles, we conclude that
\begin{equation*}
 T_\eps = T_0.
\end{equation*} 
\end{proposition}

We start by establishing the strong localization in the radial direction, which requires precise estimates and a two-step iteration process inspired by \cite{Mar24}, then turn to the vertical strong localization proven by a simpler version of the same method.

This strong localization method was introduced by Marchioro and Pulvirenti in the nineties, first published in \cite{MarPul93}, see also \cite{MP94}, in the context of $2D$-Euler flows. It has since proven to be a powerful tool for the study of concentrated vorticities and their stability and been deployed in many important articles. The approach has mainly been developed by Marchioro and his collaborators including adaptations to the $3D$ setting of vortex rings in various regimes. The arguments have been gradually refined and substantially improved in a long series of papers and noways reached an apparent level of complexity that may make it seem difficult to understand. 
Our aim in this section is to present this powerful approach in the most pedagogical and self-contained way possible, allowing the reader who has skimmed the previous sections to learn this tool.

\subsection{Strong localization in the radial direction}

To shorten the notation, let us denote for this section $r_0^{i,\eps} = X_{i,0,r}^\eps$.
Strong localization means controlling the size of the support of the vorticity, i.e. of the following quantity
\begin{equation*}
 R_t^{i,\eps} := \max \Big\{ \Big|x_r-r_0^{i,\eps}\Big| \, , \, x \in \supp \omega^\eps_{i}(t,\cdot)\Big\}.
\end{equation*}
Although the natural definition would be to use the Euclidian distance in $\R^2$, the presence of the local induction term $v_L^\eps$ in 3D axisymmetric flows, as noticed in \cite{Mar3}, makes this problem anisotropic and forces us to treat independently the radial and the vertical directions. Because of possible filamentation, we obtain a stronger control in the radial direction rather than in the vertical direction.
Let us mention that we could also define this radius by measuring the distance to the center of mass $b^\eps_{i}(t)$ rather than to $X_{i,0,r}$. However, in the current setting, it would not improve the obtained estimates but would result in longer computations.

To estimate the evolution of $R_t^{i,\eps}$, we now fix a time $t \le T_\eps$ and an index $i\in \{1,\ldots,N\}$, and focus on a fluid particle trajectory from the support of $\omega_i^\eps$ that reaches at time $t$ the farthest distance from the initial position (in the radial direction). Namely, we consider $x^{i,\eps,0} \in \R^2_+$ such that 
\begin{equation*}
 \Big|x^{i,\eps,0}_r - r_0^{i,\eps}\Big| = R_t^{i,\eps},
\end{equation*}
and $s \mapsto x^{i,\eps,t}(s)$ be the Lagrangian trajectory passing by $x^{i,\eps,0}$ at time $t$, i.e., by the transport equation \eqref{eq:Axi2D_rescaled_F}, the solution of
\begin{equation*}
 \left\{\begin{aligned}
& \der{}{s} x^{i,\eps,t}(s) =\frac1{|\ln\eps|}\Big( u^\eps_{i}(s,x^{i,\eps,t}(s)) + F^\eps_{i}(s,x^{i,\eps,t}(s))\Big), \\
& x^{i,\eps,t}(t) = x^{i,\eps,0}.
\end{aligned}\right.
\end{equation*}

The general idea of the proof is to control the velocity of this fluid particle in the radial direction (and hence the growth of $R_t^{i,\eps}$) by controlling the mass of the vorticity included in a region located between $R_t^{i,\eps}/2$ and $R_t^{i,\eps}$. Indeed, the fluid particle we focus on is relatively far from the vortex core, which gives good control over the velocity induced by the latter on the former. Therefore, we prove that estimating the vorticity mass far from the vortex core (hence potentially close to the fluid particle) induces control over its outward velocity. This is the content of the following lemma. Let us introduce for any $R >0$,
\begin{equation*}
 \mathcal{U}^{i,\eps}_R := \{ (z,r)\in \R^2_+, \; |r-r_0^{i,\eps}| < R \},
\end{equation*}
and
\begin{equation*}
 m_t^{i,\eps}(R) := \int_{(\mathcal{U}^{i,\eps}_R)^c} x_r\omega^\eps_{i}(t,x)\dd x.
\end{equation*}

\begin{lemma}\label{lemma:Rt}
There exists $C>0$ independent of $T_{0}$ such that 
\begin{equation*}
 \der{}{s} \Big[|x^{i,\eps,t}_r(s) - r_0^{i,\eps}|\Big](t) \le \frac{C}{\sqrt{|\ln\eps|}}+\frac{C}{|\ln\eps|} \left( \frac{1}{R_t^{i,\eps}} + \frac{\sqrt{m_t^{i,\eps}(R_t^{i,\eps}/2)}}{\eps}\right),
\end{equation*}
for all $t\in [0,T_{\eps}]$, $\eps>0$ and $i\in \{0,\ldots,N\}$.
\end{lemma}

\begin{proof}
As often in the previous sections, we drop during the proofs the index $i$ in the notation of $\omega_i^\eps$ and $\Omega_i^\eps$ since we perform the proofs for a fixed $i\in \{0,\ldots,N\}$ and we derive uniform estimates with respect to this parameter.

 We compute that
 \begin{equation*}
 \der{}{s} \Big[|x^{\eps,t}_r(s) - r^{\eps}_0|\Big](t) = \frac{ x^{\eps,t}_r(t)- r^{\eps}_0}{\big|x^{\eps,t}_r(t) - r^{\eps}_0\big|}\frac{e_r}{|\ln\eps|} \cdot (u^\eps+F^\eps)(t,x^{\eps,t}(t)) 
 \end{equation*}
and thus, recalling that $x^{\eps,t}(t) = x^{\eps,0}$,
\begin{align*}
 \der{}{s} \Big[|x^{\eps,t}_r(s) - r^{\eps}_0|\Big](t) & \le \frac{1}{|\ln\eps|} \Big| u_K^\eps(t,x^{\eps,0}) + u_R^\eps(t,x^{\eps,0})+F^\eps(t,x^{\eps,0})\Big|\\
 &\le \frac{C}{|\ln \eps|}\int_{\R^2_+} \frac{y_r\omega^\eps(t,y)}{\big|x^{\eps,0}-y\big|} \dd y + \frac{C}{\sqrt{|\ln\eps|}},
\end{align*}
where we have used, as often in Section~\ref{sec:weak_loc}, the decomposition of $u^\eps$ (see Lemma~\ref{lem:decomp_u}), that $u_L^\eps\cdot e_r =0$, the estimate on $F^\eps$ (see Lemma~\ref{lem:hyp_F}) and that $r^*/2\leq x^{\eps,0},y\leq 2r^*$ (see \eqref{est:supp}).
We now split the integral on $ \mathcal{U}^{\eps}_{R_t^{\eps}/2}$ and $\cal \R^2_+ \setminus \mathcal{U}^{\eps}_{R_t^{\eps}/2}.$ If $y \in \mathcal{U}^{\eps}_{R_t^{\eps}/2}$, then $|x^{\eps,0}-y| \ge R_t^{\eps}/2$, and thus, by the conservation of the total mass of vorticity \eqref{eq:consgamma} and as the vorticity is non-negative, we have
\begin{equation*}
 \int_{ \mathcal{U}^{\eps}_{R_t^{\eps}/2}} \frac{y_{r}\omega^\eps(t,y)}{|x^{\eps,0}-y|}\dd y \le \frac{C}{R_t^{\eps}}.
\end{equation*}
On the complement, we use the following property\footnote{See for instance \cite[Lemma 2.1]{ISG99} when $\Omega=\R^2$ and $f\in L^1\cap L^\infty(\R^2)$, and next apply this for $f\mathds{1}_{\Omega}$.}: for every $\Omega \subset \R^2$ and $f\in L^1\cap L^\infty(\Omega)$
\begin{equation*}
 \int_\Omega \frac{f(y)}{|x^{\eps,0}-y|}\dd y \le C \| f \|_{L^1(\Omega)}^{\frac{1}{2}} \| f \|_{L^\infty(\Omega)}^{\frac{1}{2}},
\end{equation*}
where $C$ is independent of $\Omega$, to obtain that
\begin{equation*}
 \int_{\R^2_{+}\setminus \mathcal{U}^{\eps}_{R_t^{\eps}/2}} \frac{y_{r}\omega^\eps(t,y)}{|x^{\eps,0}-y|}\dd y \le \frac{C}{\eps}\sqrt{m_t^{\eps}(R_t^{\eps}/2)},
\end{equation*}
because the transported initial vorticity is less than $C/\eps^2$, see \eqref{hyp:omega0}.
\end{proof}

It appears now clear that we need to estimate the mass of the vorticity further than $R^{i,\eps}_t /2$, namely we are looking for a weak localization estimate. Unfortunately, Corollary~\ref{cor:weak_loc} derived in Section~\ref{sec:weak_loc} is not enough because $\ln |\ln \eps|/(\eps^2 |\ln \eps|) \gg |\ln \eps|$ so we will prove a better estimate for $m_t^{i,\eps}(R_t^{i,\eps}/2)$.

As we will need to differentiate with respect to $t$, which is difficult because the domain of integration in the definition of $m_t^{i,\eps}$ depends on $t$, we introduce a regularization of $m_t^{i,\eps}$: let
\begin{equation*}
 \mu_t^{i,\eps}(R,\eta) := \int_{\R^2} \big(1-W_{R,\eta}(x_r-r_0^{i,\eps})\Big) x_r\omega^\eps_{i}(t,x) \dd x
\end{equation*}
with
\begin{equation*}
 W_{R,\eta}(s) = 
 \begin{cases}
 1 & \text{ if } |s| \le R \\
 0 & \text{ if } |s| \ge R + \eta
 \end{cases} 
\end{equation*}
satisfying that $W_{R,\eta} \in C^\infty$ and $\|W_{R,\eta}'\|_{L^\infty} \le C/\eta$ and $\| W_{R,\eta}''\|_{L^\infty} \le C/\eta^2$.

We first observe that
\begin{equation}\label{eq:encadrement_mu}
 \mu_t^{i,\eps}(R,\eta) \le m_t^{i,\eps}(R) \le \mu_t^{i,\eps}(R-\eta,\eta).
\end{equation}

We estimate in the next lemma the time derivative.

\begin{lemma}\label{lem:a_iterer}
There exists $C>0$ independent of $T_{0}$ such that, for every $t \le T_\eps$, $\eps>0$, $i\in \{1,\ldots, N\}$, $\eta>0$ and $R \in \left[\frac{d_0}{128\sqrt{|\ln \eps|}},\frac{d_0}{16\sqrt{|\ln \eps|}}\right]$, we have that
 \begin{equation*}
 \der{}{t} \mu_t^{i,\eps}(R,\eta) \le C \left( \frac{1}{\eta \sqrt{|\ln\eps|}} + \frac{m_t^{i,\eps}(R/2)}{\eta^2|\ln\eps|}\right)m_t^{i,\eps}(R)
 \end{equation*}
\end{lemma}

\begin{proof}
As in the beginning of the proof of Proposition~\ref{prop:der_energy}, we use the equation verified by $\omega^\eps$ \eqref{eq:Axi2D_rescaled_F} and the decomposition of $u^\eps$ \eqref{def:u_K,u_L,u_R}, to state that $u_L^\eps \cdot e_r = 0$ implies
\begin{align*}
 \der{}{t} \mu_t^{\eps}(R,\eta) &= \frac{-1}{|\ln\eps|} \int \big(1-W_{R,\eta}(x_r-r_0^{i,\eps})\Big) \div \Big( (u^\eps + F^\eps)(t,x) x_r \omega^\eps(t,x)\Big) \dd x \\
 &= \frac{-1}{|\ln\eps|} \int W_{R,\eta}'(x_r-r^{\eps}_0) e_r \cdot (u_K^\eps + u_R^\eps + F^\eps)(t,x) x_r \omega^\eps(t,x) \dd x = A_1 + A_2
\end{align*}
where
\begin{align*}
 A_1 & := \frac{-1}{|\ln\eps|} \int W_{R,\eta}'(x_r-r^{\eps}_0) e_r \cdot (u_R^\eps + F^\eps) x_r \omega^\eps(t,x) \dd x \\
 A_2 & := \frac{-1}{|\ln\eps|} \int W_{R,\eta}'(x_r-r^{\eps}_0) e_r \cdot u_K^\eps x_r \omega^\eps(t,x) \dd x.
\end{align*}

We start with $A_1$. Recalling from Lemmas~\ref{lem:decomp_u} and \ref{lem:hyp_F} that $\|u_R^\eps\|_{L^\infty} \le C$ and $\|F^\eps\|_{L^\infty} \le C\sqrt{|\ln \eps|}$, and that $W_{R,\eta}'= 0$ out of $[R,R+\eta]$ we have that
\begin{equation*}
 |A_1| \le \frac{C}{\eta \sqrt{|\ln\eps|}}m_t^{\eps}(R).
\end{equation*}

We now compute $A_2$. First, we expand $u_K^\eps$ and get that
\begin{equation*}
 A_2 = \frac{-1}{|\ln\eps|} \iint W_{R,\eta}'(x_r-r^{\eps}_0) e_r \cdot \sqrt{\frac{y_r}{x_r}} K(x-y) x_r \omega^\eps(t,x) y_r\omega^\eps(t,y) \dd x \dd y
\end{equation*}
then we symmetrize by exchanging $x$ and $y$ to get that
\begin{equation*}
 A_2 = \frac{-1}{|\ln\eps|} e_r \cdot \iint f(x,y) K(x-y)x_r \omega^\eps(t,x) y_r\omega^\eps(t,y) \dd x \dd y
\end{equation*}
with
\begin{equation*}
 f(x,y) = \frac12\Big(W_{R,\eta}'(x_r-r^{\eps}_0)\sqrt{\frac{y_r}{x_r}} - W_{R,\eta}'(y_r-r^{\eps}_0)\sqrt{\frac{x_r}{y_r}}\Big).
\end{equation*}

For readability, we omit the subscript $\eps$ in the set of type $\mathcal{U}^{\eps}_{R}$. As $f=0$ on $\mathcal{U}_{R}\times \mathcal{U}_{R}$, we observe that for any $\varphi\in L^\infty$
\begin{align*}
 \iint_{\R^2_+\times \R^2_+} f\varphi &= \iint_{\R^2_+\times \mathcal{U}_{R/2}} f\varphi + \iint_{\R^2_+\times \mathcal{U}^c_{R/2}} f\varphi\\
 &= \iint_{\mathcal{U}^c_{R} \times \mathcal{U}_{R/2}} f\varphi + \iint_{\mathcal{U}^c_{R} \times \mathcal{U}^c_{R/2}} f\varphi + \iint_{\mathcal{U}_{R} \times \mathcal{U}_{R/2}^c} f\varphi 
\end{align*}
and
\begin{align*}
 \iint_{\mathcal{U}_{R} \times \mathcal{U}_{R/2}^c} f\varphi + \iint_{\mathcal{U}_{R}^c \times \mathcal{U}_{R}^c} f\varphi
 &= \iint_{\mathcal{U}_{R} \times \mathcal{U}_{R}^c} f\varphi + \iint_{\mathcal{U}_{R}^c \times \mathcal{U}_{R}^c} f\varphi=\iint_{\R^2_{+} \times \mathcal{U}_{R}^c} f\varphi\\
 &=\iint_{\mathcal{U}_{R/2} \times \mathcal{U}_{R}^c} f\varphi +\iint_{\mathcal{U}_{R/2}^c \times \mathcal{U}_{R}^c} f\varphi .
\end{align*}
Summing these two lines, we get that for symmetric $f\varphi$, we may split the integral in three as
\begin{equation*}
 \iint_{\R^2_+\times \R^2_+} f\varphi= 2\iint_{ \mathcal{U}^c_R \times \mathcal{U}_{R/2}}f\varphi + 2\iint_{ \mathcal{U}^c_R \times \mathcal{U}^c_{R/2}} f\varphi- \iint_{ \mathcal{U}^c_R \times \mathcal{U}^c_R}f\varphi
\end{equation*}
and we denote accordingly $A_2 := A_{21} + A_{22} + A_{23}$.

When $(x,y) \in \mathcal{U}^c_R \times \mathcal{U}_{R/2}$, $|x-y| \ge R/2\ge \frac{d_0}{256\sqrt{|\ln \eps|}}$, hence
\begin{equation*}
 |A_{21}| \le \frac{C}{\eta \sqrt{|\ln \eps|}} m_t^{\eps}(R).
\end{equation*}
For the two other terms, we use that $|f(x,y)| \le \frac{C}{\eta^2}|x-y|$ to get that
\begin{equation*}
 |A_{22}+A_{23}| \le \frac{C}{\eta^2 |\ln\eps|} \iint_{\mathcal{U}_R^c\times \mathcal{U}_{R/2}^c} x_r\omega^\eps(t,x) y_r \omega^\eps(t,y) \dd x\dd y \le \frac{C}{\eta^2|\ln\eps|} m_t^{\eps}(R) m_t^{\eps}(R/2),
\end{equation*}
which ends the proof.
\end{proof}

We now use the previous lemma to estimate the quantity $m_t^{i,\eps}(R_t^{i,\eps}/2)$. It is not possible to reach an estimate of order $\eps^\ell$ in one step, so we use a double iterative method inspired by \cite{Mar24}. 

Starting from the weak concentration proved in Section~\ref{sec:weak_loc}, we get a weak concentration estimate of order $|\ln\eps|^{-\ell}$.

\begin{lemma}[First iteration]\label{lem:step2}
For every $\ell > 0$, there exists a $T_0$ depending on $\ell$, that we denote by $T_{0,\ell}>0$, and $\eps_\ell>0$, such that for every $\eps \in (0,\eps_\ell]$, $i\in\{1,\ldots,N\}$ and $t\in [0,T_\eps]$, we have that
\begin{equation*}
 m_t^{i,\eps}\left(\frac{d_0}{64\sqrt{|\ln \eps|}}\right) \le \frac{1}{|\ln\eps|^\ell}.
\end{equation*}
\end{lemma}
In the statement of this lemma, we read the definition of $T_\eps$ \eqref{def:Teps} as $T_\eps\leq T_{0,\ell}$.

\begin{proof}
 As the dynamics of $(\tilde Y_i)$ is independent of $\eps$ (see \eqref{eq:LP}), it is clear that we can state that
 \[
|\tilde{X}_r^\eps(t)- r^{\eps}_0| \le \frac{d_0}{300\sqrt{|\ln \eps|}}, \quad \forall t\in [0,T_\eps]
 \]
 by choosing $T_0$ small enough\footnote{This assumption may seem not very natural. A possibility to remove it is to confine the solution in a strip centered not on $r_0^\eps$, which is fixed, but on $X_{i,r}^\eps(t)$, or even on $b_{i,r}^\eps(t)$. However in the present paper, as other limitations exist on $T_0$, it does not appear have any real interest, while adding more technical difficulties.}. From Lemma~\ref{lem:X-tildeX} and up to choosing $\eps$ small enough, we have
 \begin{equation*}
 |X_r^\eps(t)- r^{\eps}_0| \le \frac{d_0}{256\sqrt{|\ln \eps|}}\quad \forall t\in [0,T_\eps].
 \end{equation*}
 For such $T_0$ and $\eps_0$, we notice that
 \[
 \cal U_{R/2}^{\eps c} \subset B\Big(X^\eps(t),C_0\frac{\ln|\ln\eps|}{|\ln\eps|}\Big)^c
 \]
 provided $R\in \left[\frac{d_0}{128\sqrt{|\ln \eps|}},\frac{d_0}{64\sqrt{|\ln \eps|}}\right]$ and $\eps_0$ are small enough. Corollary~\ref{cor:weak-X} then implies that 
 \begin{equation*}
 m_t^{\eps}(R/2) \le \gamma\frac{\ln |\ln \eps|}{|\ln \eps|}.
 \end{equation*} 
 Plugging this into the estimates of Lemma~\ref{lem:a_iterer}, we get that
 \begin{equation*}
 \der{}{t} \mu_t^{\eps}(R,\eta) \le C \left(\frac{1}{\eta \sqrt{|\ln\eps|}} + \frac{\ln |\ln \eps|}{\eta^2|\ln\eps|^2}\right)m_t^{\eps}(R),
 \end{equation*}
 for any $R\in \left[\frac{d_0}{128\sqrt{|\ln \eps|}},\frac{d_0}{64\sqrt{|\ln \eps|}}\right]$.
 Now we set $R_0 = \frac{d_0}{64\sqrt{|\ln \eps|}}$, $n = [\ln |\ln \eps|]$ and $\eta = \frac{R_0}{2n}$ and we obtain that
 \begin{align*}
 \frac{1}{\eta \sqrt{|\ln\eps|}} + \frac{\ln |\ln \eps|}{\eta^2|\ln\eps|^2} = \frac{2n}{R_0\sqrt{|\ln\eps|}} + \frac{4n^2 \ln |\ln \eps|}{R^2_0|\ln\eps|^2} \le C n + Cn \frac{(\ln |\ln \eps|)^2}{|\ln\eps|} \le Cn.
 \end{align*}
 In conclusion,
 \begin{equation*}
 \der{}{t} \mu_t^{\eps}(R,\eta) \le C\, n \, m_t^{\eps}(R).
 \end{equation*}
 Integrating, we get that
 \begin{equation*}
 \mu_t^{\eps}(R,\eta) \le \mu_0^{\eps}(R,\eta) + C \,n\int_0^t m_s^{\eps}(R)\dd s
 \end{equation*}
 and using \eqref{eq:encadrement_mu} we get that
 \begin{equation*}
 \mu_t^{\eps}(R,\eta) \le \mu_0^{\eps}(R,\eta) + Cn\int_0^t \mu_s^{\eps}(R-\eta,\eta)\dd s,
 \end{equation*}
 provided $R\in \left[\frac{R_0}2,R_0\right]$.

 Applying this inequality twice we get that
 \begin{align*}
 \mu_t^{\eps}(R_0-\eta,\eta) 
 &\le \mu_0^{\eps}(R_0-\eta,\eta) + Cn\int_0^t \mu_{s_1}^{\eps}(R_0-2\eta,\eta)\dd s_1\\
 &\le \mu_0^{\eps}(R_0-\eta,\eta) + Cnt\mu_0^{\eps}(R_0-2\eta,\eta) + (Cn)^2\int_0^t\int_0^{s_1} \mu_{s_2}^{\eps}(R_0-3\eta,\eta)\dd s_2 \dd s_1
 \end{align*}
 and iterating $n-1$ times eventually leads to
 \begin{equation*}
 \mu_t^{\eps}(R_0-\eta,\eta) \le \sum_{j=1}^{n-1} \mu_0^{\eps}(R_0-j\eta,\eta) (Cnt)^{j-1} + (Cn)^{n-1} \int_0^t\ldots\int_0^{s_{n-2}} \mu_{s_{n-1}}^{\eps}(R_0/2,\eta) \dd s_{n-1} \ldots \dd s_1
 \end{equation*}
 and therefore, since $\mu_0^{\eps}(R_0/2) = 0$ (see the initial assumption \eqref{hyp:omega0}) and $\mu_s^\eps \le \gamma$ for every $s$, we obtain that
 \begin{equation*}
 \mu_t^{\eps}(R_0-\eta,\eta) \le \frac{(Cnt)^{n-1}}{(n-1)!}.
 \end{equation*}
 Using Stirling's formula we have that for any $\ell >0$,
 \begin{align*}
 |\ln \eps|^\ell \mu_t^{\eps}(R_0-\eta,\eta) & \le \frac{1}{\sqrt{\pi (n-1)}} \exp\Big( (n-1) \big(\ln (Cnt) - \ln(n-1) + 1\big) + \ell \ln |\ln \eps|\Big) 
 \\ & \le \frac{2}{\sqrt{\pi \ln |\ln \eps|}} \exp\Big( \ln |\ln \eps| \big(\ln (Ct) + 2 + \ell \big)\Big).
 \end{align*}
 We see that if $T_{0,\ell} \le \frac{1}{Ce^{\ell +2}}$, then this expression goes to 0 as $\eps \to 0$ uniformly in $t \in [0,T_\eps]$. By \eqref{eq:encadrement_mu}, this ends the proof.
\end{proof}

We iterate this procedure, starting now from the weak concentration proved in Lemma~\ref{lem:step2}, to get a weak concentration estimate of order $\eps^{\ell}$.

\begin{lemma}[Second iteration]\label{lem:step3}
For every $\ell > 0$, there exists $T_{0,\ell}>0$ and $\eps_\ell>0$ depending on $\ell$, such that for every $\eps \in (0,\eps_\ell]$, $i\in\{1,\ldots,N\}$ and $t\in [0,T_\eps]$, we have that
\begin{equation*}
 m_t^{i,\eps}\left(\frac{d_0}{16\sqrt{|\ln \eps|}}\right) \le \eps^\ell.
\end{equation*}
\end{lemma}

\begin{proof}
 Applying Lemma~\ref{lem:step2} for $\ell =1$, we state that for every $R \in \left[ \frac{d_0}{32\sqrt{|\ln \eps|}}, \frac{d_0}{16\sqrt{|\ln \eps|}}\right]$,
 \[
 m_t^\eps(R/2)\leq m_t^\eps\Big( \frac{d_0}{64\sqrt{|\ln \eps|}}\Big)\leq \frac1{|\ln\eps|},
 \]
 so by injecting this into Lemma~\ref{lem:a_iterer}, we get
 \begin{equation*}
 \der{}{t} \mu_t^{\eps}(R,\eta) \le C \left(\frac{1}{\eta \sqrt{|\ln\eps|}} + \frac{1}{\eta^2|\ln\eps|^2}\right)m_t^{\eps}(R).
 \end{equation*}
 By taking this time $R_0=\frac{d_0}{16\sqrt{|\ln\eps|}}$, $\eta = \frac{R_0}{2n}$ and $n = [|\ln\eps|]$, we get again that
 \begin{equation*}
 \der{}{t} \mu_t^{\eps}(R,\eta) \le Cn m_t^{\eps}(R).
 \end{equation*}
 Then by the same iterative method than in Lemma~\ref{lem:step2}, we obtain that
 \begin{equation*}
 \mu_t^{\eps}(R_0-\eta,\eta) \le \frac{(Cnt)^{n-1}}{(n-1)!}.
 \end{equation*}
 We then compute that
 \begin{align*}
 \eps^{-\ell} \mu_t^{\eps}(R_0-\eta,\eta) & \le \frac{1}{\sqrt{\pi (n-1)}} \exp\Big( (n-1) \big(\ln (Cnt) - \ln(n-1) + 1\big) + \ell |\ln \eps|\Big) 
 \\ & \le \frac{2}{\sqrt{\pi |\ln \eps|}} \exp\Big( |\ln \eps| \big(\ln (Ct) + 2 + \ell \big)\Big).
 \end{align*}
 which again goes to 0 if $T_{0,\ell} \leq \frac{1}{Ce^{\ell+2}}$.
 \end{proof}

 We are now able to use Lemma~\ref{lemma:Rt} to conclude the radial part of Proposition~\ref{prop:locforte}.
 
 \begin{lemma}\label{lem:radialstrong}
 There exist $T_{0}>0$ and $\eps_0>0$, such that for every $\eps \in (0,\eps_0]$, $i\in\{1,\ldots,N\}$ and $t\in [0,T_\eps]$, we have that
 $$R_t^{i,\eps} < \frac{d_0}{2\sqrt{|\ln \eps|}}.$$
 \end{lemma}
 
 \begin{proof}
 Assume that there exist $i\in \{1,\ldots,N\}$ and a time $t_1 \in (0,T_\eps]$ such that $ R_{t_1}^{\eps} = \frac{d_0}{2\sqrt{|\ln \eps|}}$. Then by continuity of $t\mapsto R_t^{\eps}$, there exists a time $t_0 < t_1$ such that
 \begin{equation*}
 R_{t_0}^\eps = \frac{d_0}{8\sqrt{|\ln \eps|}} \qquad \text{ and } \qquad \forall t \in [t_0,t_1], \, R_t^{\eps} \ge \frac{d_0}{16\sqrt{|\ln \eps|}}. 
 \end{equation*}
 Let $C_1$ be the constant in Lemma~\ref{lemma:Rt}. Let $f : [t_0,t_1] \to \R$ be the solution of
 \begin{equation*}
 \begin{cases}
 f(t_0) = \frac{d_0}{4\sqrt{|\ln \eps|}} \vspace{1mm}\\
 \displaystyle f'(t) = 2\frac{C_1} {\sqrt{|\ln\eps|}}+2\frac{C_1}{|\ln\eps|} \left( \frac{1}{f(t)} + \frac{\sqrt{m_t^{\eps}(R_t^{\eps}/2)}}{\eps}\right). 
 \end{cases}
 \end{equation*}
 Then for every $t \in [t_0,t_1]$,
 \begin{equation*}
 R_t^{\eps} < f(t).
 \end{equation*}
 Indeed, if it isn't the case, by letting $t_2 \in [t_0,t_1]$ be the first time when $R_{t_2}^\eps = f(t_2)$, then there exists a point $x^{\eps,0} \in \supp\omega^\eps(t_2)$ such that $|x^{\eps,0} - r^{\eps}_0| = R_{t_2}^\eps$. However Lemma~\ref{lemma:Rt} is then in contradiction\footnote{For more details on this contradiction argument, see for instance \cite[Section~4.4]{DI21}.} with the fact that $t_2$ is the first such time, due to the fact that $f(t_0)>R_{t_0}^\eps$ and
 \begin{equation*}
 f'(t_2) > \frac{C_1}{\sqrt{|\ln\eps|}}+\frac{C_1}{|\ln\eps|} \left( \frac{1}{R_{t_2}^\eps} + \frac{\sqrt{m_t^{\eps}(R_{t_2}^\eps/2)}}{\eps}\right). 
 \end{equation*}

 Since for every $t \in [t_0,t_1]$, $f(t) > R_t^{\eps} \ge \frac{d_0}{16\sqrt{|\ln \eps|}}$, we can apply Lemma~\ref{lem:step3} with $\ell = 2$ to get that for $\eps>0$ small enough,
 \begin{equation*}
 f'(t) \le \frac{C}{\sqrt{|\ln\eps|}}.
 \end{equation*}
 Therefore, by integrating in time,
 \begin{align*}
 R_{t_1}^\eps < f(t_1) \le f(t_0) + (t_1-t_0) \frac{C}{\sqrt{|\ln\eps|}} \le \frac{d_0}{2\sqrt{|\ln \eps|}} 
 \end{align*}
 provided $T_0$ is small enough. This is in contradiction with the fact that $R_{t_1}^\eps = \frac{d_0}{2\sqrt{|\ln \eps|}}$. By continuity of the trajectories, the conclusion of the lemma follows from the non-existence of such time $t_1$.
 \end{proof}

\subsection{Strong localization in the vertical direction}

We now perform a very similar analysis, but in the direction $z$. We introduce this time
\begin{equation*}
 \mathcal{U}_Z^{i,\eps} = \{ (z,r)\in \R^2_+, \; |z-z^*| < Z \}
\end{equation*}
and
\begin{equation*}
 m_t^{i,\eps}(Z) = \int_{(\mathcal{U}_Z^{i,\eps})^c} x_r\omega_i^\eps(t,x)\dd x
\end{equation*}
and let
\begin{equation*}
 \mu_t^{i,\eps}(Z,\eta) = \int_{\R^2} \big(1-W_{Z,\eta}(x_z-z^*)\Big) x_r\omega_i^\eps(t,x) \dd x.
\end{equation*}
Let $t \le T_\eps$, $i\in \{1,\ldots,N\}$ and let $x^{i,\eps,0} \in \R^2_+$ such that 
\begin{equation*}
 \Big|x^{i,\eps,0}_z - z^*\Big| = Z_t^{i,\eps} := \max \Big\{ \Big|x_z-z^*\Big| \, , \, x \in \supp \omega_i^\eps(t)\Big\}.
\end{equation*}
We then have the following.
\begin{lemma}
There exists $C>0$ independent of $T_0$ such that
\begin{equation*}
 \der{}{s} \Big[|x^{\eps,t}_z(s) - z^*|\Big](t) \le C+\frac{C}{|\ln\eps|} \left( \frac{1}{Z_t^{i,\eps}} + \frac{\sqrt{m_t^{i,\eps}(Z_t^{i,\eps}/2)}}{\eps}\right),
\end{equation*}
for all $t\in [0,T_\eps]$, $\eps>0$ and $i\in \{0,\ldots,N\}$.
\end{lemma}

\begin{proof}
 We proceed similarly than in Lemma~\ref{lemma:Rt}, the main difference being that the term $u_L^\eps \cdot e_z$ does not vanish. We compute that
 \begin{equation*}
 \der{}{s} \Big[|x^{\eps,t}_z(s) - z^*|\Big](t) = \frac{ x^{\eps,t}_z(t)-z^*}{\big|x^{\eps,t}_z(t) - z^*\big|}\frac{e_z}{|\ln\eps|} \cdot (u^\eps+F^\eps)(t,x^{\eps,t}(t)) 
 \end{equation*}
 and thus
 \[
 \der{}{s} \Big[|x^{\eps,t}_z(s) - z^*|\Big](t) \le \frac{1}{|\ln\eps|} \Big| u_K^\eps(t,x^{\eps,0}) + u_L^\eps(t,x^{\eps,0}) +u_R^\eps(t,x^{\eps,0})+F^\eps(t,x^{\eps,0})\Big|.
 \]
 We already know from the proof of Lemma~\ref{lemma:Rt} that
 \begin{equation*}
 \frac{1}{|\ln\eps|} \Big| u_K^\eps(t,x^{\eps,0}) + u_R^\eps(t,x^{\eps,0})+F^\eps(t,x^{\eps,0})\Big| \le \frac{C}{\sqrt{|\ln\eps|}}+\frac{C}{|\ln\eps|} \left( \frac{1}{Z_t^\eps} + \frac{\sqrt{m_t^\eps(Z_t^\eps/2)}}{\eps}\right).
 \end{equation*}
 We now use the relation \eqref{def:u_K,u_L,u_R} between $u_L^\eps$ and $\psi^\eps$, the estimate of $\psi^\eps$ (see Lemma~\ref{lem:encadrement_psi}) and the fact that $\tilde{T}_\eps=T_\eps$ (see Section~\ref{sec:bootstrapweak}) to state that
 \begin{equation*}
 |u_L^\eps(t,x^{\eps,0})| \le C |\ln \eps|
 \end{equation*}
 which allows us to conclude.
\end{proof}

\begin{lemma}\label{lem:a_iterer_Z}
There exists $C>0$ independent of $T_{0}$ such that, for every $t \le T_\eps$, $\eps>0$, $i\in \{1,\ldots, N\}$, $\eta>0$ and $Z \in \left[\frac{h_0}{32},\frac{h_0}{16}\right]$, we have that
 \begin{equation*}
 \der{}{t} \mu_t^{i,\eps}(Z,\eta) \le C \left( \frac{1}{\eta} + \frac{m_t^{i,\eps}(Z/2)}{\eta^2|\ln\eps|}\right)m_t^{i,\eps}(Z)
 \end{equation*}
\end{lemma}

\begin{proof}
We compute that
 \begin{equation*}
 \der{}{t} \mu_t(Z,\eta) = \frac{-1}{|\ln\eps|} \int W_{Z,\eta}'(x_z-z^*) e_z \cdot (u_K^\eps + u_L^\eps + u_R^\eps + F^\eps)(t,x) x_r \omega^\eps(t,x) \dd x = A_1 + A_2 + A_3
\end{equation*}
where
\begin{equation*}
 |A_1| + |A_2| \le C \left( \frac{1}{\eta\sqrt{|\ln\eps|}} + \frac{m_t^\eps(Z/2)}{\eta^2|\ln\eps|}\right)m_t^\eps(Z)
\end{equation*}
by reproducing the arguments of Lemma~\ref{lem:a_iterer} with the only difference that $Z/2 \ge \frac{h_0}{64}$, and
\begin{equation*}
 A_3 := \frac{-1}{|\ln\eps|} \int W_{Z,\eta}'(x_z-z^*) e_z \cdot u_L^\eps(t,x) x_r \omega^\eps(t,x) \dd x
\end{equation*}
satisfies
\begin{equation*}
 |A_3| \le C/\eta.
\end{equation*}
\end{proof}

We only need one iteration process to conclude in that case.

\begin{lemma}
For every $\ell > 0$, there exists $T_{0,\ell}>0$ and $\eps_\ell>0$ depending on $\ell$, such that for every $\eps \in (0,\eps_\ell]$, $i\in\{1,\ldots,N\}$ and $t\in [0,T_\eps]$, we have that
\begin{equation*}
 m_t^{i,\eps}\left(\frac{h_0}{16}\right) \le \eps^\ell.
\end{equation*}
\end{lemma}

\begin{proof}
 As in the beginning of the proof of Lemma~\ref{lem:step2}, we remark that, up to choose $T_0$ and $\eps_0$ small enough, Corollary~\ref{cor:weak-X} implies that
 \begin{equation*}
 m_t^\eps(Z/2) \le \gamma \frac{\ln |\ln \eps|}{|\ln \eps|},
 \end{equation*}
 provided $Z \in \left[ \frac{h_0}{32}, \frac{h_0}{16}\right]$.
 Injecting this into Lemma~\ref{lem:a_iterer_Z}, we get that
 \begin{equation*}
 \der{}{t} \mu_t^\eps(Z,\eta) \le C \left( \frac{1}{\eta} + \frac{\ln |\ln \eps|}{\eta^2|\ln\eps|^2}\right)m_t^\eps(Z).
 \end{equation*}
 We take now $Z_0=\frac{d_0}{16}$ $\eta = \frac{Z_0}{2n}$ with $n = [|\ln\eps|]$, and get that
 \begin{equation*}
 \der{}{t} \mu_t^\eps(Z,\eta) \le Cn m_t^\eps(Z),
 \end{equation*}
 for all $Z\in \left[ \frac{h_0}{32}, \frac{h_0}{16}\right]$.
 Then by the same iterative method than in Lemma~\ref{lem:step2}, we obtain that
 \begin{equation*}
 \mu_t^\eps(Z_0-\eta,\eta) \le \frac{(Cnt)^{n-1}}{(n-1)!}.
 \end{equation*}
 As in Lemma~\ref{lem:step3} we conclude that
 \begin{align*}
 \eps^{-\ell} \mu_t^\eps(Z_0-\eta,\eta) & \le \frac{1}{\sqrt{\pi (n-1)}} \exp\Big( (n-1) \big(\ln (Cnt) - \ln(n-1) + 1\big) + \ell |\ln \eps|\Big) 
 \\ & \le \frac{2}{\sqrt{\pi |\ln \eps|}} \exp\Big( |\ln \eps| \big(\ln (Ct) + 2 + \ell \big)\Big)
 \end{align*}
 which gives the result as soon as $C T_{0,\ell} \leq 1/e^{2+\ell}$.
 \end{proof} 
 
 We then conclude the vertical part of Proposition~\ref{prop:locforte}.
 
 \begin{lemma}\label{lem:vertstrong}
 There exist $T_{0}>0$ and $\eps_0>0$, such that for every $\eps \in (0,\eps_0]$, $i\in\{1,\ldots,N\}$ and $t\in [0,T_\eps]$, we have that
 \begin{equation*}
 Z_t^{i,\eps} < \frac{h_0}{2}.
 \end{equation*}
 \end{lemma}
 
 \begin{proof}
 Assume that this is false, we obtain this time that $Z_t^\eps < f(t)$ for $t \in [t_0,t_1]$ where 
 \begin{equation*}
 \begin{cases}
 f(t_0) = \displaystyle \frac{h_0}{4}, \vspace{1mm} \\
 f'(t) = \displaystyle 2C_2+\frac{2C_2}{|\ln\eps|} \left( \frac{1}{f(t)} + \frac{\sqrt{m_t^\eps(Z_t^\eps/2)}}{\eps} \right)
 \end{cases}
 \end{equation*}
 and $f(t) \ge \frac{h_0}{16}$, and thus
 \begin{equation*}
 f'(t) \le C,
 \end{equation*}
 so that
 \begin{equation*}
 f(t) \le \frac{h_0}{4} +C(t_1-t_0).
 \end{equation*}
 Taking $T_0$ small enough leads to the contradiction and concludes the proof.
 \end{proof}

\subsection{Conclusion of the proof of Theorem~\ref{theo:main}}

By proving Lemmas~\ref{lem:radialstrong} and \ref{lem:vertstrong}, we have completed the proof of Proposition~\ref{prop:locforte} which implies Part \textit{(ii)} of Theorem~\ref{theo:main}, see \eqref{def:Aeta} for the definition of $\cal{A}_\eps^{r_0}(d,h_0)$.

Proposition~\ref{prop:locforte} also implies that $T_\eps = T_0$, so all the estimates we obtained throughout this article hold on the time interval $[0,T_0]$ which is independent of $\eps$. In particular, we recall that Corollary~\ref{cor:weak-X} yields Part~\textit{(i)} of Theorem~\ref{theo:main}.

\bigskip
\noindent
{\bf Acknowledgements.} This work was supported by the BOURGEONS project, grant ANR-23-CE40-0014-01 of the French National Research Agency (ANR). L.E.H. is funded by the Deutsche Forschungsgemeinschaft (DFG, German Research Foundation) – Project-
ID 258734477 – SFB 1173. C.L. also benefited of the support of the ANR under France 2030 bearing the reference ANR-23-EXMA-004 (Complexflows project of the PEPR MathsViVEs).

\appendix
\section{Two classical lemmas}\label{app:lemma}

In this section, we state two classical lemmas which are used in the article.

First, we recall a standard lemma for mass rearrangement, whose proof is available in \cite[Lemma B.1]{HLM24}.
Let $g$ be a non increasing continuous function from $(0,+\infty)$, non-negative, such that $s \mapsto sg(s) \in L^1_{\mathrm{loc}}\big([0,\infty)\big)$. Let
\begin{equation*}
 \mathcal{E}_{M,\gamma} = \left\{ f \in L^\infty_c(\R^2) \, , \, 0 \le f \le M \, , \, \int f = \gamma \right\}.
\end{equation*}

\begin{lemma}\label{lem:rearrangement}
For all $x\in\R^2$, we have
\begin{equation*}
 \max_{f \in \mathcal{E}_{M,\gamma}} \int_{\R^2} g\big(|x-y|\big) f(y)\dd y = 2\pi M \int_0^R sg(s) \dd s
\end{equation*}
where $R = \sqrt{\frac{\gamma}{\pi M}}$, namely that $f^* = M\Ind_{B(x,R)}$ is the map that maximizes this quantity on $\mathcal{E}_{M,\gamma}$.
\end{lemma}

The second lemma concerns a variant of Gronwall's inequality, whose proof is provided in \cite[Lemma A.2]{DonatiLacaveMiot}.

\begin{lemma}\label{lem:gronwall}
Let $f : \R^n \to \R^n$ such that there exists $\kappa$ such that
\begin{equation*}
 \forall \, x,y \in \R^n, \quad \big|f(x)-f(y)\big| \le \kappa|x-y|.
\end{equation*}
Let $g \in L^1(\R_+,\R_+)$ and $T \ge 0$. We assume that $z : \R_+ \to \R^n$ satisfies
\begin{equation*}
 \forall t \in [0,T], \quad z'(t) = f(z(t)),
\end{equation*}
that $y : \R_+ \to \R^n$ satisfies
\begin{equation*}
 \forall t \in [0,T], \quad |y'(t)-f(y(t))| \le g(t).
\end{equation*}
Then
\begin{equation*}
 \forall t \in [0,T], \quad |y(t) - z(t)| \le \left( \int_0^t g(s) \dd s + |y(0)-z(0)|\right)e^{\kappa t}.
\end{equation*}
\end{lemma}

\begin{proof}
The proof was already provided in \cite[Lemma A.2]{DonatiLacaveMiot} and is recalled here for the sake of completeness.

We have that for all $t \in [0,T]$,
\begin{align*}
 |y(t)-z(t)| & \leq \left|\int_0^t \big(y'(s) - z'(s) \big)\dd s \right| + |y(0)-z(0)|
 \\ & \le \int_0^t g(s)\dd s + \left| \int_0^t \big(f(y(s)) - f(z(s))\big) \dd s \right| + |y(0)-z(0)|
 \\ & \le \int_0^t g(s)\dd s + |y(0)-z(0)| + \kappa \int_0^t |y(s)-z(s)|\dd s ,
\end{align*}
so using now the classical Gronwall's inequality, since $t\mapsto \int_0^t g(s)\dd s + |y(0)-z(0)|$ is non negative and differentiable, we have that
\begin{equation*}
 |y(t)-z(t)| \le \left( \int_0^t g(s) \dd s + |y(0)-z(0)|\right)e^{\kappa t}.
\end{equation*}
\end{proof}

\section{On the suitable scaling regime}\label{appendix:regimes}
The purpose of this appendix is to provide a more detailed discussion on the suitable scaling regime leading to the \emph{leapfrogging dynamics} chosen in this paper as well as a detailed comparison to the results of \cite{Mar24}. While the authors obtain a result somewhat in the spirit of Theorem~\ref{theo:main}, the considered asymptotic regime differs significantly for the present one as does the limiting trajectory. 

\subsection{The scaling regime and respective limiting trajectory}
As mentioned in the introduction, a single vortex ring propagates in vertical direction with velocity \eqref{eq:ring_v} equaling 
\[
v_L = |\ln\eps| \frac{\gamma}{4\pi r^\ast}e_z,
\]
hence depending on the reference radius $r^{\ast}$, the intensity $\gamma$ and
in logarithmically diverging manner on its thickness $\epsilon$. When several vortex rings are considered, their motion is driven by the self-induced motion and their mutual interactions which are known to be of strength that is inversely proportional to their mutual distance. Depending on the chosen asymptotic regime, the former \cite{Mar2, Mar3} or the latter \cite{ButtaCavMar24} may be dominant. The leapfrogging regime corresponds to a regime balancing these two effects and the rigorous justification of the asymptotic dynamics is significantly more difficult. Indeed, for the 3D Euler equations it remained completely open until the recent aforementioned \cite{DavilaDelPinoMussoWei, Mar24} and our work. As a rule of thumb and outlined in the introduction, the \emph{leapfrogging dynamics} is expected to occur for vortex rings that travel at leading order at the same self-induced velocity~\eqref{eq:ring_v} of order $|\ln\eps|$, namely share the same reference radius $r^{\ast}$ and intensity $\gamma$. Further, their mutual distance is of order $1/{\sqrt{|\ln \eps|}}$ leading to a mutual interaction of order $\sqrt{|\ln \eps|}$. 
In a time of order $\mathcal{O}\left(1/|\ln \eps|\right)$ (so of order 1 in the rescaled time variable), this, in turn, generates a variation of the rings' radii of order $\mathcal{O}\left(1/{\sqrt{|\ln \eps|}}\right)$ which implies a change in their self-induced velocity~\eqref{eq:ring_v} at the order $\sqrt{|\ln \eps|}$, i.e.
\[
v_L = |\ln\eps| \frac{\gamma}{4\pi r^*}e_z + \cal{O}(\sqrt{|\ln\eps|}).
\]
We emphasize that this regime only occurs for rings of the same circulation $\gamma$. If not, the difference of their self-induced velocity \eqref{eq:ring_v} is of order $|\ln \eps|$, which cannot be balanced by their mutual interaction (of order $\sqrt{|\ln \eps|}$), leading inevitably to the rings moving quickly away from each other. This present scaling regime also corresponds to the one considered in \cite{DavilaDelPinoMussoWei, JerrardSmets}. For sufficiently small $\eps>0$, Theorem~\ref{theo:main} states that the motion of the vortex rings is governed by \eqref{def:X_i} where $t\mapsto X_i^\eps(t)$ denotes the trajectory of the $i$-th vortex ring, i.e.
\begin{equation*}
 \begin{cases}
 \displaystyle \der{}{t} X_i^\eps(t) = \frac{\gamma}{4\pi X_{i,r}^\eps(t)}e_z + \frac{\gamma}{2\pi|\ln \eps|} \sum_{j \neq i} \frac{\big(X_i^\eps(t)-X_j^\eps(t)\big)^\perp}{\big|X_i^\eps(t)-X_j^\eps(t)\big|^2},\\
 X_i^\eps(0) = X_{i,0}^\eps.
 \end{cases}
\end{equation*}
Moreover, defining the variation in the frame moving at constant speed, namely setting $X_i^\eps(t) = X^* + \frac{\gamma}{4\pi r^*} te_z + \frac{1}{\sqrt{|\ln \eps|}} Y_i^\eps(t)$, we derived in Lemma~\ref{lem:X-tildeX} the motion law for $Y_{i}^\eps(t)$ with $Y_i^\eps(0) = Y_{i,0}$:
\begin{align*}
 \der{}{t}Y_i^\eps(t) & = \frac{\gamma}{2\pi} \sum_{j\neq i} \frac{(Y_i^\eps-Y_j^\eps)^\perp}{|Y_i^\eps-Y_j^\eps|^2} + \sqrt{|\ln \eps|}\frac{\gamma}{4\pi r^*} \left( \frac{1}{1+\frac{Y_{i,r}^\eps}{r^*\sqrt{|\ln \eps|}}}- 1 \right) e_z\\
 &= \frac{\gamma}{2\pi} \sum_{j\neq i} \frac{(Y_i^\eps-Y_j^\eps)^\perp}{|Y_i^\eps-Y_j^\eps|^2}-\frac{\gamma}{4\pi (r^{\ast})^2}Y_{i,r}^\eps e_z+ \mathcal{O}\left(\frac{1}{\sqrt{|\ln\eps|}}\right).
\end{align*}
We recognize that the relative self-induced motion and interaction are of the same order. We refer to Section~\ref{sec:limiting_traj} for full details. In particular, we recover in this regime that smallest rings goes faster (see also \eqref{eq:ring_v} for this well-known fact).

Note that as $\eps\rightarrow 0$, see also Lemma~\ref{lem:X-tildeX}, one recovers the following system for the limiting trajectories 
\begin{equation}\label{eq:LP-appB}
 \left\{\begin{aligned}
 & \frac{\dd \tilde{Y}_i}{\dd t} =\frac{\gamma}{2\pi} \sum_{j\neq i} \frac{(\tilde{Y}_i-\tilde{Y}_j)^\perp}{|\tilde{Y}_i-\tilde{Y}_j|^2} - \frac{\gamma}{4\pi (r^*)^2}\tilde{Y}_{i,r} e_z, \\
 & \tilde{Y}_i(0)=Y_{i,0}.
 \end{aligned}\right.
 \end{equation}
 We refer to \cite[Equation (1.14)]{DavilaDelPinoMussoWei} and \cite[Section 1.2] {JerrardSmets} and references therein for a qualitative and quantitative analysis of solutions to \eqref{eq:LP-appB}. Further, we prove in Appendix~\ref{app:nocollision} that no collisions occur in finite time for \eqref{eq:LP-appB}. Appendix~\ref{app:different dynamics} elucidates differences in the behavior of solutions $\tilde{Y}_i$ to \eqref{eq:LP-appB} ($\eps=0$) and $Y_i^\eps$ as defined in \eqref{def:X_i} to \eqref{eq:LP} ($\eps>0$). 
 
\subsection{Comparison of scaling regimes in the literature}\label{sec:compRegimes}
Different from the present the one, see also \cite{DavilaDelPinoMussoWei, JerrardSmets}, the authors in \cite{Mar24} consider a scaling regime in which the self-induced motion of the vortex rings is of the same order as the velocity induced by the interaction of vortex rings. We also refer the reader to \cite[Section 1.1.1-1.1.2]{Meyer} for a comparison of both regimes.
More precisely, the authors in \cite{Mar24} work with vortex rings of radius $\mathcal{O}(|\ln\eps|)$, vortex circulations of order $\mathcal{O}(1)$ and the vorticity $\Omega$ that reads in terms of the potential vorticity $\omega = \frac{\Omega}{r}$, as follows 
\begin{multicols}{2}
\noindent
\begin{equation*}
\begin{cases}
 \partial_t \frac{\Omega}{r} + u \cdot \nabla \left(\frac{\Omega}{r}\right) = 0, \\
 u = \frac{1}{r}\nabla^\perp \Psi,\\
 \div \left( \frac{1}{r} \nabla \Psi \right) = \Omega \\
 \int \Omega_i \dd x = a_i, \\
 \supp \Omega_i \subset B( (0,\alpha|\ln \eps|)+x_i),
\end{cases}
\end{equation*}
\begin{equation*}
 \begin{cases}
 \partial_t \omega + u \cdot \nabla \omega = 0, \\
 u = \frac{1}{r}\nabla^\perp \Psi,\\
 \div \left( \frac{1}{r} \nabla \Psi \right) = r\omega \\
 \int r\omega_i \dd x= a_i, \\
 \supp \omega_i \subset B( (0,\alpha|\ln \eps|)+x_i).
\end{cases}
\end{equation*}
\end{multicols}
For the convenience of the reader, we translate the result of \cite{Mar24} to the scaling regime considered in the present paper. Exploiting the natural scaling invariance of the Euler equations allows one to write this problem in the regime of vortex rings of radius and circulation of order $\mathcal{O}(1)$. To this end, we define
\begin{equation*}
 \tilde{\omega}(t,x) = |\ln \eps|^3\omega(t,|\ln \eps| x), 
\end{equation*}
and check that the vortex rings' circulation is of order $\mathcal{O}(1)$, specifically
\begin{equation*}
 \int r\tilde{\omega}_i(t,x) \dd x = |\ln \eps|^2\int (r\omega_i)(t,|\ln \eps|x) \dd x = a_i.
\end{equation*}
Further, the support satisfies
\begin{equation*}
 \supp \tilde{\omega}_i \subset D\left( (0,\alpha) + \frac{x_i}{|\ln \eps|}\right),
\end{equation*}
so that the scaled rings' radii are of order $\mathcal{O}(1)$. The distances between the rings are now of order $\mathcal{O}(\frac{1}{|\ln \eps|})$, which is more singular than in \eqref{distance-regime} of our analysis. We then compute that
\begin{equation}\label{eq:M1}
 \partial_t \tilde{\omega}(t,x) = |\ln \eps|^3 \partial_t\omega(t,|\ln \eps| x), \quad \nabla \tilde{\omega}(t,x) = |\ln \eps|^4 \nabla\omega(t,|\ln \eps| x).
\end{equation}
Turning to the Biot-Savart law, we seek to determine $\tilde{\Psi}$ such that
\begin{equation*}
 \div \left( \frac{1}{r} \nabla \tilde{\Psi} \right) = r\tilde{\omega}.
\end{equation*}
Setting 
\begin{equation*}
 \tilde{\Psi}(t,x) = \frac{1}{|\ln \eps|}\Psi(t, |\ln \eps|x),
\end{equation*}
one computes that
\begin{align*}
 \div \left( \frac{1}{r} \nabla \tilde{\Psi} \right)(t,x) & = \frac{1}{|\ln \eps|}\left(-\frac{|\ln \eps|e_r}{r^2} \cdot\nabla \Psi(t, |\ln \eps|x) + \frac{|\ln \eps|^2}{r} \Delta \Psi (t, |\ln \eps|x) \right) \\
 & = \frac{|\ln \eps|^3}{|\ln \eps|} \div \left( \frac{1}{r} \nabla \Psi \right)(t,|\ln \eps|x) \\
 & = |\ln \eps|^2 (r \omega)(t,|\ln \eps|x) \\ 
 & = |\ln \eps|^3 r \omega (t,|\ln \eps|x) \\
 & = r \tilde{\omega}(t,x),
\end{align*}
so we found the right candidate for $\tilde{\Psi}$ which leads to
\begin{equation}\label{eq:M3}
 \tilde{u}(t,x) = \frac{1}{r}\nabla^\perp \tilde{\Psi}(t,x) = \frac{1}{r}\frac{|\ln \eps|}{|\ln \eps|} \nabla^\perp\Psi(t, |\ln \eps|x) 
 = |\ln \eps| u(t, |\ln \eps|x).
\end{equation}

Gathering relations~\eqref{eq:M1}
and~\eqref{eq:M3} we get from the equation for the potential vorticity $\omega$ that
\begin{equation*}
 \frac{1}{|\ln \eps|^3} \partial_t \tilde{\omega} + \frac{1}{|\ln \eps|^5}\tilde{u} \cdot \nabla \tilde{\omega} = 0
\end{equation*}
namely
\begin{equation*}
 \partial_t \tilde{\omega} + \frac{1}{|\ln \eps|^2}\tilde{u} \cdot \nabla \tilde{\omega} = 0.
\end{equation*}
In conclusion, the results of \cite{Mar24} can be interpreted as a strong confinement result (in both directions) for vortex rings of radius of order $\mathcal{O}(1)$ and circulations of order $\mathcal{O}(1)$, placed at a distance of order $\mathcal{O}\left(\frac{1}{|\ln \eps|}\right)$, for a time of order $\mathcal{O}\left(\frac{1}{|\ln \eps|^2}\right)$. Indeed, one can perform a final change of the time variable
\begin{equation*}
 \hat{\omega}(t,x) = \tilde{\omega}(t |\ln \eps|^2,x) \, \qquad \hat{u}(t,x) = \tilde{u}(t |\ln \eps|^2,x) \, ,
\end{equation*}
to recover the transport equation
\begin{equation*}
 \partial_t \hat{\omega} + \hat{u} \cdot \nabla \hat{\omega} = 0.
\end{equation*}
This regime corresponds to the second regime in \cite[Section 1.1.2]{Meyer} where it is expected that the two rings have the same local velocities. The velocities $v_L$ (local induction approximation) and $F$ (generated by the other rings) are both of order $|\ln \eps|$. This means that in this regime the vortex rings perform a rotation around each other over a distance of order $\mathcal{O}\left(\frac{1}{|\ln\eps|}\right)$ in a time of order $1/|\ln\eps|^2$. Moreover, in that regime, the vortex rings can have different circulations $\gamma_i$, and in that case the limiting trajectories were proven \cite[Theorem 2.1]{Mar24} to satisfy
\begin{equation}\label{eq:LPMAR-appB}
 \left\{\begin{aligned}
 & \frac{\dd \tilde{Y}_i}{\dd t} =\frac{\gamma_i}{2\pi} \sum_{j\neq i} \frac{(\tilde{Y}_i-\tilde{Y}_j)^\perp}{|\tilde{Y}_i-\tilde{Y}_j|^2} + \frac{\gamma_i}{4\pi r^\ast} e_z \\
 & \tilde{Y}_i(0)=Y_{i,0}.
 \end{aligned}\right.
 \end{equation}
Note that \eqref{eq:LPMAR-appB} differs from \eqref{eq:LP-appB} in terms of the vertical translation, and we remark that two vortex rings with the same circulation have the same local induction velocity, no matter of their radius. The properties of solutions to \eqref{eq:LPMAR-appB} have been investigated in \cite{MarNegrini}, in particular ruling out blow-up provided that all $\gamma_i$ are of the same sign \cite[Theorem 2.1]{MarNegrini}. See also \cite[Section 7]{Mar24} for the case of two vortex rings.

For the convenience of the reader we provide a comparison of the choice of parameters between our regime and the one studied in \cite{Mar24}.
\begin{center}
\small{
\begin{tabular}{c| c| c| c| c| c| c|}
 & radius $r^{\ast}$ & intensity $\gamma$ & velocity $v_L$ & mutual distance & interaction & time-scale \\ 
 \hline
our scaling & $\mathcal{O}(1)$ & $\mathcal{O}(1)$ & $ \frac{|\ln\eps|\gamma}{4\pi r^*}e_z + \cal{O}(\sqrt{|\ln\eps|})$ & $\mathcal{O}\left(\frac{1}{\sqrt{|\ln\eps|}}\right)$ & $\mathcal{O}(\sqrt{|\ln\eps|})$ & $\mathcal{O}\left(\frac{1}{|\ln\eps|}\right)$\\ 
\hline
re-scaling of \cite{Mar24}& $\mathcal{O}(1)$ & $\mathcal{O}(1)$ & $\mathcal{O}(|\ln\eps|)$ & $\mathcal{O}\left(\frac{1}{{|\ln\eps|}}\right)$ & $\mathcal{O}({|\ln\eps|})$ & $\mathcal{O}\left(\frac{1}{|\ln\eps|^2}\right)$ \\ 
\end{tabular}
}
\end{center}
Because of their mutual distances and the intensity of their velocities, it should be noted that both time scales are respectively suitable for observing a rotation of one vortex ring with respect to the other.

\section{Some comments on the limit dynamics}\label{app:limit dyn}
As detailed in the introduction, we have considered two limiting trajectories: $(Y_i^\eps)_{i=1,\dots,N}$ defined in \eqref{eq:LP-bis} and $(\tilde{Y}^i)_{i=1,\dots,N}$ obeying \eqref{eq:LP}, and we have proved in Lemma~\ref{lem:X-tildeX} that they remain close, namely at distance $\cal O(1/\sqrt{|\ln\eps|})$. We first remark that there is no-collision in finite time for $(\tilde{Y}^i)$ and that they are global in time. Next, we present some numerical evidence on the difference of these two limiting dynamics. Let us mention that the dynamics $(Y_i^\eps)$ also differs from the dynamics studied in \cite{JerrardSmets}, even though the two dynamics have the same limiting trajectories $(\tilde{Y}_i(t))$.

\subsection{No collision in finite time and global solutions for the second limiting system}\label{app:nocollision}

For simplicity, we rewrite System~\eqref{eq:LP} here as
\begin{equation}\label{eq:LP-annexe}
 \left\{\begin{aligned}
 & \frac{\dd q_i}{\dd t} =\frac{\gamma}{2\pi} \sum_{j\neq i} \frac{(q_i-q_j)^\perp}{|q_i-q_j|^2} - \frac{\gamma}{4\pi (r^*)^2}{q}_{i,r} e_z \\
 & q_i(0)=q_{i,0}.
 \end{aligned}\right.
 \end{equation}
Although no collisions in finite time is probably well-known, we provide the proof here since we did not find it in the literature. By Cauchy-Lipschitz theorem, there exists a unique local $C^1$ solution of \eqref{eq:LP-annexe}.

As observed in \cite{DavilaDelPinoMussoWei}, System~\eqref{eq:LP-annexe} is a Hamiltonian system for the energy
\begin{equation*}
H_0:=\frac{1}{2}\sum_{i\neq j}\ln |q_i(t)-q_j(t)|+ \frac{1}{4(r^\ast)^2}\sum_{i=1}^N q_{i,r}(t)^2.
\end{equation*}
We can indeed check that this quantity is conserved.
\begin{lemma}
 Let $t\mapsto (q_1,\dots,q_N)(t)$ a $C^1$ solution on $[0,T]$ of \eqref{eq:LP-annexe}, then
 \[
 \frac{\dd H_0}{\dd t}=0.
 \]
\end{lemma}
\begin{proof}
 We differentiate $H_0$ to obtain
 \begin{align}
 \frac{\dd H_0}{\dd t} =&
 \frac{1}{2}\sum_{i\neq j}\frac{q_i(t)-q_j(t)}{|q_i(t)-q_j(t)|^2} \cdot (\dot q_i(t)-\dot q_j(t))+ \frac{1}{2(r^\ast)^2}\sum_{i=1}^N q_{i,r}(t)\dot q_{i,r}(t) \nonumber\\
 =&
 \sum_{i\neq j}\frac{q_i(t)-q_j(t)}{|q_i(t)-q_j(t)|^2} \cdot \dot q_i(t)+ \frac{1}{2(r^\ast)^2}\sum_{i=1}^N q_{i,r}(t)\dot q_{i,r}(t). \label{eq:hamil}
 \end{align}
By using \eqref{eq:LP-annexe}, we compute the first right-hand-side sum
\begin{align*}
 \sum_{i\neq j}\frac{q_i-q_j}{|q_i-q_j|^2} \cdot \dot q_i
 =&\frac{\gamma}{2\pi} \sum_{i=1}^N\sum_{j\neq i}\sum_{k\neq i}\frac{q_i-q_j}{|q_i-q_j|^2} \cdot 
 \frac{(q_i-q_k)^\perp}{|q_i-q_k|^2} 
 - \frac{\gamma}{4\pi (r^*)^2} \sum_{i=1}^N\sum_{j\neq i}\frac{q_i-q_j}{|q_i-q_j|^2} \cdot {q}_{i,r} e_z \\
 =& \frac{-\gamma}{4\pi (r^*)^2} \sum_{i=1}^N\sum_{j\neq i} q_{i,r} \frac{q_{i,z}-q_{j,z}}{|q_i-q_j|^2} 
\end{align*}
by inverting $k$ and $j$ in the first term.

The last sum in \eqref{eq:hamil} is computed in the same way:
\begin{align*}
 \frac{1}{2(r^\ast)^2}\sum_{i=1}^N q_{i,r}\dot q_{i,r}
 =& \frac{1}{2(r^\ast)^2}\frac{\gamma}{2\pi} \sum_{i=1}^N \sum_{j\neq i} q_{i,r} e_r \cdot 
 \frac{(q_i-q_j)^\perp}{|q_i-q_j|^2} 
 = \frac{\gamma}{4\pi(r^\ast)^2} \sum_{i=1}^N \sum_{j\neq i} q_{i,r} 
 \frac{q_{i,z}-q_{j,z}}{|q_i-q_j|^2}
\end{align*}
where we recall our convention $e_z=(1,0)$, $e_r=(0,1)$ and $a^\perp =(-a_r,a_z)$.
This completes the proof of the conservation of the energy $H_0$.
\end{proof}

Moreover, we observe that the angular momentum
\begin{equation*}
 I(t):=\frac{1}{2}\sum_{i=1}^N |q_i(t)|^2
\end{equation*}
satisfies
\begin{align*}
 \frac{\dd}{\dd t}I(t)&=\frac{\gamma}{2\pi}\sum_{i\neq j} q_i(t)\cdot \frac{(q_i(t)-q_j(t))^\perp}{|q_i(t)-q_j(t)|^2}-\frac{\gamma}{4\pi (r^\ast)^2}\sum_{i=1}^N q_{i,r}(t)q_{i,z}(t)\\
 &=-\frac{\gamma}{2\pi}\sum_{i\neq j} \frac{q_i(t)\cdot q_j^\perp(t)}{|q_i(t)-q_j(t)|^2}-\frac{\gamma}{4\pi (r^\ast)^2}\sum_{i=1}^N q_{i,r}(t)q_{i,z}(t)\\
 &=-\frac{\gamma}{4\pi (r^\ast)^2}\sum_{i=1}^N q_{i,r}(t)q_{i,z}(t)\leq \frac{\gamma}{8\pi (r^\ast)^2}I(t),
 \end{align*}
where we have exchanged $i$ and $j$ in the first sum, hence
\[
I(t)\leq I(0)e^{\gamma t/(8\pi r^{\ast2})}.
\]

As a consequence, using $\ln |x-y|\leq |x-y| \leq |x|+|y|$, we get 
\begin{align*}
 \sum_{i\neq j||q_i-q_j|<1}|\ln |q_i(t)-q_j(t)||&= \sum_{i\neq j||q_i-q_j|>1}\ln |q_i(t)-q_j(t)| -2H_0 +\frac{1}{2 (r^\ast)^2}\sum_{i=1}^Nq_{i,r}(t)^2\\
 &\leq2N\sum_{i=1}^N |q_i(t)| -2H_0 +\frac{1}{2 (r^\ast)^2} I(t)\\
 &\leq 2N\sqrt{N}\sqrt{I(t)}-2H_0 +\frac{1}{2 (r^\ast)^2} I(t),
\end{align*}
for all $t\in [0,T^\ast)$, where $T^\ast$ is the maximum time of existence of \eqref{eq:LP-annexe}. The conclusion of non-collision in finite time then follows from the bound on $I(t)$. By the estimate on $I$, we also know that the particles can not escape to infinity in finite time, which implies that the unique $C^1$-solution of \eqref{eq:LP-annexe} is global.

\subsection{Difference between the two dynamics}\label{app:different dynamics}

We provide numerical evidence in the case $N=2$ that the dynamics of $Y_i^\eps(t)$ as defined in \eqref{eq:LP-bis} differs from the limit dynamics for $\tilde{Y}_i$ obeying \eqref{eq:LP}. For convenience, we choose $X^* = (X_1^\eps(0) + X_2^\eps(0))/2$ so that $Y_1(0) + Y_2(0) = 0$. Let us recall that we proved that as $\eps \to 0$, then $Y_i^\eps(t) \to \tilde{Y}_i(t)$ uniformly in $t \in [0,T_0]$, see Lemma~\ref{lem:X-tildeX}. However, for a given $\eps > 0$, the respective dynamics of the points $Y_i^\eps(t)$ and $\tilde{Y}_i(t)$ have different properties. 

First, we observe that for the second system, the center of mass $(\tilde{Y}_1 + \tilde{Y}_2)/2$ does not move:
\begin{equation*}
 \der{}{t} \frac{\tilde{Y}_1 + \tilde{Y}_2}{2} = - \frac{\gamma}{4\pi (r^*)^2}\frac{\tilde{Y}_{1,r}+\tilde{Y}_{2,r}}2 e_z,
\end{equation*}
which implies that $\tilde{Y}_{1,r} + \tilde{Y}_{2,r}$ remains constant. As we have assumed here that this is initially zero, the previous equality also gives that the vertical component is at rest, so in the end, the $(\tilde{Y}_1+\tilde{Y}_2)/2$ is constant. Then, by showing that every level set of the associated Hamiltonian is bounded, it is possible to deduce from this that the trajectories of the system $(\tilde{Y}_1,\tilde{Y}_2)$ are periodic in time for every initial data, see \cite{DavilaDelPinoMussoWei}.

\begin{figure}[ht]
 \centering
 \subfloat[Trajectories of $Y_1^\eps(t)$ and $Y_2^\eps(t)$]{
 \fbox{
 \begin{minipage}[c][4cm][c]{7.7cm}
 \centering
 \begin{tikzpicture}
 \node[inner sep=0pt] (img) at (-1,0) {\includegraphics[height=3.2cm]{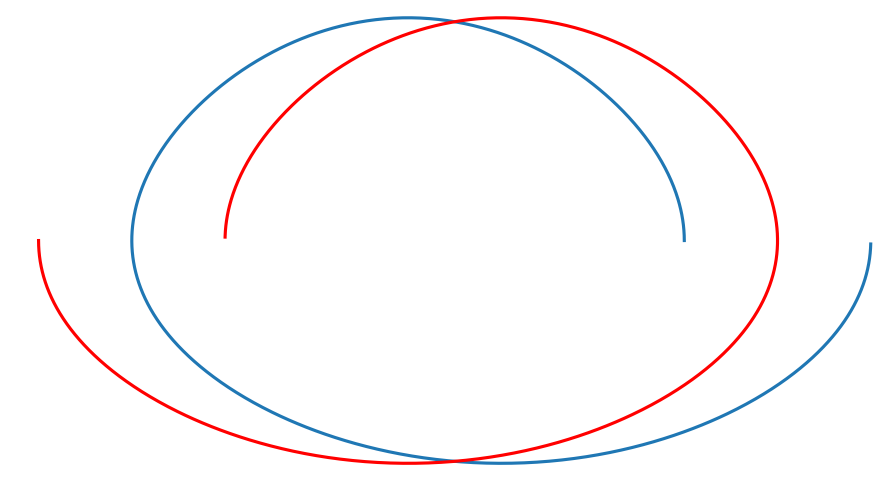}};
 \draw[->, thick] (-4,-1.8) -- (-3.5,-1.8) node[below right] {$z$};
 \draw[->, thick] (-4,-1.8) -- (-4,-1.3) node[above left] {$r$};
 \draw (0,1.8) node {};
 \draw[->,thick, red](-3.5,-0.65)--(-3.5+0.01,-0.65-0.01414);
 \draw[->,thick, blue](0.35,0.65)--(0.35-0.01,0.65+0.01414);
 \end{tikzpicture}
 \end{minipage}
 }
}
 \hfill
 \subfloat[Trajectory of $Y_1^\eps(t)$]{
 \fbox{
 \begin{minipage}[c][4cm][c]{7.7cm}
 \centering
 \includegraphics[height=3.2cm]{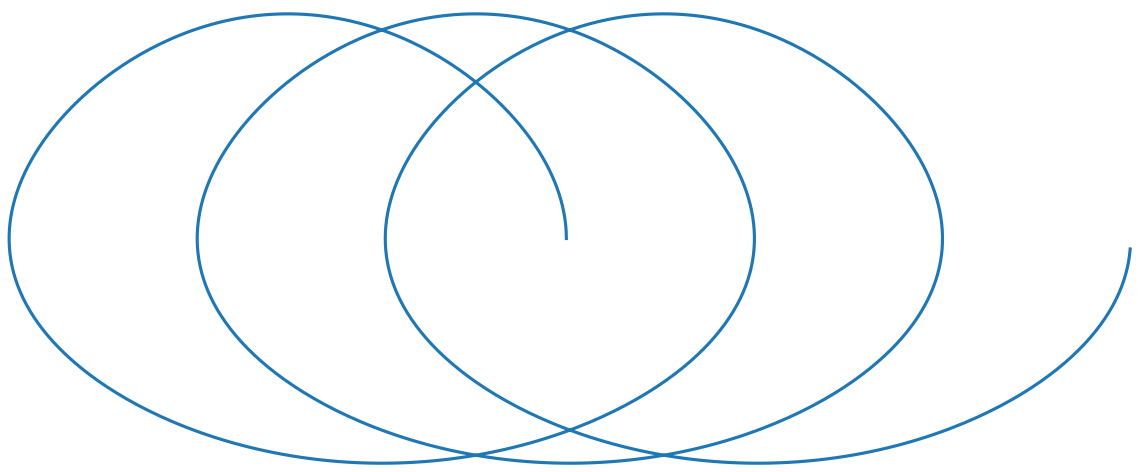}
 \end{minipage}
 }
 }
 \par\vspace{0.2em}
 \subfloat[Time evolution of $(Y_{1,z}^\eps(t)+Y_{2,z}^\eps(t))/2$]{
 \fbox{
 \begin{minipage}[c][5cm][c]{7.7cm}
 \centering
 \includegraphics[height=4.5cm]{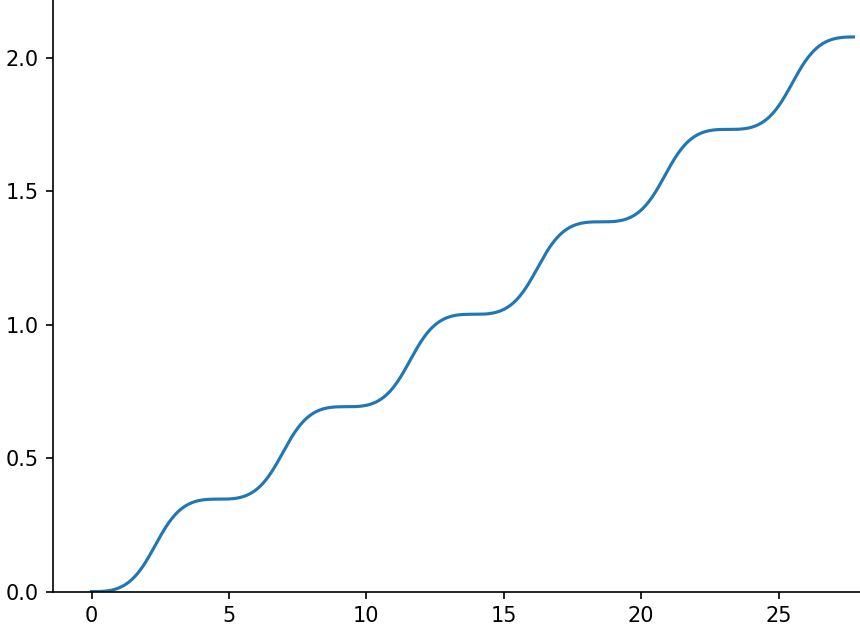}
 \end{minipage}
 }
 }
 \hfill
 \subfloat[Vertical velocity of $(Y_1^\eps(t)+Y_2^\eps(t))/2$]{
 \fbox{
 \begin{minipage}[c][5cm][c]{7.7cm}
 \centering
 \includegraphics[height=4.5cm]{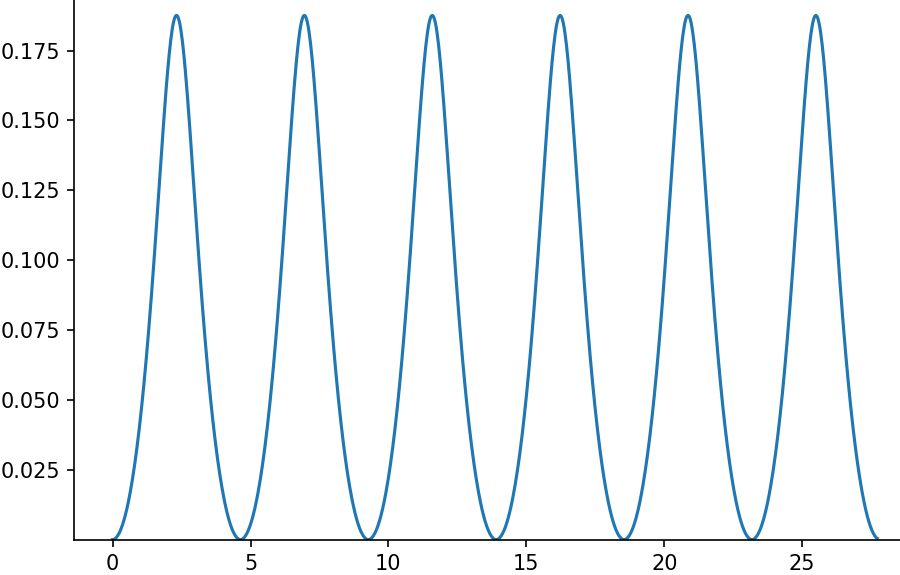}
 \end{minipage}
 }
 }
 \caption{Numerical computation of the motion in $\R^2$ of $Y_1^\eps(t)$ and $Y_2^\eps$ during a full rotation (A), and the same motion of $Y_1^\eps(t)$ during three full rotations (B). The shift towards positive $z$ is seen in (C), with a non-trivial velocity (D). The parameters used are $\eps = 0.01$, $\gamma = 2\pi$, $r^*=1$ and $Y_1(0) = (0,1.2)$.}
 \label{fig:trajs}
\end{figure}

This property on the center of mass is false for the system $(Y_1^\eps,Y_2^\eps)$ as can be seen by the following computation. By taking the derivative in time of the center of mass, we get that
\begin{equation}\label{eq:der_centre_masse}
 \der{}{t} \frac{Y_1^\eps+Y_2^\eps}{2} = \sqrt{|\ln \eps|}\frac{\gamma}{8\pi r^*} \left( \frac{1}{1+\frac{Y_{1,r}^\eps}{r^*\sqrt{|\ln \eps|}}}+\frac{1}{1+\frac{Y_{2,r}^\eps}{r^*\sqrt{|\ln \eps|}}}- 2 \right) e_z
\end{equation}
which has again no radial part, but cannot vanish at all times. Moreover, by noticing that the radial components are bounded and opposite at every time, we perform a Taylor expansion of \eqref{eq:der_centre_masse} up to the second order to get that
\begin{align*}
 \der{}{t} \frac{Y_1^\eps+Y_2^\eps}{2} & = \sqrt{|\ln \eps|}\frac{\gamma}{8\pi r^*} \left( -\frac{Y_{1,r}^\eps + Y_{2,r}^\eps}{r^*\sqrt{|\ln \eps|}} + \frac{(Y_{1,r}^\eps)^2 + (Y_{2,r}^\eps)^2}{(r^*)^2|\ln \eps|} + \mathcal{O}\left(\frac{1}{|\ln \eps|^{3/2}}\right) \right)e_z \\
 & = \frac{\gamma}{8\pi}\frac{(Y_{1,r}^\eps)^2 + (Y_{2,r}^\eps)^2}{(r^*)^3 \sqrt{|\ln \eps|}}e_z+ \mathcal{O}\left(\frac{1}{|\ln \eps|}\right).
\end{align*}
This quantity being non-negative, we expect the dynamics of $(Y_1^\eps,Y_2^\eps)$ to display a non-linear shift towards positive $z$. This fact is observed numerically as shown in Figure~\ref{fig:trajs}.

\begin{figure}[ht]
 \centering
 \includegraphics[width=0.8\linewidth]{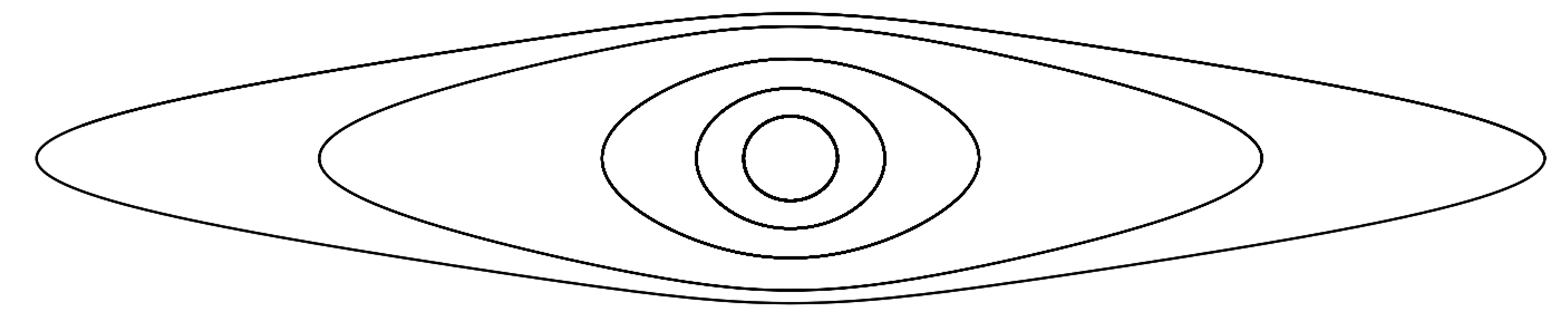}
 \caption{Phase portrait of the centered dynamics $(\hat{Y}_1^\eps,\hat{Y}_2^\eps)$.}
 \label{fig:portrait_phase}
\end{figure}

 To understand more in detail the dynamics, one could introduce the corrected positions made such that the center of mass of the configuration is preserved, namely the points
\begin{equation*}
 \hat{Y}_i^\eps(t) = Y_i^\eps(t) - \frac{Y_1^\eps(t)+Y_2^\eps(t)}{2},
\end{equation*}
which satisfy by construction that $\hat{Y}_1^\eps = - \hat{Y}_2^\eps$, meaning that the dynamics reduces to the study of a single point in $\R^2$, allowing to obtain the phase portrait type Figure~\ref{fig:portrait_phase}, displaying different trajectories of $\hat{Y}_1^\eps(t)$ (and thus of $\hat{Y}_2^\eps(t)$) for various initial data. In particular, we recover a periodic motion, with a phase portrait very similar to that of the $\tilde{Y}_i$, see \cite{DavilaDelPinoMussoWei}. In conclusion, we claim that the dynamics of $(Y_1^\eps,Y_2^\eps)$ differs from the dynamics of $(\tilde{Y}_1,\tilde{Y}_2)$ essentially by a non-trivial shift towards the positive $z$.

%\bibliographystyle{abbrv}
%\bibliography{reference} 

\end{document}